# ON FOLIATIONS IN $\mathrm{PSL}_4(\mathbb{R})$-TEICHMÜLLER THEORY

ALEXANDER NOLTE

ABSTRACT. We carry out a detailed study of the structure of domains of discontinuity $\Omega_\rho$ in $\mathbb{RP}^3$ of $\mathrm{PSL}_4(\mathbb{R})$-Hitchin representations $\rho$. We then prove the foliated component $\Omega_\rho^1$ of $\Omega_\rho$ has exactly two group-invariant foliations by properly embedded projective line segments and has a unique foliation by properly embedded convex domains in projective planes. This gives a finiteness counterpart to work of Guichard and Wienhard. We also prove analogues for the non-foliated component $\Omega_\rho^2$ and deduce a rigidity of projective equivalences of properly convex foliated projective structures on unit tangent bundles of surfaces.

Domains of discontinuity for $\mathrm{PSL}_n(\mathbb{R})$-Hitchin representations in flag manifolds $X$ have rich structure. We systematically study the case $n = 4$ and $X = \mathbb{RP}^3$. In contrast to the much-studied case $n = 3$ and $X = \mathbb{RP}^2$, our domains are neither convex nor $C^1$ (Figure 1).

The maximal open domain $\Omega_\rho = \Omega_\rho^1 \sqcup \Omega_\rho^2$ of discontinuity in $\mathbb{RP}^3$ for a $\mathrm{PSL}_4(\mathbb{R})$ Hitchin representation $\rho$ of the fundamental group $\Gamma$ of a closed, oriented hyperbolic surface $S$ was first studied in celebrated work of Guichard-Wienhard [11]. A central role was played in [11] by a curious pair $(\mathcal{F}_{\mathrm{pcf}}, \mathcal{G}_{\mathrm{pcf}})$ of nested $\rho(\Gamma)$-invariant foliations of the *foliated component* $\Omega_\rho^1$, the leaves of which are properly convex domains in projective planes and lines, respectively.

We classify all geometrically similar foliations here. We construct a new invariant foliation $\mathcal{G}_{\mathrm{tcf}}$ of $\Omega_\rho^1$ by projective line segments (Figure 2 and §3). Our central rigidity theorem is:

**Theorem A** (Exactly Two). *Let $\rho : \Gamma \to \mathrm{PSL}_4(\mathbb{R})$ be Hitchin. Then $\mathcal{G}_{\mathrm{pcf}}$ and $\mathcal{G}_{\mathrm{tcf}}$ are the only $\rho(\Gamma)$-invariant foliations of $\Omega_\rho^1$ by properly embedded projective line segments.*

*$\mathcal{F}_{\mathrm{pcf}}$ is the unique foliation of $\Omega_\rho^1$ by properly embedded convex domains in projective planes.*

Thm. A is a finiteness counterpart to Guichard-Wienhard's work. Indeed, in [11] the structure induced on $\Omega_\rho^1/\rho(\Gamma)$ by the foliations $\mathcal{F}_{\mathrm{pcf}}$ and $\mathcal{G}_{\mathrm{pcf}}$ is the basis for a theory of projective structures on the unit tangent bundle $\mathrm{T}^1 S$ with Hitchin holonomy. We classify pairs of foliations of $\Omega_\rho^1/\rho(\Gamma)$ that can be the basis of similar theories: there are exactly two.

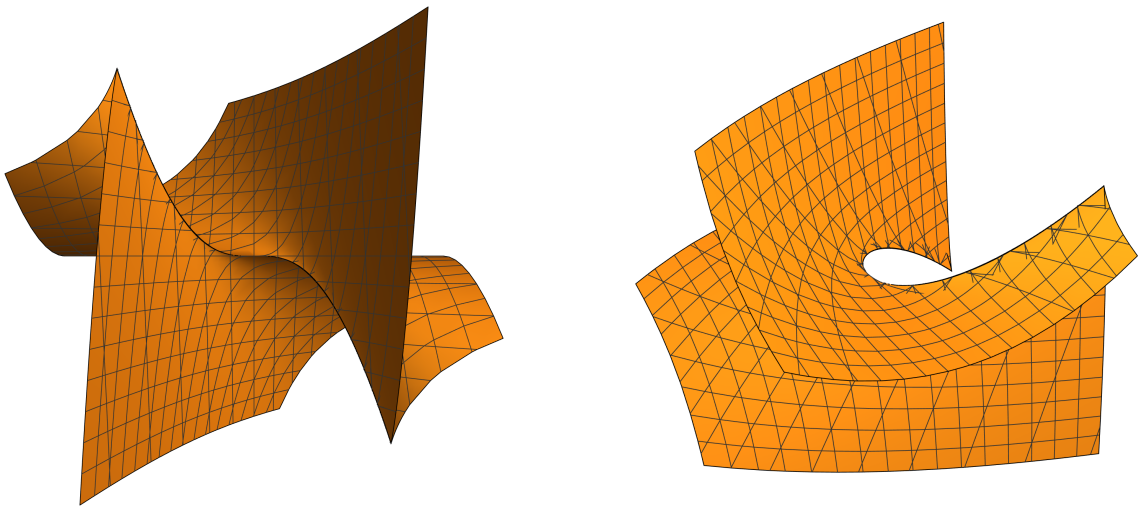

FIGURE 1. The domain of discontinuity for a Fuchsian $\mathrm{PSL}_4(\mathbb{R})$ Hitchin representation in an affine chart for $\mathbb{RP}^3$. The analogue for $n = 3$ is an ellipse.

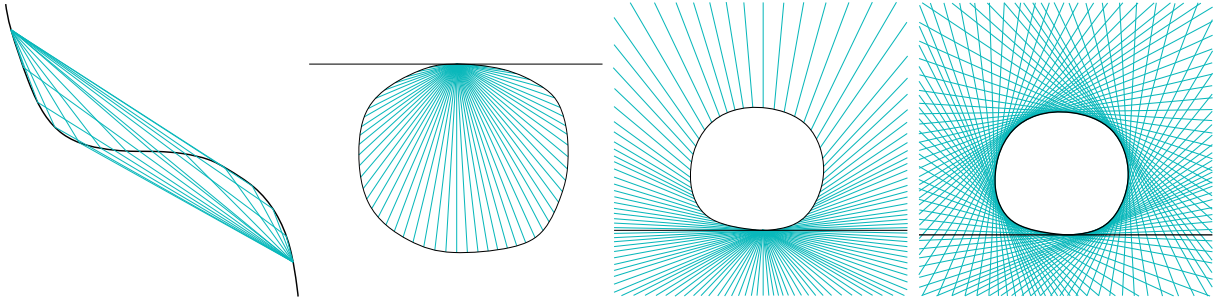

FIGURE 2. Left to right: $\mathcal{G}_{\text{tcf}}$, $\mathcal{G}_{\text{pcf}}$, $\mathcal{G}_{\text{ctrf}}$, and $\mathcal{G}_{\text{ctaf}}$. The abbreviations are for *transverse convex foliation*, *properly convex foliation*, *concave transverse trifoliation* and *concave tangent trifoliation*. Constructions are in §3.

A central question in higher Teichmüller theory [12, 13, 29] is to give a *qualitative* interpretation of $\text{PSL}_n(\mathbb{R})$-Hitchin representations in terms of analogues of hyperbolic surfaces in Teichmüller theory. Guichard-Wienhard's landmark work in $\text{PSL}_4(\mathbb{R})$ is the highest-rank known description of this type. Thm. A suggests that there should be only finitely many ways to carry out similar descriptions for any fixed $n$ and flag manifold in general. We only discuss the case $n = 4$ here. No analogues of Thm. A are known in similar settings.

Thm. A has sharp hypotheses. We do not assume $\rho(\Gamma)$-invariance or proper convexity in our result on $\mathcal{F}_{\text{pcf}}$. We allow once-punctured projective lines in our result on $\mathcal{G}_{\text{pcf}}$ and $\mathcal{G}_{\text{tcf}}$; the result is false if projective lines are allowed as leaves or if $\rho(\Gamma)$-invariance is not assumed (§6). We also prove an analogue of Thm. A for the other component $\Omega_\rho^2$ of $\Omega_\rho$ (Thm. F).

The source of rigidity in Thm. A is rather subtle. Our proof relies on a detailed understanding of the structure of $\partial \Omega_\rho$ and the precise nature of its singularities. What we prove goes well beyond what is known in general about similar representations. A sample statement is below. For notation, let $\partial \Gamma$ be the Gromov boundary of $\Gamma$, which is a circle with a natural $\Gamma$-action, let $\mathcal{F}(\mathbb{R}^4)$ be the manifold of full flags $(V_1, V_2, V_3)$ of nested subspaces of $\mathbb{R}^4$, and let $\xi_\rho$ be the unique continuous $\rho$-equivariant map $\partial \Gamma \to \mathcal{F}(\mathbb{R}^4)$ (see §2, or [10, 16]).

**Theorem B** (Intersection Classification). *Let $\rho \in \text{Hit}_4(S)$ have limit map $\xi_\rho = (\xi^1, \xi^2, \xi^3)$. Then every projective plane $P$ in $\mathbb{RP}^3$ has one of the four mutually exclusive forms $\xi^3(x)$, or $\xi^2(x) \oplus \xi^1(y)$, or $\xi^1(x) \oplus \xi^1(y) \oplus \xi^1(z)$, or $(\xi^3(x) \cap \xi^3(z)) \oplus \xi^1(y)$ for some $x, y, z \in \partial \Gamma$.*

*In all cases, $P \cap \partial \Omega_\rho$ consists of at most three properly convex $C^1$ arcs and at most one projective line. Figure 3 accurately depicts the configurations of arcs and lines in each case.*

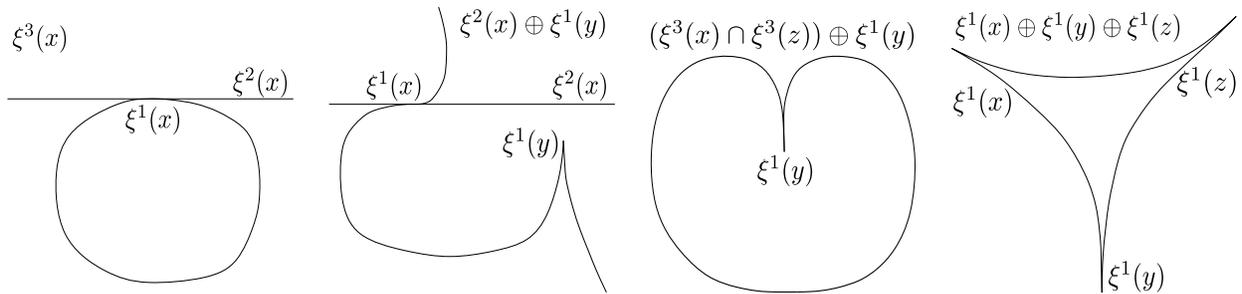

FIGURE 3. All types of intersections of $\partial \Omega_\rho$ with projective planes in $\mathbb{RP}^3$.

Thm. B is an instance of the following theme that we develop: though $\Omega_\rho^1$ is neither convex nor $C^1$, convexity and $C^1$ regularity are ubiquitous in objects associated to $\Omega_\rho^1$, and the failures of $\Omega_\rho^1$ to be convex and $C^1$ may be precisely understood.

The simplest case of intersections in Thm. B, of $\xi^3(x)$, is in Guichard-Wienhard's work [11]; the other cases are new. Other qualitative features of the shape of $\Omega_\rho$ are established in §4. Notably, the structure of a "front view" of $\partial\Omega_\rho$ that generalizes the shape of Figure 1, Left is a key point in the proof of Thm. A.

These results go well beyond the state of knowledge on domains of discontinuity for Anosov representations, and what was known about $\Omega_\rho$. For instance, even the following is new:

**Proposition C** (Locally Convex off Frenet). *A point $p \in \partial\Omega_\rho$ is a $C^1$ point if and only if $p$ is not in the image of $\xi^1$. For any point $p \notin \mathrm{Im}(\xi^1)$, there is a neighborhood $U \subset \partial\Omega_\xi$ of $p$ that is contained in the boundary of a properly convex domain in $\mathbb{RP}^3$.*

1.0.1. *Properly Convex Foliated Projective Structures.* A consequence of Thm. A is a rigidity of projective equivalences of properly convex foliated projective structures. These are refinements of Thurston-Klein $(\mathrm{PSL}_4(\mathbb{R}), \mathbb{RP}^3)$-structures on the unit tangent bundle $\mathrm{T}^1 S$, defined below and in [11], which carry the extra data of distinguished foliations. In [11], Guichard-Wienhard proved that a holonomy map from a moduli space $\mathcal{P}_{\mathrm{pcf}}(\mathrm{T}^1 S)$ of such structures was a homeomorphism onto the $\mathrm{PSL}_4(\mathbb{R})$-Hitchin component $\mathrm{Hit}_4(S)$.

In detail, a properly convex foliated projective structure $(A, \mathcal{F}, \mathcal{G})$ is a $(\mathrm{PSL}_4(\mathbb{R}), \mathbb{RP}^3)$-structure $A$ on $\mathrm{T}^1 S$ and two nested foliations $\mathcal{F}$ and $\mathcal{G}$ that are isotopic as a pair to the stable foliation and geodesic foliation[1] of $\mathrm{T}^1 S$, respectively. Lifts of leaves of $\mathcal{F}$ to the universal cover of $\mathrm{T}^1 S$ are required to develop to properly convex domains in copies of $\mathbb{RP}^2$ and lifts of leaves of $\mathcal{G}$ are required to develop to properly convex domains in copies of $\mathbb{RP}^1$.

The relevant point here is that Guichard-Wienhard's moduli space $\mathcal{P}_{\mathrm{pcf}}(\mathrm{T}^1 S)$ of properly convex foliated projective structures is defined in terms of a *refined* equivalence relation. Specifically, two properly convex foliated projective structures $(A, \mathcal{F}_1, \mathcal{G}_1)$ and $(B, \mathcal{F}_2, \mathcal{G}_2)$ are equivalent if there is a projective equivalence $\varphi$ of $A$ and $B$ isotopic to the identity that satisfies the further requirement of respecting foliations, i.e. $\varphi_* \mathcal{F}_1 = \mathcal{F}_2$ and $\varphi_* \mathcal{G}_1 = \mathcal{G}_2$.

It is unclear from first principles how much this refined equivalence relation differs from marked projective equivalence. A corollary of our main results are that the two are the same:

**Corollary D** (Equivalence Rigidity). *Let $(A, \mathcal{F}, \mathcal{G})$ and $(B, \mathcal{F}', \mathcal{G}')$ be properly convex foliated projective structures on $\mathrm{T}^1 S$. Then any identity-isotopic projective equivalence $\varphi: A \to B$ is an equivalence of properly convex foliated projective structures: $\varphi^*(\mathcal{F}, \mathcal{G}) = (\mathcal{F}', \mathcal{G}')$.*

Cor. D does not use most of the strength of Thm. A, because the nesting of foliations in the definition of properly convex foliated projective structures is a major constraint that is not available in the generality of the setting of Thm. A.

A matter related to this refinement of equivalence relation was if the natural map from $\mathcal{P}_{\mathrm{pcf}}(\mathrm{T}^1 S)$ to the moduli space $\mathcal{P}(\mathrm{T}^1 S)$ of Thurston-Klein $(\mathrm{PSL}_4(\mathbb{R}), \mathbb{RP}^3)$-structures on $\mathrm{T}^1 S$ was injective (it is a case of [11, Rmk. 2.4]). Comparing the results of [11] and Guichard-Wienhard's later work [12, Theorem 11.3] allows for one to conclude that the natural map $\mathcal{P}_{\mathrm{pcf}}(\mathrm{T}^1 S) \to \mathcal{P}(\mathrm{T}^1 S)$ is in fact injective (see e.g. [25, §3] for a similar deduction). Cor. D also implies, much more directly, injectivity of the natural map $\mathcal{P}_{\mathrm{pcf}}(\mathrm{T}^1 S) \to \mathcal{P}(\mathrm{T}^1 S)$.

---

[1]The stable foliation and the geodesic foliation of $\mathrm{T}^1 S$ are specified by a choice of reference hyperbolic metric, but are well-known to be independent of the choice as topological foliations (e.g. [11] or §2 below).

**Remark.** *The proofs of Thms. A above and F below rely on the construction of two families of projective structures on $\mathrm{T}^1 S$ of a similar flavor to those of* [11]. *These projective structures have explicit developing maps and map canonical foliations of the space $\partial \Gamma^{(3)+}$ of positively oriented triples in $\partial \Gamma$ to distinguished domains in projective lines and planes.*

*It seems likely that analogues of Guichard-Wienhard's main theorem in* [11] *hold for moduli spaces of "decorated" projective structures that formalize the features of these examples. This would be interesting to prove. To maintain a focus on the geometry of representations that are already known to be Hitchin in this paper, we do not pursue this direction here.*

1.1. **Remarks and Further Results.** First, we remark on some conventions. Throughout the paper, all foliations are topological and are considered equivalent if and only if all of their leaves agree. We take the phrase *projective line segment* to mean any open convex subset of a projective line other than the entire projective line.

1.1.1. *Universality.* Next, we note that our developing maps, constructions of foliations, and analysis of the shape of $\partial \Omega_\rho$ do not use the group invariance of limit maps and hold in a "universal" sense. For instance, for an arbitrary hyperconvex Frenet curve $\xi$ from the circle $\partial \Gamma$ to $\mathcal{F}(\mathbb{R}^4)$ (§2), there is an analogous region $\Omega_\xi$ to $\Omega_\rho$ in the group-invariant case, and a connected component $\Omega_\xi^1$ of $\Omega_\xi$ is foliated by convex domains. Call this foliation $\mathcal{F}_{\mathrm{pcf}}$. The codimension-1 uniqueness in Thm. A holds in this level of generality:

**Theorem E** (Universal Uniqueness). *Let $\xi$ be a hyperconvex Frenet curve in $\mathcal{F}(\mathbb{R}^4)$. Then $\mathcal{F}_{\mathrm{pcf}}$ is the only foliation of $\Omega_\xi^1$ by properly embedded convex domains in projective planes.*

We emphasize that the lack of a group action hypothesis is a novelty of, and the main difficulty of, Thm. E and the uniqueness of $\mathcal{F}_{\mathrm{pcf}}$ in Thm. A. Statements in this level of generality are desirable when possible, because group invariant hyperconvex Frenet curves induce hyperconvex Frenet curves of interest in lower-dimensional flag manifolds that do not come equipped with natural group actions. These curves are typically obtained with versions of the Frenet Restriction Lemma [24, Prop. 3.2] or the Frenet Projecton Lemma 4.9 below. The boundaries of the convex domains foliating $\Omega_\rho^1$ are a first example of this phenomenon.

1.1.2. *Trifoliations of $\Omega_\rho^2$.* The natural adaptations of $\mathcal{F}_{\mathrm{pcf}}, \mathcal{G}_{\mathrm{pcf}}$, and $\mathcal{G}_{\mathrm{tcf}}$ to the other component $\Omega_\rho^2$ of $\Omega_\rho$ are 3-to-1 analogues of foliations that we call *trifoliations* (see §3.2.1). In short, a trifoliation of $\Omega_\rho^2$ is given by a 3-sheeted covering $\pi : M' \to \Omega_\rho^2$ and a foliation $\mathcal{F}$ of the cover $M'$ that is compatible with the covering in the sense that the restriction of $\pi$ to any leaf of $\mathcal{F}$ is injective. The point is that a trifoliation is like a foliation but there are three leaves, not one, though any point of $\Omega_\rho^2$.

The need for such a modification can be seen in the Fuchsian case in the model of $\mathbb{RP}^3$ as projective classes of homogeneous degree 3 polynomials on $\mathbb{R}^2$ (e.g. [11, §3] or [23, §3.2.3]). In this model, $\Omega_\rho^1$ consists of the polynomials with one real projective root of multiplicity one. A leaf of $\mathcal{F}_{\mathrm{pcf}}$ consists of all $[f] \in \mathbb{RP}^3$ with a common unique real projective root of multiplicity one. On the other hand, $\Omega_\rho^2$ consists of the $[f] \in \mathbb{RP}^3$ that take 3 distinct real projective roots. The structures induced by $\xi$ on $\Omega_\rho^2$ have a 3-to-1 nature as a consequence.

In general, if $\xi : \partial \Gamma \to \mathcal{F}(\mathbb{R}^4)$ is a hyperconvex Frenet curve, then every point in $\Omega_\xi^2$ is contained in exactly three planes of the form $\xi^3(x)$ for $x \in \partial \Gamma$. The intersection of each such plane $\xi^3(x)$ with $\Omega_\rho^2$ is the complement of the union of the closure $\overline{C_x}$ of a properly convex domain in $\xi^3(x)$ and a supporting line to $C_x$. Accordingly, following [25], we define:

**Definition** (Concave Domain). *A domain $\Omega$ in $\mathbb{RP}^2$ is* concave *if $\mathbb{RP}^2 - \Omega$ is the union of the closure of a properly convex domain $\Psi$ and a supporting line to $\partial\Psi$.*

Concave domains also arise in the study of flag-manifold structures for SL$_3(\mathbb{R})$-Hitchin representations, and are central objects in [25]. The intersections $\xi^3(x) \cap \Omega_\rho^2$ form a trifoliation $\mathcal{F}_{\text{ccf}}$ of $\Omega_\rho^2$ by concave domains. In §3 we construct two codimension-2 trifoliations $\mathcal{G}_{\text{ctaf}}$ and $\mathcal{G}_{\text{ctrf}}$ of $\Omega_\rho^2$ by projective line segments. Of these, $\mathcal{F}_{\text{ccf}}$ and $\mathcal{G}_{\text{ctaf}}$ are in [11] and $\mathcal{G}_{\text{ctrf}}$ is new.

**Theorem F** (Exactly Two: $\Omega_\rho^2$). *Let $\rho : \Gamma \to \text{PSL}_4(\mathbb{R})$ be Hitchin. Then $\mathcal{G}_{\text{ctaf}}$ and $\mathcal{G}_{\text{ctrf}}$ are the only $\rho(\Gamma)$-invariant trifoliations of $\Omega_\rho^2$ by properly embedded projective line segments.*

*Let $\xi$ be an arbitrary hyperconvex Frenet curve in $\mathcal{F}(\mathbb{R}^4)$. Then $\mathcal{F}_{\text{ccf}}$ is the only trifoliation of $\Omega_\rho^2$ by properly embedded concave domains in projective planes.*

1.2. **Context and Related Work.** This work falls into the project of understanding the geometry of Hitchin representations. This question was suggested by Hitchin in the paper [13] that initiated higher Teichmüller theory and has seen a great deal of development, e.g. [1, 5, 12, 15]. General results in this direction have notable non-qualitative aspects. For instance, even for Hit$_5(S)$ the topological type of the manifolds where the relevant geometric structures are defined is unknown.

Our focus is on the less developed study of the *qualitative* theory of geometric structures associated to Hitchin representations. For Hitchin representations in PSL$_n(\mathbb{R})$, which are fundamental examples of distinguished interest, the basic theory of the case $n = 3$ in $\mathbb{RP}^2$ is well-understood [3, 8] and the theory in the full-flag manifold $\mathcal{F}(\mathbb{R}^3)$ is the subject of [25]. The highest-rank case that has seen progress is $n = 4$, which was the subject of [11]. The only other works focusing on the qualitative projective geometry of all PSL$_4(\mathbb{R})$-Hitchin representations are [11] and [23].

Much of the author's motivation for this work is to develop techniques and understanding that could be useful in approaching the qualitative projective geometry of general Hitchin representations in PSL$_{2n}(\mathbb{R})$. In particular, when one attempts to describe the structure induced by general Hitchin representations, subtle questions arise on the limits of the convexity and general position properties of objects associated to hyperconvex Frenet curves. Sections 3 and 4 here study such questions in the case $2n = 4$.

We also remark that analytic methods have been used to give geometric descriptions of representations in higher Teichmüller spaces for *rank 2* Lie groups, such as in [2, 4, 6, 20, 22]. These methods rely on uniqueness phenomena for equivariant minimal surfaces the rank 2 setting [4, 17, 18] that is false in all greater ranks [19, 27]. Because of this non-uniqueness of minimal surfaces in higher rank, a contrast between our setting and the results expected from adaptations of these analytic methods to higher rank is the *rigidity* we establish of the foliations that appear in our setting (c.f. the discussion in §6 below).

The foliations and domains that we study are connected to other objects of interest, which we now describe. First, the leaves of $\mathcal{F}_{\text{pcf}}$ for non-Fuchsian representations $\rho$ exhibit rich projective-geometric phenomena [23]. It is appealing that they are intrinsically distinguished objects from the geometry of $\Omega_\rho^1$, as is established by the codimension-1 result of Thm. A.

Next, the structure of a family of objects, including the boundary $\partial\Omega_\rho$ of $\Omega_\rho$ here, is salient in the only known construction of canonical complete metrics on Hit$_n(S)$ for $n > 3$ [24]. It seems likely that fine structural results on $\partial\Omega_\rho$, such as those proved here, will be useful in understanding these metrics.

Furthermore, foliations of the form we consider are connected to the dynamics of Hitchin representations through Sambarino's refraction flows [28], as discovered in [25]. From the finiteness we prove here, we expect that constructions similar to those in [25] will only be able to produce refraction flows for length functionals that are *positive roots* of $\mathfrak{sl}_n(\mathbb{R})$.

We finally mention that Farre-Pozzetti-Viaggi recently initiated the study of similar foliated theories to those in $\mathrm{PSL}_4(\mathbb{R})$ for hyperconvex representations into $\mathrm{PSL}_n(\mathbb{C})$ in [7].

1.3. **Outline, Sketch of Central Proof.** In this subsection we outline the structure of the paper then give a sketch of the technical core of the paper.

We begin with an outline. Our proofs of Thms. A and F hinge on features of the structure of $\partial\Omega_\rho$. We spend §3 and §4 developing this structure after recalling the essential background on hyperconvex Frenet curves to our methods in §2.

In §3 we record four explicit developing maps associated to Hitchin representations in $\mathrm{PSL}_4(\mathbb{R})$. Two are due to Guichard-Wienhard in [11] and two are new. Basic properties of these developing maps give proofs that each of the foliations and trifoliations that appear in Thms. A and F are genuine foliations and trifoliations, respectively. Properties of these developing maps are essential to later arguments. Many analogues of the basic facts we prove here are plainly false for Hitchin representations in $\mathrm{SL}_3(\mathbb{R})$. The data obtained from this section is useful in understanding the contrasts of our setting to that of $\mathrm{SL}_n(\mathbb{R})$ for odd $n$.

We then analyze the intersections of $\partial\Omega_\rho$ with planes in $\mathbb{RP}^3$ and some projections of $\xi^1(\partial\Gamma)$ onto planes in affine charts in §4. The case $P = (\xi^3(x) \cap \xi^3(z)) \oplus \xi^1(y)$ of Thm. B is the most difficult to analyze because of the mixing of intersection and sum in its definition. In addressing this case, we prove an analogue (Lemma 4.9) of the Frenet Restriction Lemma [24, Prop. 3.2], which shows certain projections of general hyperconvex Frenet curves to subspaces remain hyperconvex Frenet curves. This leads to new distinguished foliations of the leaves of $\mathcal{F}_{\mathrm{pcf}}$ by convex arcs, which provide the structure used in our arguments.

Thm. B and the rigidity of $\mathcal{F}_{\mathrm{pcf}}$ and $\mathcal{F}_{\mathrm{ccf}}$ claimed in Thms. A, E, and F follow directly from this section. The main tools in §4 are Guichard's work [9] and the extra data provided by our knowledge of developing maps. In particular, [9] is often, but not always, useful in proving that the arcs mentioned in Thm. B are convex. The data obtained from developing maps is then used to determine how the convex arcs that appear glue together. A range of behaviors occur, as is seen in Figure 3, for instance.

1.3.1. *Sketch of Codimension-2 Classification for $\Omega_\rho^1$.* We now give a detailed sketch of the proof of the codimension-2 claim in the first paragraph of Thm. A, which is the heart of the paper. In this sketch, we suppress standard background and numerous technical subtleties.

Fix a $\rho(\Gamma)$-invariant foliation $\mathcal{G}$ of $\Omega_\rho^1$ whose leaves are properly embedded line segments. We show $\mathcal{G} = \mathcal{G}_{\mathrm{pcf}}$ or $\mathcal{G} = \mathcal{G}_{\mathrm{tcf}}$. The outline of the proof is to examine the endpoints of leaves on $\partial\Omega_\rho$ and to apply $\rho(\Gamma)$-invariance of $\mathcal{G}$ to force $\mathcal{G}$ to contain enough leaves of $\mathcal{G}_{\mathrm{pcf}}$ or $\mathcal{G}_{\mathrm{tcf}}$ to deduce that $\mathcal{G} = \mathcal{G}_{\mathrm{pcf}}$ or $\mathcal{G} = \mathcal{G}_{\mathrm{tcf}}$.

There are two immediately relevant peices of basic structure of $\partial\Omega_\rho$. First, $\partial\Omega_\rho = \bigsqcup_{x\in\partial\Gamma} \xi^2(x)$. Second, a homeomorphism $\partial\Gamma^2 \to \partial\Omega_\rho$ is given by sending $(x, y)$ to $\xi^2(x) \cap \xi^3(y)$ off the diagonal in $\partial\Gamma \times \partial\Gamma$ and $(x, x) \mapsto \xi^1(x)$ on the diagonal. Write $\xi^1(\partial\Gamma) = \Xi^1$ and $\Xi^2 = \partial\Omega_\rho - \Xi^1$. The singular locus in $\partial\Omega_\rho$ that is visible in Figure 1 is exactly $\Xi^1$ (§4).

The simple way we have to force leaves of $\mathcal{G}_{\mathrm{pcf}}$ or $\mathcal{G}_{\mathrm{tcf}}$ into $\mathcal{G}$ is as follows. Suppose $\mathcal{G}$ contains a leaf with an endpoint in $\Xi^2$ of the form $\xi^2(x) \cap \xi^3(y)$ with $x = \gamma^-$ a repelling fixed-point of some $\gamma \in \Gamma$. Then an application of $\rho(\Gamma)$-invariance shows that $\mathcal{G}$ also contains a leaf of $\mathcal{G}_{\text{pcf}}$ with endpoints $\xi^1(\gamma^+)$ and $\xi^2(\gamma^-) \cap \xi^3(\gamma^+)$. Something similar happens when $\mathcal{G}$ contains a leaf with an endpoint of the form $\xi^1(\gamma^-)$.

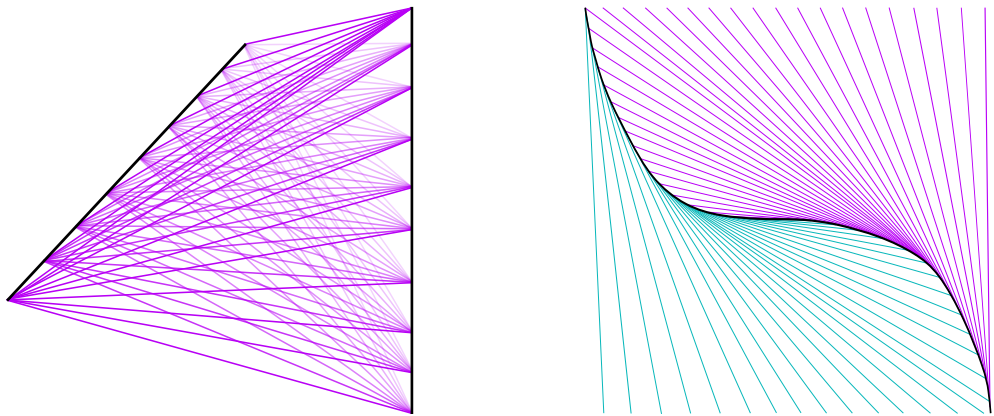

FIGURE 4. Left a peice of a group-invariant foliation that provides no leverage in the proof of Thm. A. Right: the canonical foliation of $\partial\Omega_\rho$ by projective lines, as seen by a line pointing directly into the singular locus.

We now describe what useful constraints we get from these inclusions of leaves and why this is not enough to conclude the main result on its own. The above is a moderately effective way to force $\mathcal{G}$ to contain leaves of $\mathcal{G}_{\text{pcf}}$ or $\mathcal{G}_{\text{tcf}}$. When combined with standard facts on Gromov boundaries, it gives the following. Let $\mathbb{D}$ be a cross-section of a foliation chart for $\mathcal{G}$ and $\phi_F : \mathbb{D} \to \partial\Omega_\rho$ be the map given by taking forward endpoints of leaves. Let $\mathscr{C}$ be a connected subset of $\mathbb{D}$ that is not a point and on which $\phi_F$ is continuous. Then if $\phi_F(\mathscr{C}) \subset \Xi^1$ and $\phi_F(\mathscr{C})$ is not a single point, then $\mathcal{G} = \mathcal{G}_{\text{pcf}}$ or $\mathcal{G} = \mathcal{G}_{\text{tcf}}$. Similarly, if $\phi_F(\mathscr{C}) \subset \Xi^2$ and is not contained in a single line of the form $\xi^2(x)$, then $\mathcal{G} = \mathcal{G}_{\text{pcf}}$.

Unfortunately, this is not quite enough to conclude the desired rigidity. Namely, no inclusions of leaves of $\mathcal{G}_{\text{pcf}}$ or $\mathcal{G}_{\text{tcf}}$ in $\mathcal{G}$ may be deduced by the above means from portions of $\mathcal{G}$ where endpoints are in $\Xi^2$ and have constant $\xi^2$-coordinates. This can occur locally, as illustrated in Figure 4, Left. There is an important asymmetry in the roles of $x$ and $y$ in endpoints of leaves that have the form $\xi^2(x) \cap \xi^3(y)$ here. In particular, the computation used to ensure the inclusion of a leaf of $\mathcal{G}_{\text{pcf}}$ when $x = \gamma^-$ completely breaks when instead only $y = \gamma^-$, and no similar conclusion follows (§5.1, Rmk. 5.5).

This phenomenon obstructs any local argument for Thm. A. So we must seek a more involved source of rigidity.

We now describe this rigidity. To focus on the main case of endpoints of leaves, suppose $\mathcal{G}$ has a leaf $\ell$ with both endpoints in $\Xi^2$. We will show that this is impossible.

So to begin, let $\mathcal{G}$ be a purported $\rho(\Gamma)$-invariant foliation with such a leaf $\ell$ and seek a contradiction. Because of the possibility that $\mathcal{G}$ has leaves with locally-constant $\xi^2$-coordinates near $\ell$, we examine the portion of $\mathcal{G}$ near a different, carefully chosen leaf. To produce this improved leaf, we use the uniform covergence group property of $\Gamma$ to produce a sequence $\gamma_n \in \Gamma$ so $\ell_0 = \lim_{n\to\infty} \gamma_n \ell$ is a leaf of $\mathcal{G}_{\text{pcf}}$ contained in $\mathcal{G}$. This use of the uniform convergence group property of $\Gamma$ is the global structure of the action of $\Gamma$ used in this case.

The structure of $\partial\Omega_\rho$ and its ruling, as seen from an observer on this leaf $\ell_0$ facing towards $\Xi^1$, is illustrated in Figure 4, Right. The point of changing perspective to $\ell_0$ is that the leaf $\ell_0$ cannot be inside of a problematic local model as illustrated in Figure 4 Left, though it can be arbitrarily close to portions of $\mathcal{G}$ with this problematic model.

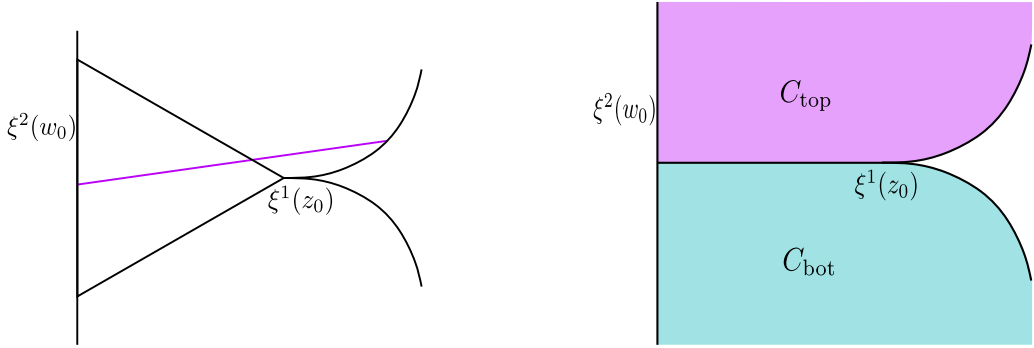

FIGURE 5. Left: if $\mathcal{G}$ contains the leaves at the top and bottom of this triangle and restricts to a local foliation of this plane, the triangle traps leaves with endpoints on $\xi^2(w_0)$. This substantially constrains the shape of $\mathcal{K}$ in $\mathbb{D}$. Right: the structure of the cusp of in $\partial\Omega_\rho \cap (\xi^2(w_0) \oplus \xi^1(z_0))$ at $\xi^1(z_0)$ leads to the key observation that certain unions of $\ell_0$ with segments of $\xi^2(w_0)$ and $\partial\Omega_\rho \cap (\xi^2(w_0) \oplus \xi^1(z_0))$ bound convex domains.

The core of our argument is to over-constrain $\mathcal{G}$ near $\ell_0$. To set notation, from the definition of $\mathcal{G}_{\mathrm{pcf}}$ in §3, there are $z_0, w_0 \in \partial\Gamma$ so that $\ell_0$ is inside the projective line $(\xi^2(w_0) \cap \xi^3(z_0)) \oplus \xi^1(z_0)$. Take a cross-section $\mathbb{D}$ of a foliation chart for $\mathcal{G}$ that intersects $\ell_0$. Refer to *forward endpoints* as those in the direction of $\Xi^1$ on $\ell_0$.

Concerns over continuity of endpoints of leaves through $\mathbb{D}$ come up when endpoints leave or enter $\Xi^1$. We do, however, have *some* structure from the constraints imposed by our basic construction. Namely, $\mathbb{D}$ has a decomposition $\mathbb{D} = \mathscr{O} \cup \mathscr{K}$ by whether forward endpoints are in $\Xi^2$ (for $\mathscr{O}$) or $\Xi^1$ (for $\mathscr{K}$). Our basic structural constraints imply that $\mathscr{K}$ is closed and every leaf through a given connected component of $\mathscr{K}$ shares a common forward endpoint in $\Xi^1$. Similarly, every leaf through a given connected component $\mathscr{C}$ of $\mathscr{O}$ lies in the part in $\Xi^2$ of a common line $\xi^2(x_\mathscr{C})$ in the ruling of $\partial\Omega_\rho$.

In contrast, backwards endpoints of leaves through $\mathbb{D}$ are much better controlled. They are continuous on $\mathbb{D}$, and the leaf inclusions implied by our basic considerations force all backwards endpoints of leaves through $\mathbb{D}$ to be contained in the *single* line $\xi^2(w_0)$.

This last characteristic is a major constraint. This constraint forces all leaves of $\mathcal{G}$ through $\mathbb{D}$ to be in the one-parameter family of planes containing $\xi^2(w_0)$. All such planes have the form $\xi^2(w_0) \oplus \xi^1(y)$ or $\xi^3(w_0)$, and distinct such planes intersect only in $\xi^2(w_0)$. The key point is that *this forces $\mathcal{G}$ to have a local product structure near $\ell_0$, in which $\mathcal{G}$ restricts to local foliations of planes of the form $\xi^2(w_0) \oplus \xi^1(y)$ ($y \in \partial\Gamma - \{w_0\}$)*. These planes have a consistent and highly-structured form, illustrated in Figure 3, second from Left. Intersections of the cross-section $\mathbb{D}$ with this one-parameter family of planes are leaves of an auxiliary foliation $\mathscr{D}$ of $\mathbb{D}$ by line segments that is useful for book-keeping.

This geometrically-adapted local product structure allows us to deduce some features of $\mathcal{G}$ through considering foliations of *plane* regions by line segments. Namely, working in this lower-dimensional setting, geometric considerations lead to two further constraints on the set $\mathscr{K}$ of points in $\mathbb{D}$ with forward endpoints in $\Xi^1$:

(1) $\mathscr{K}$ is a union of segments of leaves of $\mathscr{D}$, with at most one segment per leaf of $\mathscr{D}$. This is because if two leaves in $\xi^2(w_0) \oplus \xi^1(y)$ have forward-endpoint $\xi^1(y)$, to avoid intersections of leaves, the foliation $\mathcal{G}$ must contain the unique foliation of the entire triangle they bound by segments ending on $\xi^1(y)$. See Figure 4, Left.
(2) The forward endpoints on the leaf $L_0$ of $\mathscr{D}$ that intersects $\ell_0$ are continuous, unlike the possible behavior of forward endpoints on general leaves of $\mathscr{D}$. This is because $\ell_0$ cuts the nearby region in $\Omega_\rho^1 \cap (\xi^2(z_0) \oplus \xi^1(w_0))$ into two convex domains, as it is tangent to the cusp at $\xi^1(z_0)$. See Figure 4, Right. This implies $L_0 \subset \mathscr{K}$.

We next consider the subset $\mathscr{K}'$ of $\mathscr{K}$ consisting of the leaves of $\mathscr{D}$ that are entirely contained in $\mathscr{K}$. Note $\mathscr{K}'$ is nonempty as $L_0 \subset \mathscr{K}'$. From basic considerations, $\mathscr{K}'$ is closed and has empty interior. Classical line and plane topology ensures we may change $\ell_0$ to a nearby leaf so that $L_0$ is an endpoint of an interval in the open subset of the leaf space of $\mathscr{D}$ complementary to $\mathscr{K}'$. We may do this while retaining all relevant properties of $\ell_0$. Let $\mathscr{U}_0$ be an adjacent open set in $\mathscr{O}$ to $L_0$, with forward endpoints all in a single line $\xi^2(x_{\mathscr{U}_0})$. The structure of $\xi^3(z_0) \cap \partial\Omega_\rho$ forces $x_{\mathscr{U}_0} \neq z_0$: the line segments that would need to be included in $\mathcal{G}$ in order for $x_{\mathscr{U}_0}$ to be $z_0$ cannot escape the convex domain in $\xi^3(z_0)$.

Now $\mathcal{G}$ is at last over-constrained near $L_0$: the structure of the "front view" in Figure 4, Right and continuity considerations on $\mathscr{U}_0$ force $\mathcal{G}$ to contain a leaf through $L_0$ with an endpoint in $\xi^2(x_{\mathscr{U}_0})$, which contradicts that $L_0 \in \mathscr{K}'$.

**Acknowledgements.** I am grateful to Olivier Guichard and Max Riestenberg for stimulating conversations, to Parker Evans for conversations on the relation of Collier-Toulisse-Tholozan's work to this setting and comments on an earlier draft, to Mike Wolf for his support and interest, and to Katherine Booth for the colorblind-friendly color palette used in illustrations. This material is based on work supported by the National Science Foundation under Grants No. 1842494 and 2005551.

## 2. Notations & Reminders

**2.1. Gromov Boundaries of Surface Groups.** Let $\Gamma$ be the fundamental group of a closed, orientable surface $S$ of genus at least 2. Let us recall some standard features (e.g. [14]) of the Gromov boundary $\partial\Gamma$ of $\Gamma$, which is homeomorphic to the circle and has a natural action of $\Gamma$.

Let $\partial\Gamma^{(2)}$ be the collection of pairs $(x, y)$ of distinct elements of $\partial\Gamma$. Fix an orientation, and with it a circular order, on $\partial\Gamma$ throughout the rest of the paper. Let $\partial\Gamma^{(3)+}$ be the collection of positively oriented triples in $\partial\Gamma$.

We recall that the action of $\Gamma$ on $\partial\Gamma$ has *North-South dynamics*, meaning that every $\gamma \in \Gamma - \{e\}$ has two fixed-points $\gamma^+, \gamma^- \in \partial\Gamma$ so that for every $x \in \partial\Gamma - \{\gamma^-\}$ we have $\lim_{n\to\infty} \gamma^n x = \gamma^+$. For $\gamma \in \Gamma - \{e\}$, we call the associated fixed-point pair $(\gamma^-, \gamma^+)$ the *pole-pair* associated to $\gamma$. Pole-pairs are dense in $\partial\Gamma^{(2)}$.

A useful global feature of the action of $\Gamma$ on $\partial\Gamma$ is that $\Gamma$ is a *uniform convergence group*. The formulation of this we shall use is that every point in $\partial\Gamma$ is a conical limit point of $\Gamma$ in the sense that for any $x \in \partial\Gamma$ there is a sequence $\gamma_n \in \Gamma$ and distinct $a, b \in \partial\Gamma$ so that $\lim_{n\to\infty} \gamma_n x = a$ and $\lim_{n\to\infty} \gamma_n y = b$ for all $y \neq x$ in $\partial\Gamma$.

The product structure of $\partial\Gamma^3$ induces natural foliations of $\partial\Gamma^{(3)+}$ that are mapped to our foliations of $\Omega_\rho^1$ and trifoliations of $\Omega_\rho^2$ by developing maps in §3. To this end, let $\mathcal{F}, \mathcal{H}$, and $\mathcal{G}$ be the foliations of $\partial\Gamma^{(3)+}$ arising from the product structure of $\partial\Gamma^3$ with:

(1) Leaves $f_{x_0}$ of $\mathcal{F}$ given by $\{(x,y,z) \in \partial\Gamma^{(3)+} \mid x = x_0\}$ $(x_0 \in \partial\Gamma)$,
(2) Leaves $h_{y_0}$ of $\mathcal{H}$ given by $\{(x,y,z) \in \partial\Gamma^{(3)+} \mid y = y_0\}$ $(y_0 \in \partial\Gamma)$,
(3) Leaves $g_{x_0 y_0}$ of $\mathcal{G}$ given by $\{(x,z,y) \in \partial\Gamma^{(3)+} \mid x = x_0, y = y_0\}$ $((x_0, y_0) \in \partial\Gamma^{(2)})$,

Though we shall not directly use it in the following, it is useful for motivation to recall the well-known fact that fixing a hyperbolic metric on $S$ induces an identification of $\partial\Gamma^{(3)+}$ with the unit tangent bundle $\mathrm{T}^1\mathbb{H}^2$ of the hyperbolic plane $\mathbb{H}^2$, under which $\mathcal{G}$ is mapped to the geodesic foliation and $\mathcal{F}$ is mapped to the stable foliation (e.g. [11, §2]).

2.2. **Hitchin Representations and Hyperconvex Frenet Curves.** We recall the essential background on hyperconvex Frenet curves for the following. In particular, we review results of Guichard on improvements that can be made to the general position and limit compatibility properties of hyperconvex Frenet curves [9].

For a vector space $V$, let $\mathcal{F}(V)$ denote the manifold of full flags of nested proper, nontrivial subspaces $(V_1, V_2, ..., V_{n-1})$ in $V$ and let $\mathrm{Gr}_k(V)$ be the Grassmannian of $k$-planes in $V$.

The hyperconvex Frenet curve condition is phrased in terms of sums of subspaces:

**Definition 2.1.** *Let $I$ be a closed interval, open interval, or $\partial\Gamma$. A continuous map $\xi : I \to \mathcal{F}(\mathbb{R}^n)$ is a* hyperconvex Frenet curve *if:*

(1) (General Position). *For any integers $k_1, ..., k_j$ with $\sum_{l=1}^{j} k_l = p \leq n$, and distinct points $x_1, ..., x_j \in I$, the sum $\xi^{k_1}(x_1) + ... + \xi^{k_j}(x_j)$ is direct.*
(2) (Limit Compatibility). *For any $x \in I$, indices $k_1, ..., k_j$ as above, and any sequence $(x_1^m, ..., x_j^m)$ of $j$-tuples of distinct points in $I$ that converge to $(x, ..., x)$,*

$$\xi^p(x) = \lim_{m \to \infty} \bigoplus_{l=1}^{j} \xi^{k_l}(x_l^m).$$

The general position condition (1) is also referred to as *hyperconvexity*. To prove that a curve is a hyperconvex Frenet curve, it suffices to verify the general position condition for the case $p = n$ and the general limit compatibility condition.

The standard examples of hyperconvex Frenet curves are those in $\mathcal{F}(\mathbb{R}^3)$. In this case, a continuous map $(\xi^1, \xi^2)$ from $\partial\Gamma$ to $(\mathbb{RP}^2, \mathrm{Gr}_2(\mathbb{R}^3))$ is a hyperconvex Frenet curve if and only if $\xi^1$ is a homeomorphism of $\partial\Gamma$ onto the boundary $\partial C_\xi$ of a properly convex, strictly convex $C^1$ domain $C_\xi$ in $\mathbb{RP}^2$ and $\xi^2(x)$ is the tangent line to $\partial C_\xi$ at $\xi^1(x)$ for all $x \in \partial\Gamma$.

The relevance of hyperconvex Frenet curves to our setting is the following characterization of Hitchin representations in terms of hyperconvex Frenet curves:

**Theorem 2.2** (Labourie [16, Thm. 1.4], Guichard [10, Thm. 1]). *A representation $\rho : \Gamma \to \mathrm{PSL}_n(\mathbb{R})$ is Hitchin if and only if there exists a $\rho$-equivariant hyperconvex Frenet curve $\xi_\rho : \partial\Gamma \to \mathcal{F}(\mathbb{R}^n)$.*

Thm. 2.2 shall be our working definition of Hitchin representations.

The hyperconvex Frenet curve condition behaves well with respect to projective duality, which is the map $\perp: \mathrm{Gr}_k(\mathbb{R}^n) \to \mathrm{Gr}_{n-k}((\mathbb{R}^n)^*)$ that takes a subspace $V$ of $\mathbb{R}^n$ to the subspace $V^\perp$ of $(\mathbb{R}^n)^*$ consisting of linear functionals that vanish on $V$. Given a map $\xi = (\xi^1, ..., \xi^{n-1}) : I \to \mathcal{F}(\mathbb{R}^n)$, duality produces a map $\xi^* = ((\xi^{n-1})^\perp, ..., (\xi^1)^\perp) : I \to \mathcal{F}((\mathbb{R}^n)^*)$, that we shall write as $\xi^* = (\xi^{*1}, \xi^{*2}, ..., \xi^{*n-1})$.

**Theorem 2.3** (Guichard [9, Thm. 2]). *A continuous map $\xi : I \to \mathcal{F}(\mathbb{R}^n)$ is a hyperconvex Frenet curve if and only if the dual map $\xi^*$ is a hyperconvex Frenet curve.*

It is often useful to consider dual hyperconvex Frenet curves in order to interchange the roles of sum and intersection. In particular, the hyperconvex Frenet condition can be phrased (and verified) in terms of intersections:

**Proposition 2.4** (Guichard [9, Lemme 5]). *Let $I$ be a closed interval, open interval, or $S^1$. A continuous map $\xi : I \to \mathcal{F}(\mathbb{R}^n)$ is a hyperconvex Frenet curve if and only if:*

(1) (General Position). *For any integers $k_1, ..., k_j$ with $\sum_{l=1}^{j} k_l = p \leq n$, and distinct points $x_1, ..., x_j \in I$, the intersection $\xi^{n-k_1}(x_1) \cap ... \cap \xi^{n-k_j}(x_j)$ has dimension $n - p$.*
(2) (Limit Compatibility). *For any $x \in I$, indices $k_1, ..., k_j$ as above, and sequence $(x_1^m, ..., x_j^m)$ of $j$-tuples of distinct points in $I^j$ that converge to $(x, ..., x)$, we have $\xi^{n-p}(x) = \lim_{m \to \infty} \bigcap_{l=1}^{j} \xi^{n-k_l}(x_l^m)$.*

In proving that a curve is a hyperconvex Frenet curve from this perspective, it suffices to verify the general position condition for $p = n$ and the general limit compatibility condition.

As it happens, the hyperconvex Frenet curve condition implies that far more combinations of subspaces than those mentioned in its definition are in as general position as possible. The following is such an improvement that we shall use repeatedly.

**Lemma 2.5** (Guichard [9, Prop. 6]). *If $I$ is a closed interval, open interval, or $S^1$, and $\xi = (\xi^1, ..., \xi^{n-1})$ is a hyperconvex Frenet curve $\xi : I \to \mathcal{F}(\mathbb{R}^n)$, then:*

(1) (Improved General Position) *For all $k$-tuples $(m_1, ..., m_k)$ and $l$-tuples $(p_1, ..., p_l)$ of natural numbers so that $\sum_{i=1}^{k} m_i = \sum_{i=1}^{l} p_i \leq n$ and points $x_1, ..., x_k, y_1, ..., y_l$ so that for some index $g$, we have $y_1 < ... < y_g < x_1 < ... < x_k < y_{g+1} < ... < y_l$, then the following sum is direct:*

$$\bigoplus_{i=1}^{l} \xi^{p_i}(y_i) \oplus \left(\bigcap_{i=1}^{k} \xi^{n-m_i}(x_i)\right) = \mathbb{R}^n.$$

(2) (Improved Limit Compatibility) *For all $x \in I$, $k$-tuples $m_1, ..., m_k$ and $l$-tuples $p_1, ..., p_l$ in $\mathbb{N}$ so $\sum_{i=1}^{l} p_i = p \leq \sum_{i=1}^{k} m_i = m \leq n$, and sequences of the form $(y_1^j, ..., y_g^j, x_1^j, ..., x_k^j, y_{g+1}^j, ..., y_l^j)$ of positively ordered $(l+k)$-tuples in $I$ that converge to $(x, ..., x)$,*

$$\lim_{j \to \infty} \bigoplus_{i=1}^{l} \xi^{p_i}(y_i^j) \oplus \left(\bigcap_{i=1}^{k} \xi^{n-m_i}(x_i^j)\right) = \xi^{n-m+p}(x).$$

The order conditions in Lemma 2.5 are crucial. We shall encounter situations where consequences of Lemma 2.5 are true exactly as far as the order hypotheses above hold, and false any further.

**Remark 2.6.** *There is a rather efficient notation in much of the literature, in which when a hyperconvex Frenet curve $\xi$ is understood from context, $\xi^k(x)$ is written as $x^k$ or $x_\rho^k$ (e.g. [26]). This notation improves the readability of some expressions and computations, but is not well-suited to our uses because we consider a large number of different hyperconvex Frenet curves that are induced by an original hyperconvex Frenet curve. We do not use this more efficient notation because it is often ambiguous in our setting.*

## 3. Developing Maps and Invariant Foliations

In this section we give four constructions of projective structures on $\partial \Gamma^{(3)+}$ induced by a hyperconvex Frenet curve $\xi : \partial \Gamma \to \mathcal{F}(\mathbb{R}^4)$ whose developing images are connected components of $\mathbb{RP}^3 - \bigsqcup_{x \in \partial \Gamma} \xi^2(x)$. Each construction maps canonical foliations of $\partial \Gamma^{(3)+}$ onto objects in $\mathbb{RP}^3$ that are natural under projective transformations. These give rise to the distinguished foliations and trifoliations of our main results.

Throughout the next two sections we consider hyperconvex Frenet curves $\xi : \partial \Gamma \to \mathcal{F}(\mathbb{R}^4)$ that are not necessarily invariant under Hitchin representations. For any such hyperconvex Frenet curve $\xi$, define $\partial \Omega_\xi = \bigsqcup_{x \in \partial \Gamma} \xi^2(x)$ and $\Omega_\xi = \mathbb{RP}^3 - \partial \Omega_\xi$.

The first two constructions are due to Guichard-Wienhard in [11]; the second two are new. We briskly review Guichard-Wienhard's constructions for the sake of self-containment and completeness, then discuss the new constructions in more detail. That these developing maps are covering maps of one or three sheets onto the connected components $\Omega_\xi^1$ and $\Omega_\xi^2$ of $\Omega_\xi$ is crucial in our study of the structure of $\partial \Omega_\xi$.

### 3.1. Properly Convex Foliated Projective Structures.

Let $\xi : \partial \Gamma \to \mathcal{F}(\mathbb{R}^4)$ be a hyperconvex Frenet curve. For any $x \in \partial \Gamma$, define the maps $\xi_x^1 : \partial \Gamma \to \mathbb{P}(\xi^3(x))$ and $\xi_x^2 : \partial \Gamma \to \mathrm{Gr}_2(\xi^3(x))$ by

$$\xi_x^1(y) = \begin{cases} \xi^2(y) \cap \xi^3(x) & y \neq x \\ \xi^1(x) & y = x \end{cases}, \qquad \xi_x^2(y) = \begin{cases} \xi^3(y) \cap \xi^3(x) & y \neq x \\ \xi^2(x) & y = x. \end{cases}$$

The two argument map $\xi_\cdot^1(\cdot) : \partial \Gamma^2 \to \partial \Omega_\xi$ is an equivariant homeomorphism onto $\partial \Omega_\rho$ [11, §4]. Guichard-Wienhard note in [11] that $\xi_x = (\xi_x^1, \xi_x^2)$ is a hyperconvex Frenet curve in $\xi^3(x)$. This is a reflection of a general phenomenon exhibited by hyperconvex Frenet curves that we shall use in a number of settings.

**Lemma 3.1** ([24, Prop. 3.2]: Frenet Restriction). *Let $\xi : \partial \Gamma \to \mathcal{F}(\mathbb{R}^n)$ be a hyperconvex Frenet curve and fix $1 < D < n$ and $x_0 \in \partial \Gamma$. Then $\xi_{x_0,D} = (\xi_{x_0,D}^1, ..., \xi_{x_0,D}^{D-1}) : \partial \Gamma \to \mathcal{F}(\xi^D(x_0))$ defined by*

$$\xi_{x_0,D}^k(y) = \begin{cases} \xi^{n-D+k}(y) \cap \xi^D(x_0) & y \neq x_0 \\ \xi^k(x_0) & y = x_0 \end{cases}$$

*is a hyperconvex Frenet curve.*

Next, define $\mathrm{dev}_{\mathrm{pcf}} : \partial \Gamma^{(3)+} \to \mathbb{RP}^3$ by

$$\mathrm{dev}_{\mathrm{pcf}}(x, y, z) = (\xi^1(x) \oplus \xi_x^1(z)) \cap (\xi_z^1(x) \oplus \xi_x^1(y)).$$

See Figure 6, Left. Guichard-Wienhard show that $\mathrm{dev}_{\mathrm{pcf}}$ is a homeomorphism onto a connected component, that we call $\Omega_\xi^1$, of $\Omega_\xi$ in [11, §4]. It has the properties that:

(1) Every leaf $g_{xy}$ of $\mathcal{G}$ maps to a properly embedded properly convex subset of a projective line in $\mathbb{RP}^3$. A foliation $\mathcal{G}_{\mathrm{pcf}}$ of $\Omega_\rho^1$ is specified by $\mathrm{dev}_{\mathrm{pcf}}$, with leaves $\mathrm{dev}_{\mathrm{pcf}}(g_{xy})$ for $(x, y) \in \partial \Gamma^{(2)}$. See Figure 7, Left.
(2) Every leaf $f_x$ of $\mathcal{F}$ maps to a properly embedded properly convex domain, that we denote by $C_x$ throughout the following, in $\xi^3(x) \subset \mathbb{RP}^3$. A foliation $\mathcal{F}_{\mathrm{pcf}}$ of $\Omega_\rho^1$ is induced by $\mathrm{dev}_{\mathrm{pcf}}$ with leaves $\mathrm{dev}_{\mathrm{pcf}}(f_x)$ for $x \in \partial \Gamma$.

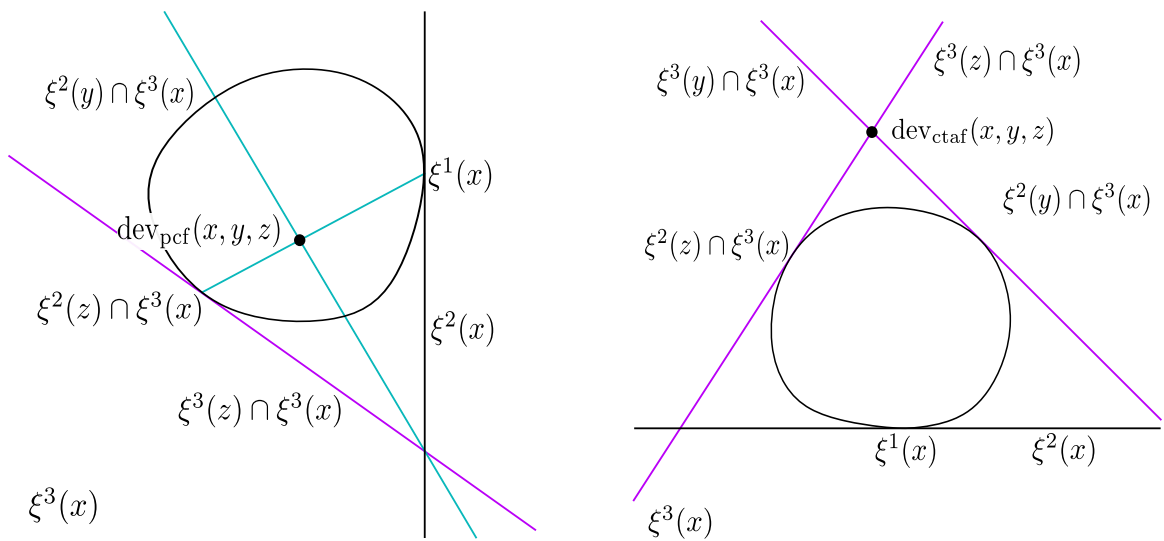

FIGURE 6. From left to right: illustrations of dev$_{\text{pcf}}$ and dev$_{\text{ctaf}}$. These are the developing maps appearing in Guichard-Wienhard's work [11].

When $\xi$ is the limit map of a Hitchin representation $\rho$, then dev$_{\text{pcf}}$ is equivariant with respect to the actions of $\Gamma$ on $\partial \Gamma^{(3)+}$ and $\rho(\Gamma)$ on $\mathbb{RP}^3$, and $(\text{dev}_{\text{pcf}}, \rho)$ defines a Thurston-Klein $(\text{PSL}_4(\mathbb{R}), \mathbb{RP}^3)$-structure on $T^1 S$.

**Remark 3.2.** *In the following, since the entries of $\xi^1_\cdot(\cdot)$ play highly asymmetric roles, we often favor less succinct notation that emphasizes the different roles of entries, with the hope that it makes the arguments clearer.*

### 3.2. Concave Tangent Foliated Projective Structures.
The developing map is simpler in this case. Define dev$_{\text{ctaf}} : \partial \Gamma^{(3)+} \to \mathbb{RP}^3$ by

$$\text{dev}_{\text{ctaf}}(x, y, z) = \xi^3(x) \cap \xi^3(y) \cap \xi^3(z).$$

See Figure 6. As dev$_{\text{pcf}}$ is a homeomorphism onto $\Omega^1_\xi$, every point in $\Omega^1_\xi$ is in exactly one $\xi^3(x)$ for $x \in \partial \Gamma$. Consequently, dev$_{\text{ctaf}}$ has image contained in $\overline{\Omega^2_\rho}$. It is then a direct consequence of hyperconvexity that dev$_{\text{ctaf}}$ has image contained in $\Omega^2_\rho$, i.e. does not intersect $\xi^2(w)$ for any $w \in \partial \Gamma$.

In [11], this map is recorded and a brief indication of how to prove the following is given. We present a detailed argument here for completeness.

**Proposition 3.3.** *The map* dev$_{\text{ctaf}}$ *is a 3-sheeted covering of $\Omega^2_\rho$.*

*Proof.* The claim follows from proving that dev$_{\text{ctaf}}$ is proper, continuous, 3-to-1, and locally injective. Then by invariance of domain, dev$_{\text{ctaf}}$ is a proper local homeomorphism and hence a 3-sheeted covering map.

Continuity is immediate from hyperconvexity of $\xi$. To see local injectivity note that by hyperconvexity of $\xi$, for all $x, y, z, w$ distinct in $\partial \Gamma$,

$$(\xi^3(x) \cap \xi^3(y) \cap \xi^3(z)) \cap \xi^3(w) = \{0\},$$

so that

$$\text{dev}^{-1}_{\text{ctaf}}(\text{dev}_{\text{ctaf}}(x, y, z)) = \{(x, y, z), (z, y, z), (y, z, x)\}.$$

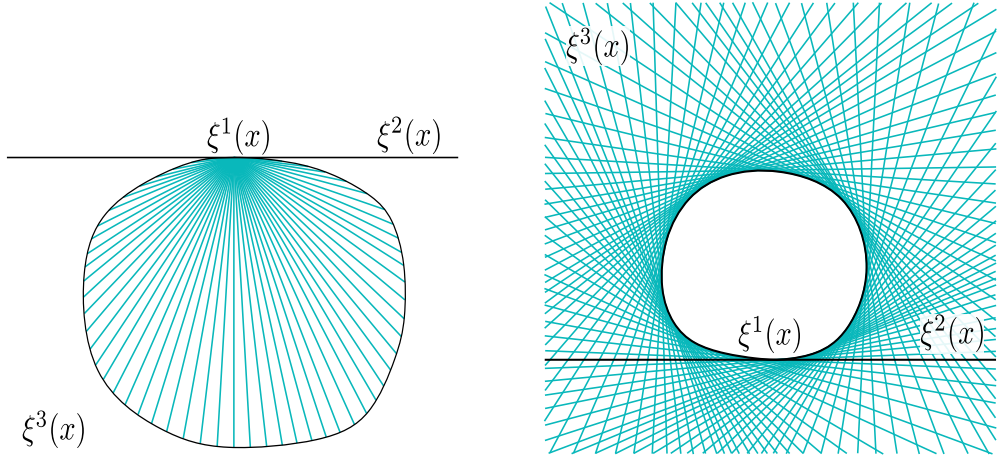

FIGURE 7. Leaves of $\mathcal{G}_{\mathrm{pcf}}$ (left) and $\mathcal{G}_{\mathrm{ctaf}}$ (right) within subspaces of the form $\xi^3(x)$ for some $x \in \partial\Gamma$. The unions of the illustrated leaves give single leaves of the codimension-1 (tri-)foliations $\mathcal{F}_{\mathrm{pcf}}$ and $\mathcal{F}_{\mathrm{ccf}}$, respectively.

We next claim $\mathrm{dev}_{\mathrm{ctaf}}$ is proper. So let $(x_n, y_n, z_n)$ be a sequence in $\partial\Gamma^{(3)+}$ that leaves all compact subsets. To show properness, it suffices to demonstrate that $\mathrm{dev}_{\mathrm{ctaf}}(x_n, y_n, z_n)$ has a subsequence leaving all compact subsets of $\Omega_\xi^2$. Up to subsequence, $(x_n, y_n, z_n)$ converges in $\partial\Gamma^3$ to a limit $(x_\infty, y_\infty, z_\infty)$ with at least two entries equal. We break into cases on $(x_\infty, y_\infty, z_\infty)$.

**Case 1.** $x_\infty = y_\infty = z_\infty$.

In this case, hyperconvexity shows
$$\lim_{n \to \infty} \xi^3(x_n) \cap \xi^3(y_n) \cap \xi^3(z_n) = \xi^1(x_\infty) \in \partial\Omega_\rho^2.$$

**Case 2.** *Exactly two of $x_\infty, y_\infty$, and $z_\infty$ are equal.*

Then without loss of generality take $x_\infty = y_\infty \neq z_\infty$. Then $\xi^3(z_\infty)$ is transverse to $\xi^2(x_\infty)$, and by hyperconvexity $\lim_{n \to \infty} \xi^3(x_n) \cap \xi^3(y_n) = \xi^2(x_\infty)$. We conclude
$$\lim_{n \to \infty} \xi^3(x_n) \cap \xi^3(y_n) \cap \xi^3(z_n) = \xi^2(x_\infty) \cap \xi^3(z_\infty) \in \partial\Omega_\xi^2.$$

So $\mathrm{dev}_{\mathrm{ctaf}}$ is a 3-sheeted covering of $\Omega_\xi^2$. □

**Remark 3.4.** *In the group invariant case, equivariance of developing maps and properness of the action of $\rho(\Gamma)$ on $\Omega_\rho$ may be used to give alternative proofs of properness of the developing maps in this section. One of the points of this section is that the developing maps we construct are valid without the assumption of group-invariance and so work, for instance, for hyperconvex Frenet curves obtained through applications of the Frenet Restriction Lemma.*

3.2.1. *Trifoliations.* That $\mathrm{dev}_{\mathrm{ctaf}}$ is a covering of three sheets leads to the relevant analogues of the foliations of $\Omega_\xi^1$ being 3-to-1 adaptations of foliations. A modification of the definition of a foliation that includes the objects of interest is as follows.

**Definition 3.5** (*j*-foliation). *A codimension-$k$ $j$-foliation $(M', \mathcal{F}, \pi)$ of a manifold $M$ is:*
  (1) *A (not assumed to be connected) manifold $M'$,*
  (2) *A codimension-$k$ foliation $\mathcal{F}$ of $M'$,*

(3) *A $j$-sheeted cover $\pi : M' \to M$ so that the restriction of $\pi$ to any leaf of $\mathcal{F}$ is injective.*

*A $j$-foliation chart for a $j$-foliation $\mathcal{T}$ on an evenly covered neighborhood $U \subset M$ is a map of the form $\phi \circ p$ where $p$ is a local inverse for the covering and $\phi$ is a foliation chart for $p(U)$. Given a $j$-foliation $\mathcal{T} = (M', \mathcal{F}, \pi)$ and a leaf $\ell$ of $\mathcal{F}$ we say $\pi(\ell)$ is a leaf of $\mathcal{T}$. We consider two $j$-foliations of $M$ equivalent if all of their leaves coincide.*

*A $j$-foliation with $j = 3$ is a* trifoliation.

The notion of a $\rho(\Gamma)$-invariant trifoliation of our main results is defined in terms of the following basic notion of pushforward. Namely, given a manifold $M$ the group Homeo($M$) acts on $j$-foliations of $M$ by post-composition of the coverings, e.g. $\varphi(M', \mathcal{F}, \pi) = (M', \mathcal{F}, \varphi \circ \pi)$ for all $j$-foliations $(M', \mathcal{F}, \pi)$ of $M$ and $\varphi \in$ Homeo($M$). So for a Hitchin representation $\rho : \Gamma \to \mathrm{PSL}_4(\mathbb{R})$ and a trifoliation $\mathcal{T}$ of $\Omega_\rho^2$, we say that $\mathcal{T}$ is $\rho(\Gamma)$-*invariant* if $\rho(\gamma)\mathcal{T} = \mathcal{T}$ for all $\gamma \in \Gamma$.

3.2.2. *Consequences.* The above shows that a codimension-2 trifoliation $\mathcal{G}_{\mathrm{ctaf}}$ of $\Omega_\rho^2$ is specified by dev$_{\mathrm{ctaf}}$ with leaves dev$_{\mathrm{ctaf}}(g_{xy})$ for $g_{xy}$ a leaf of $\mathcal{G}$. Each leaf of $\mathcal{G}_{\mathrm{ctaf}}$ is a properly embedded properly convex subset of a projective line in $\mathbb{RP}^3$. See Figure 7, Right.

Similarly, dev$_{\mathrm{ctaf}}$ specifies a codimension-1 trifoliation $\mathcal{F}_{\mathrm{ccf}}$ of $\Omega_\xi^2$ with leaves dev$_{\mathrm{ctaf}}(f_x)$ for $f_x$ a leaf of $\mathcal{F}$. Each leaf of $\mathcal{F}_{\mathrm{ccf}}$ is a concave region in a projective plane. See Figure 7. Again, if $\xi$ is the limit map of a Hitchin representation $\rho$, then dev$_{\mathrm{ctaf}}$ is equivariant and (dev$_{\mathrm{ctaf}}, \rho$) defines a $(\mathrm{PSL}_4(\mathbb{R}), \mathbb{RP}^3)$-structure on T$^1 S$.

The following characterizations of $\Omega_\rho^1$ and $\Omega_\rho^2$ are immediate consequences of the surjectivity of the developing maps in the last two subsections and are useful in the sequel.

**Lemma 3.6** (Domain Characterizations). *We have:*
$$\Omega_\xi^1 = \{p \in \Omega_\xi \mid \text{there is a unique } x \in \partial\Gamma \text{ so } p \in \xi^3(x)\},$$
$$\Omega_\xi^2 = \{p \in \Omega_\xi \mid \text{there are distinct } x, y, z \in \partial\Gamma \text{ so } p \in \xi^3(x) \cap \xi^3(y) \cap \xi^3(z)\}.$$

*Furthermore, for all $(x, y) \in \partial\Gamma^{(2)}$, the intersection $\xi^3(x) \cap \xi^3(y) \cap \Omega_\xi^2$ is nonempty.*

3.3. **Properly Convex Transversely Foliated Projective Structures.** This is the first new construction. Define dev$_{\mathrm{pctf}} : \partial\Gamma^{(3)+} \to \mathbb{RP}^3$ by
$$\mathrm{dev}_{\mathrm{pctf}}(x, y, z) = (\xi^1(x) \oplus \xi^1(z)) \cap \xi^3(y).$$

See Figure 8, Left.

**Proposition 3.7.** *The map* dev$_{\mathrm{pctf}}$ *is a homeomorphism of $\partial\Gamma^{(3)+}$ onto $\Omega_\xi^1$.*

*Proof.* As in the other cases, we shall prove that dev$_{\mathrm{pctf}}$ is a proper continuous injection $\partial\Gamma^{(3)+} \to \Omega_\xi^1$, which suffices for the claim.

Note first that for any distinct $x, w, z \in \partial\Gamma$ that $(\xi^1(x) \oplus \xi^1(z)) \cap \xi^2(w) = \{0\}$ by Lemma 2.5, and $(\xi^1(x) \oplus \xi^1(z)) \cap \xi^3(y)$ is neither $\xi^1(y)$ nor $\xi^2(w) \cap \xi^3(y)$ for any $w \in \partial\Gamma - \{y\}$. We conclude that the image of dev$_{\mathrm{pctf}}$ does not intersect $\partial\Omega_\xi$, and so is contained in $\Omega_\xi$.

We next claim that the image of dev$_{\mathrm{pctf}}$ is contained in $\Omega_\xi^1$. To see this, suppose otherwise and observe that by Lemma 3.6 if dev$_{\mathrm{pctf}}(x, y, z) \in \Omega_\rho^2$ then there are distinct $r, w \in \partial\Gamma - \{y\}$ so that dev$_{\mathrm{pctf}}(x, y, z) = \xi^3(r) \cap \xi^3(w) \cap \xi^3(y)$. Note that $(\xi^1(x) \oplus \xi^1(z)) \cap \xi^3(x) = \xi^1(x) \in \partial\Omega_\xi$, so that $r$ and $w$ are also not $x$. Similarly, $r$ and $w$ are not $z$. Counting, at least two of $r, w$, and $y$ are in the same connected component of $\partial\Gamma - \{x, z\}$. Call these points $y_1$ and $y_2$.

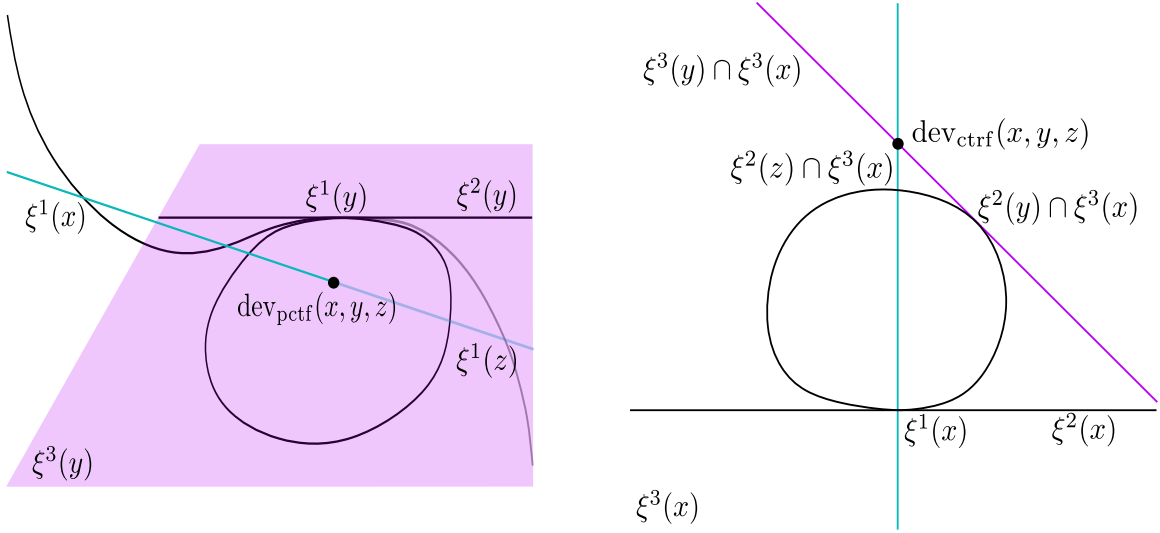

FIGURE 8. From left to right: illustrations of $\text{dev}_{\text{pctf}}$ and $\text{dev}_{\text{ctrf}}$.

Then $x, z, y_1, y_2$ may be arranged in a cyclic order in which $x$ and $z$ are adjacent, so that Lemma 2.5 applies and
$$(\xi^1(x) \oplus \xi^1(z)) \cap \xi^3(y_1) \cap \xi^3(y_2) = \{0\},$$
which is incompatible with the assumption $\text{dev}_{\text{pctf}}(x, y, z) = \text{dev}_{\text{ctaf}}(r, w, y)$. We conclude that the image of $\text{dev}_{\text{pcf}}$ is contained in $\Omega^1_\rho$.

We next turn towards injectivity. Note by the general position property of hyperconvex Frenet curves that for distinct $x, y, z, w \in \partial\Gamma$ the sum $\xi^1(x) + \xi^1(y) + \xi^1(z) + \xi^1(w)$ is direct, so that

(3.1)
$$(\xi^1(x) \oplus \xi^1(y)) \cap (\xi^1(z) \oplus \xi^1(w)) = \{0\},$$
$$(\xi^1(x) \oplus \xi^1(y)) \cap (\xi^1(x) \oplus \xi^1(z)) = \xi^1(x) \in \partial\Omega_\xi.$$

Now assume $\text{dev}_{\text{pctf}}(x_1, y_1, z_1) = \text{dev}_{\text{pctf}}(x_2, y_2, z_2)$. Then Eq. (3.1) shows $\{x_1, z_1\} = \{x_2, z_2\}$. On the other hand, because the image of $\text{dev}_{\text{pctf}}$ is contained in $\Omega^1_\rho$, Lemma 3.6 and the definition of $\text{dev}_{\text{pctf}}$ imply that $y_1 = y_2$. This forces $(x_1, y_1, z_1) = (x_2, y_2, z_2)$, so that $\text{dev}_{\text{pctf}}$ is injective. So $\text{dev}_{\text{pctf}}$ is a continuous injection, hence a local homeomorphism by invariance of domain.

We next claim that $\text{dev}_{\text{pctf}}$ is proper. So let $(x_n, y_n, z_n)$ leave all compact subsets of $\partial\Gamma^{(3)+}$. After taking a subsequence we may arrange for this sequence to have a limit $(x_\infty, y_\infty, z_\infty)$ in $\partial\Gamma^3$. Symmetry reduces our consideration to three cases, discussed below.

**Case 1.** $x_\infty = y_\infty = z_\infty$.

Here,

(3.2)
$$[(\xi^1(x_n) \cap \xi^1(z_n)) \oplus \xi^3(y_n)]^\perp = (\xi^{*3}(x_n) \cap \xi^{*3}(z_n)) \oplus \xi^{*1}(y_n),$$

where $(\xi^{*1}, \xi^{*2}, \xi^{*3})$ is the dual hyperconvex Frenet curve (defined in §2). The order condition of Lemma 2.5 is trivially satisfied and so
$$\lim_{n \to \infty} (\xi^{*3}(x_n) \oplus \xi^{*3}(z_n)) \cap \xi^{*1}(y_n) = \xi^{*3}(x_\infty).$$

Hence the expression in Eq. (3.2) converges to $\xi^1(x_\infty) \in \partial\Omega_\xi$.

**Case 2.** $x_\infty = y_\infty \neq z_\infty$.

Then $\lim_{n\to\infty} \xi^1(x_n) \oplus \xi^1(z_n) = \xi^1(x_\infty) \oplus \xi^1(z_\infty)$, which is transverse to $\xi^3(y_\infty) = \xi^3(x_\infty)$. So by Eq. (3.1),
$$\lim_{n\to\infty} (\xi^1(x_n) \oplus \xi^1(z_n)) \cap \xi^3(y_n) = (\xi^1(x_\infty) \oplus \xi^1(z_\infty)) \cap \xi^3(x_\infty) = \xi^1(x_\infty) \in \partial\Omega_\xi.$$

**Case 3.** $x_\infty = z_\infty \neq y_\infty$.

From hyperconvexity, $\lim_{n\to\infty} \xi^1(x_n) \oplus \xi^1(z_n) = \xi^2(x_\infty)$, which is transverse to $\xi^3(y_\infty)$. So
$$\lim_{n\to\infty} (\xi^1(x_n) \oplus \xi^1(z_n)) \cap \xi^3(y_n) = \xi^2(x_\infty) \cap \xi^3(y_\infty) \in \partial\Omega_\xi.$$

We conclude that dev$_{\text{pctf}}$ is an injective proper local homeomorphism of topological 3-manifolds and hence a homeomorphism onto $\Omega^1_\xi$. □

Our proof above shows that dev$_{\text{pctf}}$ maps leaves $g_{xz}$ of $\mathcal{G}$ to connected components of the twice-punctured projective lines
$$\ell_{xz} = \xi^1(x) \oplus \xi^1(z) - \{\xi^1(x), \xi^1(z)\}.$$
and that $\Omega^1_\rho$ is the disjoint union of these twice-punctured projective lines. So, if $\xi$ is induced by a Hitchin representation $\rho$ the twice-punctured projective lines $\ell_{xz}$ are the unions of two leaves of a $\rho(\Gamma)$-invariant foliation of $\Omega^1_\rho$ by properly embedded line segments. See Figure 9, Left. The leaves are all secant lines of $\xi^1(\partial\Gamma)$ in affine charts.

**Definition 3.8.** *The codimension-2 foliation of $\Omega^1_\xi$ specified by* dev$_{\text{pctf}}$ *is $\mathcal{G}_{\text{tcf}}$. The leaves of $\mathcal{G}_{\text{tcf}}$ are connected components of $\ell_{xz}$ ($(x,z) \in \partial\Gamma^{(2)}$).*

Furthermore, dev$_{\text{pctf}}$ maps leaves of $\mathcal{H}$ to leaves of $\mathcal{F}_{\text{pcf}}$. As before, in the presence of a group action (dev$_{\text{pctf}}, \rho$) defines a (PSL$_4(\mathbb{R}), \mathbb{RP}^3$) structure on T$^1 S$.

3.4. **Concave Transversely Foliated Projective Structures.** This is our final construction of developing maps. Define dev$_{\text{ctrf}} : \partial\Gamma^{(3)+} \to \Omega^2_\xi$ by
$$(x,y,z) \mapsto (\xi^1(x) \oplus (\xi^2(z) \cap \xi^3(x))) \cap \xi^3(y).$$
See Figure 8, Right. Note that the image of dev$_{\text{ctrf}}$ is disjoint from $\partial\Omega_\xi$ because $\xi^3(x) \cap \xi^3(y) \cap \xi^2(w) = \{0\}$ for all $(x,y,w) \in \partial\Gamma^{(3)+}$. The edge case to check here that dev$_{\text{ctrf}}(x,y,z) \notin \xi^2(w)$ for $w \in \{x,y\}$ follows from the structure of $\partial\Omega_\xi \cap \xi^3(x)$. Furthermore, the image of dev$_{\text{ctrf}}$ is contained in $\Omega^2_\xi$ since every image point is contained in at least two projective planes of the form $\xi^3(x)$ ($x \in \partial\Gamma$) by the definition of dev$_{\text{ctrf}}$.

**Proposition 3.9.** *The map* dev$_{\text{ctrf}}$ *is a 3-sheeted covering* $\partial\Gamma^{(3)+} \to \Omega^2_\rho$.

*Proof.* The strategy is a slight modification of the previous cases. We first show that restrictions of dev$_{\text{ctrf}}$ to leaves of $\mathcal{F}$ are homeomorphisms onto their images, then show dev$_{\text{ctrf}}$ is a local homeomorphism that exactly 3-to-1 surjects $\Omega^2_\rho$. The modification of initially restricting to leaves of $\mathcal{F}$ is made to handle the fact that preimages of points in $\Omega^2_\rho$ under dev$_{\text{ctrf}}$ do not seem to admit a simple and explicit description of a similar flavor to that of dev$_{\text{ctaf}}$.

Let $x \in \partial\Gamma$ be fixed and let $f_x$ be the leaf of $\mathcal{F}$ corresponding to $x$. Let $\phi_x : f_x \to \xi^3(x) \cap \Omega^2_\xi$ be the restriction of dev$_{\text{ctrf}}$. We show $\phi_x$ is a homeomorphism onto $\Omega^2_\rho \cap \xi^3(x)$ by verifying that it is a proper injective local homeomorphism between topological surfaces.

Let $(y, z)$ and $(r, w)$ be so $(x, y, z)$ and $(x, r, w)$ are in $\partial \Gamma^{(3)+}$. We wish to show $\phi_x(y, z) = \phi_x(r, w)$ if any only if $(y, z) = (r, w)$. We note $\text{dev}_{\text{ctrf}}(x, y, z) = (\xi^1_{x,3}(x) \oplus \xi^1_{x,3}(z)) \cap \xi^2_{x,3}(y)$ in the notation of Lemma 3.1. Note that $\phi_x(y, z)$ is contained in $\xi^1_{x,3}(x) \oplus \xi^1_{x,3}(z) - \{\xi^1(x)\}$ and $\phi_x(r, w)$ is contained in $\xi^1_{x,3}(x) \oplus \xi^1_{x,3}(w) - \{\xi^1(x)\}$. If $w \neq z$ then these once-punctured projective lines are disjoint and so $\phi_x(y, z) \neq \phi_x(r, w)$. The case $z = w$ is handled by noting that for distinct $y, r \in (x, w)$ that the dual formulation of Lemma 2.5 in terms of intersections implies

$$\xi^2_{x,3}(y) \cap \xi^2_{x,3}(r) \cap (\xi^1_{x,3}(x) \oplus \xi^1_{x,3}(z)) = \{0\}.$$

So $\phi_x(y, z) \neq \phi_x(r, w)$ in this case. This implies $\phi_x$ is injective and continuous, hence a local homeomorphism by invariance of domain.

We next claim $\phi_x$ is proper. Let $(y_n, z_n)$ leave all compact subsets of $f_x$. After subsequence, $(y_n, z_n)$ converges to a point $(y_\infty, z_\infty) \in \partial \Gamma^2$. We break into four cases. In all cases, only three points occur and so the order hypotheses required to take limits are trivial.

**Case 1.** $y_\infty = z_\infty = x$.

Then $\lim_{n \to \infty} \phi_x(y_n, z_n) = \lim_{n \to \infty} (\xi^1_{x,3}(x) \oplus \xi^1_{x,3}(z_n)) \cap \xi^2_{x,3}(y_n) = \xi^1_{x,3}(x) \in \partial \Omega_\xi$.

**Case 2.** $z_\infty = x$ and $y_\infty \neq z_\infty$.

Then $\lim_{n \to \infty} \xi^1(x) \oplus \xi^1_{x,3}(z_n) = \xi^2(x)$, which is transverse to $\xi^2_{x,3}(y_\infty)$, so that

$$\lim_{n \to \infty} (\xi^1(x) \oplus \xi^1_{x,3}(z_n)) \cap \xi^2_{x,3}(y_n) = \xi^2(x) \cap \xi^2_{x,3}(y_\infty) \in \partial \Omega_\xi.$$

**Case 3.** $y_\infty = x$ and $z_\infty \neq x$.

Then $(\xi^1(x) \oplus \xi^1_{x,3}(z_\infty)) \cap \xi^2(x) = \xi^1(x)$ so that

$$\lim_{n \to \infty} \phi(y_n, z_n) = \lim_{n \to \infty} (\xi^1_{x,3}(x) \oplus \xi^1_{x,3}(z_n)) \cap \xi^2_{x,3}(y_n) = \xi^1(x) \in \partial \Omega_\xi.$$

**Case 4.** $y_\infty = z_\infty$ and $z_\infty \neq x$.

Here, $\xi^1_{x,3}(x) \oplus \xi^1_{x,3}(z_\infty)$ is transverse to $\xi^2_{x,3}(y_\infty)$ and so

$$\lim_{n \to \infty} \phi_x(y_n, z_n) = \lim_{n \to \infty} (\xi^1_{x,3}(x) \oplus \xi^1_{x,3}z_n)) \cap \xi^2_{x,3}(y_n)$$
$$= (\xi^1_{x,3}(x) \oplus \xi^1_{x,3}(z_\infty)) \cap \xi^2_{x,3}(y_\infty) = \xi^1_{x,3}(z_\infty) \in \partial \Omega_\xi.$$

So $\phi_x$ is proper, hence a homeomorphism. We next observe that any $p \in \Omega^2_\xi$ has exactly three preimages under $\text{dev}_{\text{ctrf}}$. This is because every $p \in \Omega^2_\xi$ is contained in $\xi^3(x)$ for exactly three $x \in \partial \Gamma$, and for all $x \in \partial \Gamma$, the restricted developing map $\phi_x$ is a homeomorphism of $f_x$ onto $\xi^3(x) \cap \Omega^2_\xi$.

We next claim that $\text{dev}_{\text{ctrf}}$ is locally injective. Then by invariance of domain it follows that $\text{dev}_{\text{ctrf}}$ is a local homeomorphism.

The observation here is that for $x'$ near $x$ the line $\xi^3(x') \cap \xi^3(x)$ is close to $\xi^2(x)$, and hence not near $\text{dev}_{\text{ctrf}}(x, y, z)$ for $y$ and $z$ less close to $x$. To make this precise, note the map

$$\Psi : \partial \Gamma^2 \to \text{Gr}_2(\mathbb{R}^4)$$
$$(x, y) \mapsto \begin{cases} \xi^3(x) \cap \xi^3(y) & x \neq y \\ \xi^2(x) & x = y \end{cases}$$

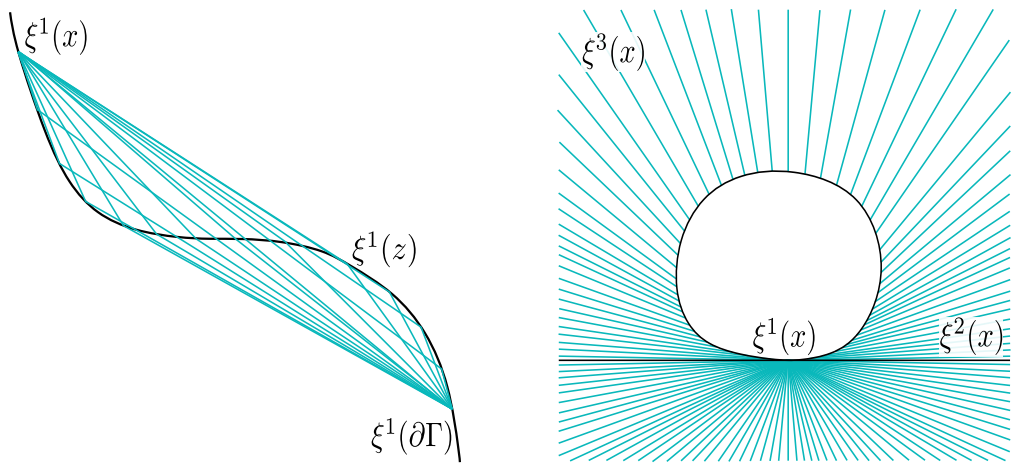

FIGURE 9. Leaves of $\mathcal{G}_{\text{tcf}}$ (left) and $\mathcal{G}_{\text{ctrf}}$ (right). In the left illustration, the black curve is a segment of $\xi^1(\partial\Gamma)$. The right illustration is contained within a subspace of the form $\xi^3(x)$ for some $x \in \partial\Gamma$. In the right illustration, points in this plane are contained in two other leaves of $\mathcal{G}_{\text{ctrf}}$ that do not lie in this plane.

is continuous, and the image of the diagonal of $\partial\Gamma^2$ under $\Psi$ is $\partial\Omega_\rho$. So for every neighborhood $U$ of $\partial\Omega_\rho$ in $\mathbb{RP}^3$ there is a neighborhood $V$ of the diagonal in $\partial\Gamma^2$ so that if $(x, y) \in V$ then $\Psi(x, y) \subset U$.

Now, let $p = (x, y, z) \in \partial\Gamma^{(3)+}$. Take a neighborhood $U$ of $\partial\Omega_\rho$ with closure disjoint from $\text{dev}_{\text{ctrf}}(p)$ and let $V$ be the neighborhood of the diagonal in $\partial\Gamma^2$ so that $\Psi(x,') \subset U$ for all $(x, x') \in V$. Now suppose $q = (x', y', z') \in \partial\Gamma^{(3)+}$ has $(x, x') \in V$. If $x \neq x'$ then $\text{dev}_{\text{ctaf}}(q) \cap \text{dev}_{\text{ctaf}}(p)$ is contained in $\xi^3(x) \cap \xi^3(x')$, which is contained in $U$. So $\text{dev}_{\text{ctaf}}(q) \neq \text{dev}_{\text{ctaf}}(p)$ by construction of $U$. Combining this with the injectivity of $\phi_x$ in the case $x = x'$ shows that $\text{dev}_{\text{ctaf}}$ is locally injective.

It is a standard fact from elementary topology that a surjective local homeomorphism $f : X \to Y$ between Hausdorff spaces with the property that there is a $k \in \mathbb{N}$ so for every $p \in Y$, the preimage of $p$ under $\phi$ has cardinality $k$ is a $k$-sheeted covering. So $\text{dev}_{\text{ctrf}}$ is a 3-sheeted covering. □

Our proof showed that $\text{dev}_{\text{ctrf}}(g_{xz})$ is a properly embedded properly convex subset of a projective line for every leaf $g_{xz}$ of $\mathcal{G}$ and that $\text{dev}_{\text{ctrf}}(f_x)$ is a leaf of $\mathcal{F}_{\text{ccf}}$ for any leaf $f_x \in \mathcal{F}$. As before, $\text{dev}_{\text{ctrf}}$ specifies a trifoliation $\mathcal{G}_{\text{ctrf}}$ with leaves $\text{dev}_{\text{ctrf}}(g_{xz})$ for $g_{xz}$ a leaf of $\mathcal{G}$. See Figure 9, Right. If $\xi$ is induced by a Hitchin representation $\rho$, then $\mathcal{G}_{\text{ctrf}}$ is $\rho(\Gamma)$-invariant. Additionally, from $\rho$-equivariance $(\text{dev}_{\text{ctrf}}, \rho)$ gives a $(\text{PSL}_4(\mathbb{R}), \mathbb{RP}^3)$-structure on $\text{T}^1 S$.

## 4. THE SHAPE OF THE DOMAIN

Our proofs of our main theorems rely on the precise qualitative behavior of the shape of $\partial\Omega_\xi$. This section establishes these matters. We begin by showing that $\partial\Omega_\xi$ is locally convex except on the image of $\xi^1$. Then we classify intersections of projective planes with $\partial\Omega_\xi$, describe the view an observer on the boundary of $\partial\Omega_\xi$ outside of the image of $\xi^1$ has of $\partial\Omega_\rho$, and describe some projections of $\xi^1(\partial\Gamma)$ to projective planes of the form $\xi^3(x)$ ($x \in \partial\Gamma$).

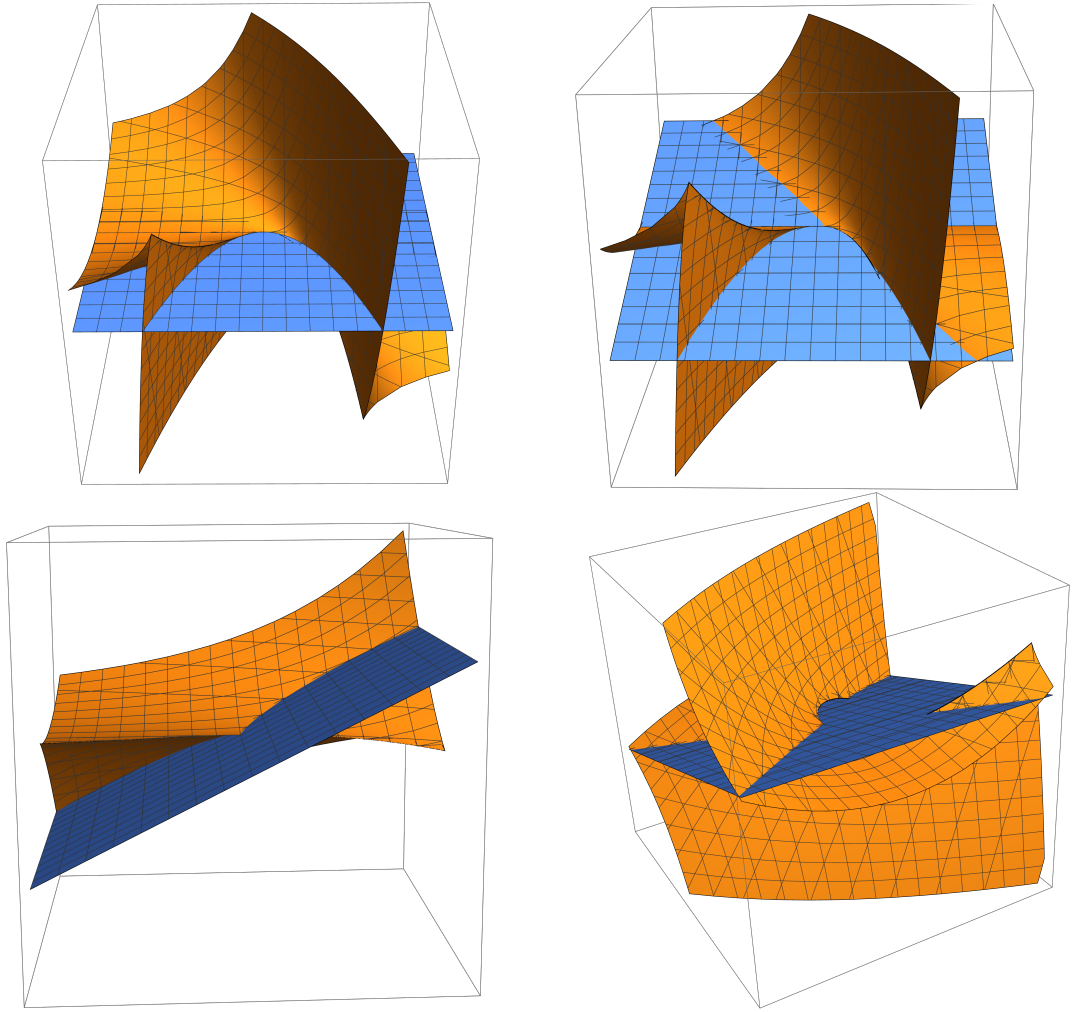

FIGURE 10. The four types of plane intersections with the Fuchsian domain, with views chosen to highlight relevant features. Top, left to right: planes of the form $\xi^3(x)$ and $\xi^2(x) \oplus \xi^1(y)$. Bottom, left to right: planes of the form $(\xi^3(x) \cap \xi^3(z)) \oplus \xi^1(y)$ and $\xi^1(x) \oplus \xi^1(y) \oplus \xi^1(z)$.

4.0.1. *Remarks on Figures.* We have found accurate pictures indispensable in developing intuition for the results proved in this section. We describe how to produce such figures.

All computer renderings in this paper are of the *Fuchsian* example. That is, this domain of discontinuity is preserved by a $\mathrm{PSL}_2(\mathbb{R})$-subgroup of $\mathrm{PSL}_4(\mathbb{R})$ that acts irreducibly on $\mathbb{R}^4$. All domains of discontinuity for $\mathrm{PSL}_4(\mathbb{R})$ Hitchin representations contained in such irreducible $\mathrm{PSL}_2(\mathbb{R})$ subgroups are projectively equivalent to this example. This example has exceptional smoothness behaviors that are never true in other examples (e.g. [23, 30]). One may compute using discriminants and the description of the Veronese embedding in §4 of [11] that the intersection of $\partial \Omega_\xi$ with an appropriate affine chart is the zero locus of $F(x, y, z) = 18xyz - 4x^3z + x^2y^2 - 4y^3 - 27z^2$.

The hand illustrations throughout the paper emphasize convexity properties less clearly visible in computer renderings. A serious attempt has been made to make the illustrations consistent with the many qualitative constraints placed by the results proved below.

4.1. **Local Convexity.** The following decomposition of $\partial\Omega_\xi$ tends to be useful, as the two pieces behave quite differently (see also [24, §3]).

**Definition 4.1.** *Let $\Xi^1 = \text{Im}(\xi^1)$ and $\Xi^2 = \partial\Omega_\xi - \Xi^1$.*

Away from $\Xi^1$, the boundary $\partial\Omega_\xi$ is locally convex and regular:

**Lemma 4.2** (Locally Convex Off Frenet). *Let $p \in \Xi^2$. Then $\partial\Omega_\xi$ is $C^1$ in a neighborhood $V$ of $p$. We may take $V$ so that there is a properly convex domain $\Psi_V$ so $V \subset \partial\Psi_V$ and $V \subset \partial(\Psi_V \cap \Omega_\xi^1)$.*

**Remark 4.3.** *A different and more general argument for $C^1$ regularity, inspired by knowledge of this case, appears in other work of the author ([24, §3.2]). In contrast to the argument there, the present proof also gives the conclusion of local convexity, which is essential in the sequel. The conclusion of $C^1$ regularity cannot be improved beyond $C^{1+\alpha}$ regularity (which can be seen to follow from the argument in [24]).*

*Proof.* We first show that $\overline{\Omega_\xi^1}$ is locally cut out as the intersection of half-spaces except at points in $\Xi^1$, which gives the desired local convexity. Fix a reference metric $d$ on $\mathbb{RP}^3$. Since $p \in \Xi^2$, there are $(x_0, y_0) \in \partial\Gamma^{(2)}$ so that $p = \xi^2(x_0) \cap \xi^3(y_0) = \xi^1_{y_0}(x_0)$. The point of the argument is the observation, from the qualitative structure of $\xi^3(x_0) \cap \partial\Omega_\rho$, that close-by lines $\xi^2(z)$ for $z \in \partial\Gamma$ near $x_0$ have uniform lower bounds on the distance between $p = \xi^2(x_0) \cap \xi^3(y_0)$ and their points of intersection $\xi^2(z) \cap \xi^3(x_0)$ with $\xi^3(x_0)$. Note the interchange of the position of $x_0$ between the two expressions.

First, we restrict neighborhoods appropriately. From the continuity of $\xi^1_{x'}(w)$ and the fact that $\lim_{x \to x'} \xi^1_{x'}(x) = \xi^1(x)$,

(1) There is a neighborhood $U_2$ of $y_0$ and a neighborhood $U_1$ of $x_0$ and a positive $\delta > 0$ so that $d(\xi^1_a(b), \xi^1_c(d)) > \delta$ for all $a \in U_2$ and $b, c, d \in U_1$. Take $U_1$ and $U_2$ to be connected.

(2) After reducing the size of the connected neighborhoods $U_1$ and $U_2$, there are neighborhoods $U_3$ of $x_0$ and $U_4$ of $y_0$ so that with $V_1 = U_4 \times U_3$, the image $V$ of $V_1$ under the two-argument map $\xi^1(\cdot)$ is entirely contained in an affine chart $\mathcal{A}$, which is fixed throughout the following.

To record the local one-sidedness used below, pick a small line segment $\ell_{y_0}$ in the closure of the convex domain $C_{y_0} \subset \xi^3(y_0)$ so that $\ell_{y_0}$ is entirely contained in $\mathcal{A}$, one endpoint of $\ell_{y_0}$ is on $\partial C_{y_0}$, and the other endpoint of $\ell_{y_0}$ is in $C_{y_0}$.

Now let $z \in U_3$. We consider the function $d_z$ on $\mathcal{A}$ defined as signed distance from $\xi^3(z)$, with sign determined by being positive on the open half-space $\Pi_z$ containing $\ell_{y_0}$. For any $w \in U_3 - \{z\}$, our choice of neighborhoods implies that the line segment $\xi^1_{U_4}(w)$ does not intersect $\xi^3(z)$, and so $d_z$ restricted to $\xi^1_{U_4}(w)$ has constant sign. Furthermore, $d_z$ restricted to the segment of $\partial C_{y_0}$ given by $\xi^1_{y_0}(U_3)$ is non-negative and zero only at $\xi^1_{y_0}(z)$ by the strict convexity of $C_{y_0}$.

We conclude that for all $r \in U_4$ and $z \in U_3$, the function $d_z$ evaluated on $\xi^1_r(U_3)$ is non-negative and zero only at $\xi^1_r(z)$. From continuity and non-negativity of $d_{y_0}$ we conclude that every $d_z$ is non-negative on $\partial C_z \cap V_2$.

A point $p \in \mathcal{A}$ is contained in the properly convex set $\Psi_V = \bigcap_{z \in U_3} \Pi_z$ if and only if $d_z(p) > 0$ for all $z \in U_3$. Similarly, $p$ is contained in $\partial\Psi_V$ if and only if $d_z(p) \geq 0$ for all $z \in U_3$ and there is a $z \in U_3$ so $d_z(p) = 0$. It follows that $V$ is contained in $\partial\Psi_V$, and all claims except that $V$ is $C^1$ have been established.

Now, as $V$ is contained in the boundary of a properly convex domain, to show that $V$ is $C^1$ it suffices to prove that every point in $V$ has a unique supporting plane. We show that $\xi^3(b)$ is the unique supporting plane to $V$ at $\xi_a^1(b)$ for all $a \in U_4$ and $b \in U_3$. Indeed, the above implies $\xi^3(b)$ is a supporting plane at any such point. On the other hand, because a segment of $\xi^2(b)$ containing $\xi_a^1(b)$ is contained in $\partial V$, if $P$ is any supporting plane to $V$ at $\xi_a^1(b)$ then $P$ must contain $\xi^2(b)$. Furthermore, $P$ must support $V \cap \xi^3(a)$ at $\xi_a^1(b)$, and hence contain $\xi^3(a) \cap \xi^3(b)$ since $\partial C_a$ is $C^1$ at $\xi_a^1(b)$ (by [23, Prop. 3.7] or the Frenet Restriction Lemma). This forces $P = \xi^3(b)$, and concludes the proof. □

4.2. **Sum-only Plane Intersections.** We now turn to giving detailed descriptions of intersections of $\partial \Omega_\xi$ with projective planes, and descriptions of views of $\partial \Omega_\xi$ from appropriate vantage points. The motivation for our choice of planes to intersect $\partial \Omega_\xi$ with is that propositions from the previous section applied to the dual hyperconvex Frenet curve imply:

**Lemma 4.4** (Tetrachotomy). *Every projective plane in $\mathbb{RP}^3$ has exactly one of the four mutually exclusive forms $\xi^3(x)$ or $\xi^2(x) \oplus \xi^1(y)$ or $(\xi^3(x) \cap \xi^3(z)) \oplus \xi^1(y)$ or $\xi^1(x) \oplus \xi^1(y) \oplus \xi^1(z)$ for $x, y, z$ distinct points in $\partial \Gamma$.*

*Proof.* Indeed, from §3, the maps $\partial \Gamma^{(3)+} \to \mathbb{RP}^3$

$$(x, y, z) \mapsto (\xi^1(x) \oplus \xi^1(z)) \cap \xi^3(y),$$
$$(x, y, z) \mapsto \xi^3(x) \cap \xi^3(y) \cap \xi^3(z)$$

surject $\mathbb{RP}^3 - \partial \Omega_\xi$. As every point in $\partial \Omega_\xi$ has the form of either $\xi^2(x) \cap \xi^3(y)$ ($x \neq y$) or $\xi^1(x)$, considering duals shows that every projective plane in $\mathbb{RP}^3$ has exactly one of the four desired forms. □

The simplest case of the tetrachotomy of Lemma 4.4 is $\xi^3(x)$, which is handled by Guichard-Wienhard's work [11] and [23, Prop. 3.7] or the Frenet Restriction Lemma 3.1.

This subsection handles the cases of planes that are formed using only sums of elements of $\xi$, namely $\xi^2(x) \oplus \xi^1(y)$ and $\xi^1(x) \oplus \xi^1(y) \oplus \xi^1(z)$. Analysis of the final case, $(\xi^3(x) \cap \xi^3(z)) \oplus \xi^1(y)$ is made considerably more complicated by the essential mixing of sum and intersection in the subspace's description, which makes Lemma 2.5 not apply at a key step. We postpone this case until we have developed the necessary machinery, which relies on the analysis of the case $\xi^2(x) \oplus \xi^1(y)$.

Let us adopt an elementary definition that is relevant in the sequel. See Figure 11, Right for an example of the objects the below formalizes.

**Definition 4.5.** *A continuous injection $\gamma : [-\epsilon, \epsilon] \to \mathbb{RP}^2$ has a* balanced cusp *at $\gamma(0)$ if:*

*(1) The restrictions of $\gamma$ to $(-\epsilon, 0)$ and $(0, \epsilon)$ are $C^1$,*
*(2) $\gamma([-\epsilon, \epsilon])$ has one-sided tangent lines at $\gamma(0)$, which agree,*
*(3) $\gamma$ crosses its one-sided tangent line $L$ at $\gamma(0)$.*
*(4) If $L'$ is any line transverse to $L$ through $\gamma(0)$, then $\gamma$ does not cross $L'$ at $\gamma(0)$.*

4.2.1. *Intersections with $\xi^2(x) \oplus \xi^1(y)$.* This is the most important new intersection pattern for our later arguments. We will prove that Figure 11, Left is qualitatively accurate. More precisely, we have:

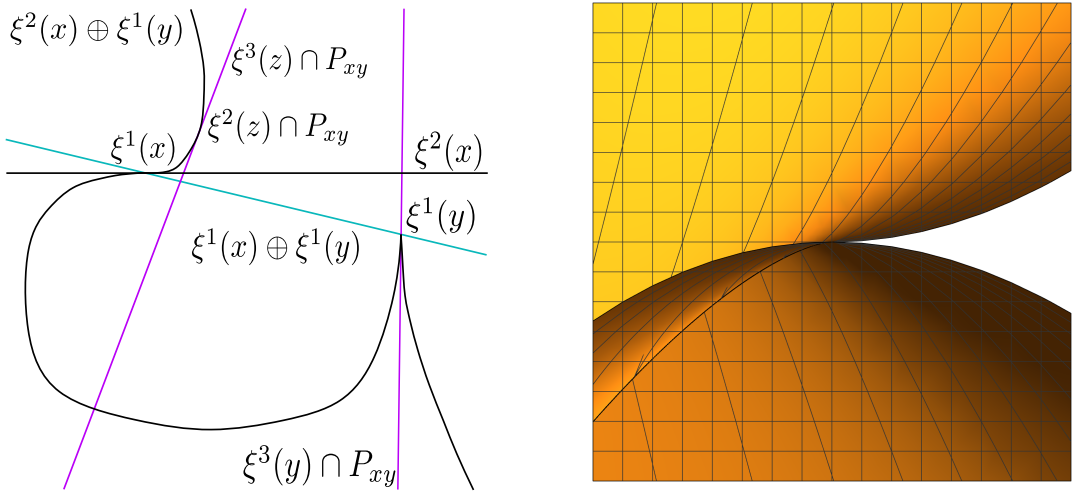

FIGURE 11. Left: $(\xi^2(x)\oplus\xi^1(y))\cap\partial\Omega_\xi$ in black. Intersections with subspaces of the form $\xi^3(z)$ are in purple and the line $\ell_{xy}$ is in teal. Right: a view of the Fuchsian domain in $\mathbb{RP}^3$ where the presence of balanced cusps is clear.

**Proposition 4.6** ($\xi^2(x)\oplus\xi^1(y)$ Pattern)**.** *Let $(x,y)\in\partial\Gamma^{(2)}$. Introduce the notation $P_{xy} = \xi^2(x)\oplus\xi^1(y)$ and $C_{xy}=\partial\Omega_\xi\cap P_{xy}$. Let $\eta_{xy}=(\eta^1_{xy},\eta^2_{xy}):\partial\Gamma\to\mathcal{F}(P_{xy})$ be given by*

$$(4.1) \qquad \eta^1_{xy}(z) = \begin{cases} \xi^2(z)\cap P_{xy} & z\neq x,y \\ \xi^1(z) & z=x \text{ or } z=y \end{cases}, \qquad \eta^2_{xy}(z) = \xi^3(z)\cap P_{xy}.$$

*Then:*
  *(1) $C_{xy} = \xi^2(x) \cup \eta^1_{xy}(\partial\Gamma)$,*
  *(2) The map $\eta_{xy}$ and its image satisfy:*
    *(a) The restriction of $\eta_{xy}$ to $[y,x]$ or $[x,y]$ is a hyperconvex Frenet curve,*
    *(b) Both $\eta^1_{xy}$ and $\eta^2_{xy}$ are injective and the full map $\eta_{xy}$ is continuous.*
    *(c) The image under $\eta^1_{xy}$ of any open set containing $x$ or $y$ is not convex.*
    *(d) $\eta^1(\partial\Gamma)$ is $C^1$ at $\xi^1(x)$ and crosses its tangent line $\xi^2(x)$ at $\xi^1(x)$,*
    *(e) $\eta^1(\partial\Gamma)$ has a balanced cusp at $\xi^1(y)$,*
    *(f) For every $z\in\partial\Gamma-\{y,x\}$, let $\mathcal{U}_z$ be the connected component of $\partial\Gamma-\{y,x\}$ not containing $z$. Then $\eta^2_{xy}(z)$ intersects $\eta^1_{xy}(\mathcal{U}_z)$ exactly once.*
  *(3) The relative positions of $\Omega^1_\xi$ and $\Omega^2_\xi$ in $P_{xy}$ satisfy:*
    *(a) All cusps in $\partial\Omega_\xi\cap P_{xy}$ have their non-convex sides in $\Omega^2_\xi$.*
    *(b) $\Omega^1_\xi$ is on the convex sides of $\eta^1_{xy}([x,y])$ and $\eta^1_{xy}([y,x])$.*

We remark that Claim (2f) is useful in establishing qualitative constraints on various constructions below, including the "top view" projections of $\xi^1(\partial\Gamma)$ in §4.4.1.

*Proof.* Claim (1) is immediate from transversality of $\xi^2(z)$ ($z\neq x$) to $P_{xy}$, which follows from hyperconvexity of $\xi$ and the definition of $\partial\Omega_\xi$.

We next prove that the restrictions of $\eta_{xy}$ to the half-open intervals $[y,x)$ and $(x,y]$ are hyperconvex Frenet curves. First note that transversality of relevant line and plane fields to $P_{xy}$ implies $\eta_{xy}$ is continuous on $[y,x)$. We verify hyperconvexity through its formulation in terms of intersections (Lemma 2.4).

So,[2] let $p_1, ..., p_k$ be so $p_1 + ... + p_k = 3$ and let $z_1, ..., z_k$ be in $(y, x)$. Order $z_1, ..., z_k$ so that $x < y < z_1 < ... < z_k$. Then by Lemma 2.5,

$$\left(\bigcap_{i=1}^{k} \eta_{xy}^{3-p_i}(z_i)\right)^{\perp} \supset \left[(\xi^2(x) \oplus \xi^1(y)) \cap \left(\bigcap_{i=1}^{k} \xi^{4-p_i}(z_i)\right)\right]^{\perp}$$

$$= (\xi^{*2}(x) \cap \xi^{*3}(y)) + \left(\bigoplus_{i=1}^{k} \xi^{*p_i}(z_i)\right)$$

$$= \mathbb{R}^4,$$

so that the general position requirement is satisfied on $(y, x)$.

To see limit compatibility on $(y, x)$, let $p_1, ..., p_k$ have $p = \sum_{i=1}^{k} p_i \leq 3$ and let $(x_1^n, ..., x_k^n) \in (y, x)$ be strictly increasing in the cyclic order on $\partial\Gamma$ for every fixed $n$ and $\lim_{n \to \infty}(x_1^n, ..., x_k^n) = (z, ..., z)$ for some $z \in (y, x)$. Then using hyperconvexity as proved above,

$$\lim_{n \to \infty} \left(\bigcap_{i=1}^{k} \eta_{xy}^{3-p_i}(x_i^n)\right)^{\perp} = \lim_{n \to \infty} \left((\xi^2(x) \oplus \xi^1(y)) \cap \left(\bigcap_{i=1}^{k} \xi^{4-p_i}(x_i^n)\right)\right)^{\perp}$$

$$= \lim_{n \to \infty} (\xi^{*2}(x) \cap \xi^{*3}(y)) \oplus \left(\bigoplus_{i=1}^{k} \xi^{*p_i}(x_i^n)\right)$$

$$= (\xi^{*2}(x) \cap \xi^{*3}(y)) \oplus \xi^{*p}(z)$$

$$= (\eta_{xy}^{3-p}(z))^{\perp},$$

so that that limit compatibility holds on $(y, x)$. So the restriction of $\eta_{xy}$ to $(y, x)$ is a hyperconvex Frenet curve. The restriction to $(x, y)$ is analogous. Observe also that $\eta_{xy}^2$ is continuous on all of $\partial\Gamma$ as a consequence of transversality of all $\xi^3(z)$ to $P_{xy}$ for $z \in \partial\Gamma$.

We next claim that $\eta_{xy}^1$ is continuous on $[y, x]$. From this, it is an exercise in plane geometry and calculus to see that the restriction of $\eta_{xy}$ to all of $[y, x]$ is a hyperconvex Frenet curve. Continuity at $y$ is an immediate consequence of the observation that $\xi^2(y)$ is transverse to $P_{xy}$ and that $\xi^2(y) \cap P_{xy} = \xi^1(y)$. The continuity of $\eta_{xy}^1$ at $x$ requires some more care.

So let $x_n \in [y, x)$ be a sequence with limit $x$. We claim $\eta_{xy}^1(x_n)$ has limit $\xi^1(x)$. After taking a subsequence, we may at least arrange for $\eta_{xy}^1(x_n)$ to converge to some point $p_\infty$. Fixing a background metric on $\mathbb{RP}^3$, since $\xi^2$ is continuous as a map to $\text{Gr}_2(\mathbb{R}^4)$, the projective lines $\xi^2(x_n)$ converge to $\xi^2(x)$ in the Hausdorff topology. In particular, any limit point of $\xi^2(x_n) \cap P_{xy}$ must be in $\xi^2(x)$. So $p_\infty \in \xi^2(x)$.

The basic observation is now that hyperconvexity of $\xi$ shows $\lim_{n \to \infty} \xi^3(x_n) \cap \xi^2(x) = \xi^1(x)$, and after further subsequence the convergence is monotone. The considerations that:

(1) $\eta_{xy}^1$ restricted to $(y, x)$ traces out a properly convex $C^1$ arc in $P_{xy}$ with tangent lines given by $\eta_{xy}^2$,
(2) $\xi^2(z)$ never intersects $\xi^1(x) \oplus \xi^1(y)$ for $z \notin \{x, y\}$

trap $p_\infty$ in the closed projective line segment

$$[x_1, x] := (\xi^3([x_1, x)) \cap \xi^2(x)) \cup \{\xi^1(x)\}.$$

---

[2]The following notation appears general, but of course only handles a small number of actual cases. We shall adopt such general-seeming notations in the following when it avoids redundant case analysis.

By the same reasoning, for all $n$ because the projective line $\xi^3(x_n) \cap P_{xy}$ is a supporting line to $\eta^1_{xy}([y,x))$ at $\eta^1_{xy}(x_n)$, the point $p_\infty$ is trapped in the analogously defined segment $[x_n, x]$. The Frenet Restriction Lemma 3.1 implies the intersection of these segments is exactly $\{\xi^1(x)\}$. We conclude that $p_\infty = \xi^1(x)$. The independence of the limit point $\xi^1(x)$ of the choices of subsequence in the argument implies the desired continuity of $\eta^1_{xy}$ restricted to $[y, x]$ at $x$.

All analogous statements for $[x, y]$ are similar. This proves (2a) and the continuity claim in (2b). The injectivity of $\eta^1_{xy}$ is an immediate consequence of hyperconvexity. The injectivity of $\eta^2_{xy}$ follows from Prop. 3.7. Indeed, Prop. 3.7 implies that the map $(\partial\Gamma - \{x, y\}) \to \ell_{xy}$ given by $z \mapsto \xi^3(z) \cap \ell_{xy}$ is a homeomorphism, from which injectivity of $\eta^2_{xy}$ as a map to $\mathrm{Gr}_2(P_{xy})$ follows.

We now examine the local structure of $C_{xy}$ at $\xi^1(x)$ and $\xi^1(y)$. At each such point, $\mathrm{Im}(\eta^1_{xy})$ is locally obtained by gluing together two strictly convex, properly convex $C^1$ arcs that are disjoint except at their endpoints and which have agreeing tangent lines at their endpoints. There are four possible such configurations, up to renormalization, as depicted in Figure 12.

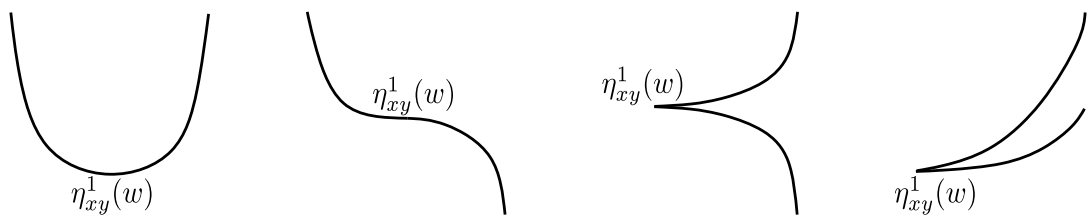

FIGURE 12. Possible configurations of the image of $\eta^1_{xy}$ at points in $\Xi^1$. These are depicted around a point of the form $\xi^1(w)$ for $w \in \{x, y\}$. In the below, we refer to these configurations as type (1)-(4) in accordance with their arrangement from left to right above.

We first address $\xi^1(y)$. The behavior of $\mathrm{dev}_{\mathrm{pcf}}$ recorded in §3.1 and the qualitative structure of $\xi^3(y) \cap \partial\Omega_\xi$ implies that the common tangent line $\xi^3(y) \cap P_{xy}$ intersects $\partial\Omega_\xi$ exactly at the two points $\xi^1(y)$ and $\xi^3(y) \cap \xi^2(x)$, and contains nonempty segments in both $\Omega^1_\xi$ and $\Omega^2_\xi$ in any neighborhood of either point. Arrangements (1) and (4) are inconsistent with this because no transition of connected components of $\Omega_\xi$ occurs on the tanent line at the center point. Counting intersections of tangent lines and applying Lemma 3.6 shows that if arrangement (2) occcurs then each local side of $P_{xy} - C_{xy}$ would need to be contained in $\Omega^2_\xi$, which is impossible by the same considerations on the tangent line $\xi^3(y) \cap P_{xy}$. So configuration (3) occurs at $\xi^1(y)$.

We next consider $\xi^1(x)$. Since $P_{xy}$ is transverse to $\xi^3(x)$, which is the tangent plane to $\partial\Omega_\xi$ at $\xi^2(x) \cap \xi^3(z)$ for any $z \neq x$, the component of $\Omega_\xi$ that $P_{xy}$ intersects changes as $\xi^2(x)$ is crossed at any point other than $\xi^1(x)$. Arrangement (4) is ruled out by applying Lemma 3.6 to the two segments of $\eta^1_{xy}(\partial\Gamma)$ meeting at $\xi^1(x)$, as it would imply a transverse crossing of $\partial\Omega_\xi$, both sides of which intersect the same connected component $\Omega^2_\xi$. Arrangements (1) and (3) are similarly ruled out by our remark on crossings at $\xi^2(x)$ and Lemma 3.6.

This proves (2c)-(2d). Both claims in (3) follow from applications of Lemma 3.6.

To prove the final claim (2f), let $z \in (x, y)$. Note that the union of $\eta^1_{xy}([y, x])$ and a connected component $\mathcal{C}$ of $\ell_{xy} = \xi^1(x) \oplus \xi^1(y) - \{\xi^1(x), \xi^1(y)\}$ is the boundary of a properly convex domain $D_{[y,x]}$ in $P_{xy}$. There is an analogous domain $D_{[x,y]}$ for $[x, y]$. By hyperconvexity of $\xi$, the line $\eta^2_{xy}(z) = \xi^3(z) \cap P_{xy}$ must intersect $\xi^1(x) \oplus \xi^1(y)$ in a unique point $p_z \in \ell_{xy}$. By convexity of $D_{[x,y]}$, we have $p_z \in \partial D_{[y,x]}$; by the convexity of $D_{[y,x]}$, there is exactly one other point of intersection of $\eta^2_{xy}(z)$ with $\partial D_{[y,x]}$, which must be in $\eta^1_{xy}((y,x))$. □

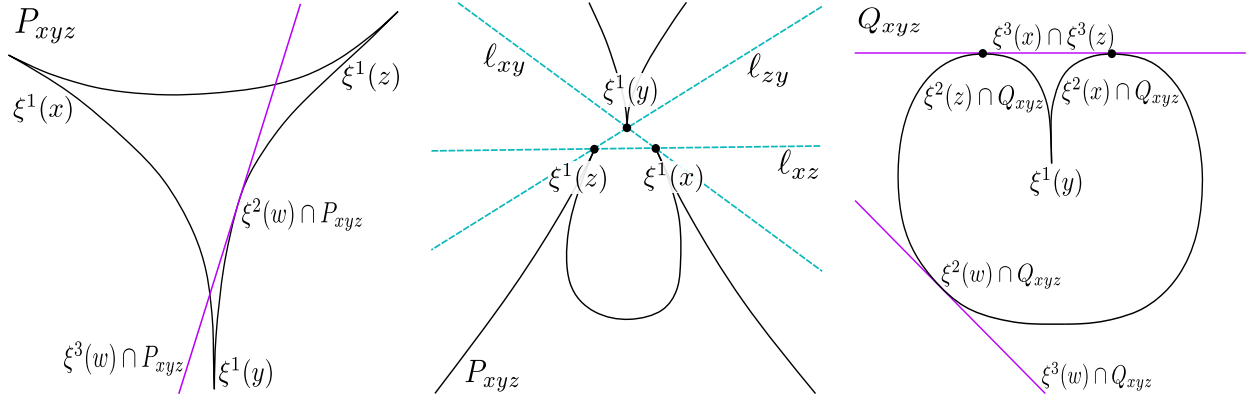

FIGURE 13. Left: annotated intersection of $\xi^1(x) \oplus \xi^1(y) \oplus \xi^1(z)$ with $\partial \Omega_\xi$ in an affine chart that emphasizes symmetric roles of entries. Middle: sketch of the same intersection in an affine chart that emphasizes behavior under degenerations. Here, $(x, y, z)$ is a positively oriented triple of close-together points. Right: annotated intersection of $(\xi^3(x) \cap \xi^3(z)) \oplus \xi^1(y)$ with $\partial \Omega_\xi$.

Note that the order hypothesis in Guichard's work [9] is crucial in applications of Lemma 2.5 in the preceding proof. It plays similarly essential roles in the rest of the section.

4.2.2. *Intersections with $\xi^1(x) \oplus \xi^1(y) \oplus \xi^1(z)$.* Let us call $P_{xyz} = \xi^1(x) \oplus \xi^1(y) \oplus \xi^1(z)$ and $C_{xyz} = \partial \Omega_\xi \cap P_{xyz}$ for $(x, y, z) \in \partial \Gamma^{(3)+}$.

Analysis of this case is simplified by the transversality of $P_{xyz}$ to all projective lines and planes of the forms $\xi^2(x)$ and $\xi^3(x)$, and the symmetry of the roles of $x, y$, and $z$. The picture of $C_{xyz} \subset P_{xyz}$ is given in two different affine charts in Figure 13, Left and Middle.

**Proposition 4.7** ($\xi^1(x) \oplus \xi^1(y) \oplus \xi^1(z)$ Pattern). *Let $x, y, z$ be distinct points in $\partial \Gamma$ and $\eta_{xyz} = (\eta^1_{xyz}, \eta^2_{xyz}) : \partial \Gamma \to \mathcal{F}(P_{xyz})$ be given by*

(4.2) $$\eta^1_{xyz}(w) = \xi^2(w) \cap P_{xyz}, \qquad \eta^2_{xyz}(w) = \xi^3(w) \cap P_{xyz}.$$

*Then:*
  (1) $C_{xyz} = \eta^1_{xyz}(\partial \Gamma)$,
  (2) *The map $\eta_{xyz}$ and its image satisfy:*
    (a) *The restriction of $\eta_{xyz}$ to $[x, y]$, $[y, z]$, or $[z, x]$ is a hyperconvex Frenet curve,*
    (b) *$\eta^1_{xyz}$ and $\eta^2_{xyz}$ are continuous injections,*
    (c) *The image of $\eta^1_{xyz}$ has balanced cusps at $\xi^1(x)$, $\xi^1(y)$, and $\xi^1(z)$,*
    (d) *The image of $\eta^1_{xyz}$ is contained in the closed simplex formed as the closure of a connected component of the complement of the projective lines containing $\ell_{xy}, \ell_{yz}$, and $\ell_{xz}$.*

*Proof.* Claim (1) follows from the definition of $C_{xyz}$. That $\eta^1_{xyz}(w)$ and $\eta^2_{xyz}(w)$ are continuous follows readily from hyperconvexity, in particular from the transversality of $\xi^2(w)$ and $\xi^3(w)$ to $P_{xyz}$ for all $w \in \partial \Gamma$. Injectivity of $\eta^1_{xyz}$ follows from the consequence of hyperconvexity that $\xi^2(u) \cap \xi^2(v) = \{0\}$ for all $u \neq v$.

To show that $\eta^2_{xyz}$ is injective it suffices to show that $\xi^3(u) \cap \xi^3(v) \cap P_{xyz}$ is not a projective line for any $u \neq v$. To see this, observe that since $\ell_{xy}$ is entirely contained in $\Omega^1_\xi$ by Prop. 3.7, Lemma 3.6 shows for any such $u$ and $v$ that $\xi^3(u) \cap \xi^3(v) \cap (\xi^1(x) \oplus \xi^1(y)) = \emptyset$. Dimension counting, we conclude that $\xi^3(u) \cap \xi^3(v) \cap P_{xyz}$ can not be a projective line. This gives (2b).

We prove (2a) for $[x, y]$. The other cases follow symmetrically. We verify the hyperconvex Frenet curve condition through its formulation in terms of intersections (Prop. 2.4). For general position on $(x, y)$, let $w_1, ..., w_k$ be ordered so $x < w_1 < ... < w_k < y$ and $n_1, ..., n_k$ be natural numbers so $\sum_{i=1}^k n_i = 3$. Then by Lemma 2.5,

$$\left(\bigcap_{i=1}^k \eta^{3-n_i}_{xyz}(w_i)\right)^\perp = \sum_{i=1}^k (\xi^{*n_i}(w_i) + P^\perp_{xyz})$$
$$\supset \left(\xi^{*3}(x) \cap \xi^{*3}(y) \cap \xi^{*3}(z)\right) + \xi^{*n_1}(w_1) + ... + \xi^{*n_k}(w_k)$$
$$= \mathbb{R}^4.$$

We conclude that the general position condition holds on $(x, y)$. Extensions to endpoints are handled similarly to the extension of hyperconvexity to the (straightforward) endpoint $\xi^1(y)$ in Lemma 4.6. Limit compatibility is verified essentially as at $\xi^1(y)$ in Lemma 4.6 as well, using transversality of all projective lines $\xi^2(x)$ and projective planes $\xi^3(x)$ to $P_{xyz}$. This proves (2a).

We now turn to analyzing how the segments $\eta^1_{xyz}([x,y])$, $\eta^1_{xyz}([y,z])$, and $\eta^1_{xyz}([z,x])$ glue together. We shall analyze $\eta^1_{xyz}([x,y])$ and $\eta^1_{xyz}([y,z])$ at $\xi^1(y)$; the other endpoints follow identically due to the symmetric roles of $x, y,$ and $z$.

Up to rotation, the possible configurations of $\eta^1_{xyz}([x,y])$ and $\eta^1_{xyz}([y,z])$, together with $\ell_{xy}$ and $\ell_{yz}$ are shown in Figure 14.

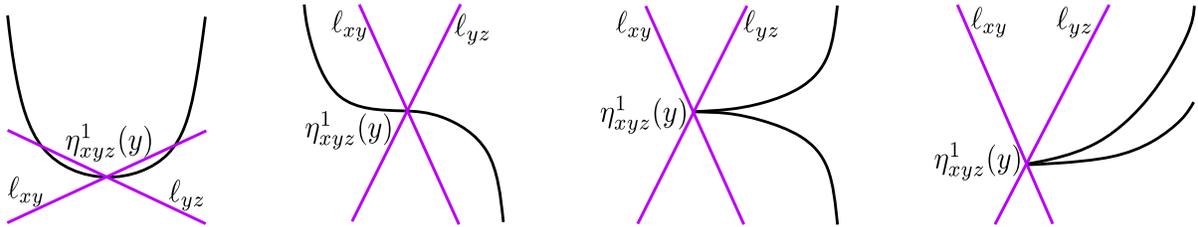

FIGURE 14. Possible configurations of the image of $\eta^1_{xyz}$ at $\eta^1_{xyz}(y)$, with the lines $\ell_{xy}$ and $\ell_{yz}$ in purple. The containment of all but 2 points of each of $\ell_{xy}$ the purple lines in $\Omega^1_\rho$ is how we deduce local structure of $C_{xyz}$ at singular points. We refer to these configurations as type (1)-(4) from left to right.

Observe that as $P_{xyz}$ is transverse to $\Xi^2 \subset \partial\Omega_\xi$ at every point of intersection, every transverse crossing of $\eta^1_{xyz}(\partial\Gamma - \{x,y,z\})$ must change connected component of $\Omega_\xi$. This observation, together with the consequence of Prop. 3.7 that $\ell_{xy}$ and $\ell_{yz}$ are entirely contained in $\Omega^1_\rho$ rules out cases (1) and (2). Case (4) is ruled out by noting from tangent vector counting that if case (4) were to occur, then a transverse crossing of $\Xi^2 \cap P_{xyz}$ could be made, both sides of which are in $\Omega^2_\xi$. We conclude that a balanced cusp occurs at $\xi^1(y)$, finishing (2c). The final claim, (2d), is a consequence of (2c). □

### 4.3. Projections of Frenet Curves.
The assertions on subspaces of vector spaces associated to hyperconvex Frenet curves in Lemma 2.5 concern only subspaces with quite restricted forms. We shall need to work with slightly more general subspaces. In this subsection, we first present a basic linear-algebraic lemma on subspaces of vector spaces in general position that can be used to re-organize some expressions of subspaces into forms where hyperconvexity may be applied. We then prove a general lemma on certain "projections" of hyperconvex Frenet curves to subspaces that remain hyperconvex Frenet curves.

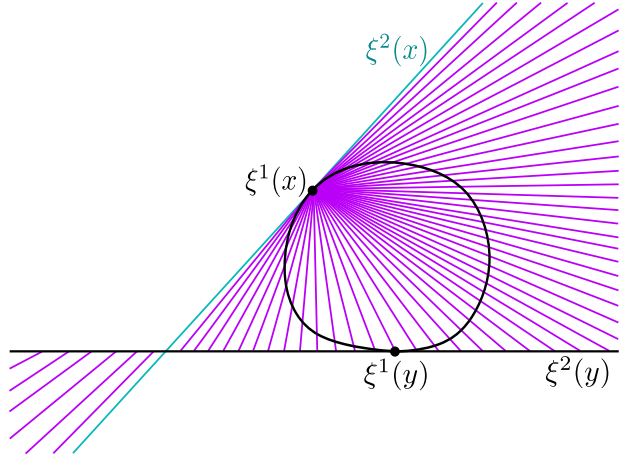

FIGURE 15. The map $\tau_{yx}^1$ for $l = 1$ and $k = 2$ is obtained by drawing lines from $\xi^1(x)$ to $\xi^2(y)$ and projecting points $p$ in $\xi^1(\partial\Gamma)$ to the intersection of the line $\ell_p$ between $\xi^1(x)$ and $\xi^2(y)$ through $p$ to the intersection of $\ell_p$ with $\xi^2(y)$. The picture in this case is that of stereographic projection.

4.3.1. *A Linear Algebra Lemma.* We shall have repeated use of the following lemma on subspaces of vector spaces *in general position*. We give a careful proof in the appendix because counting dimensions of subspaces formed with both intersections and sums is delicate.

**Lemma 4.8** (General Position Projection)**.** *Let $V$ be an $n$-dimensional vector space. Let $X, Y, W_1, ..., W_k$ be subspaces of $V$ so that $V = X \oplus Y = X \oplus W_1 \oplus W_2 \oplus \cdots \oplus W_k$.*
*Then for $j = 1, ..., k$,*

$$(4.3) \qquad Y \cap (X \oplus W_1 \oplus \cdots \oplus W_j) = \bigoplus_{i=1}^{j} ((X \oplus W_i) \cap Y),$$

*and the sum on the right is direct. In particular, $Y = \bigoplus_{i=1}^{k}((X \oplus W_i) \cap Y)$, and the sums in Eq. 4.3 have dimension $\dim(Y \cap (X \oplus W_1 \oplus \cdots \oplus W_j)) = \sum_{i=1}^{j} \dim W_j$ for $j = 1, ..., k$.*

The name of Lemma 4.8 is due to the corresponding picture in projective space.

4.3.2. *The Projection Lemma.* The following gives a second way, in addition to the Frenet Restriction Lemma 3.1, that a hyperconvex Frenet curve $\xi$ induces other hyperconvex Frenet curves in subspaces $\xi^k(y)$. Namely, one obtains such curves by placing a viewer on the image of $\xi$ and visually projecting an appropriate curve $\xi^j(\partial\Gamma)$ to $\xi^k(y)$. This construction is useful in accessing subtle general position consequences of the hyperconvex Frenet condition.

**Lemma 4.9** (Frenet Projection)**.** *Let $\xi : \partial\Gamma \to \mathcal{F}(\mathbb{R}^n)$ be a hyperconvex Frenet curve, let $(x, y) \in \partial\Gamma^{(2)}$ and let $k + l = n$. Then the following map $\tau_{yx} = (\tau_{yx}^1, ..., \tau_{yx}^{k-1}) : \partial\Gamma \to \mathcal{F}(\xi^k(y))$* *is a hyperconvex Frenet curve:*

$$\tau_{yx}^j(z) = \begin{cases} (\xi^l(x) \oplus \xi^j(z)) \cap \xi^k(y) & z \neq x, y \\ \xi^{l+j}(x) \cap \xi^k(y) & z = x, \\ \xi^j(y) & z = y \end{cases} \quad (z \in \partial\Gamma, j = 1, ..., k-1)$$

See Figure 15.

*Proof.* We verify the hyperconvex Frenet curve condition in terms of its formulation in terms of sums. Let $n_1, ..., n_p$ be so that $\sum_{i=1}^p n_i = k$ and let $(w_1, ..., w_p) \in \partial\Gamma$ be distinct. We do case analysis on if $x$ or $y$ is in $\{w_1, ..., w_p\}$.

Suppose first that $x, y \notin \{w_1, ..., w_p\}$. Then the General Position Projection Lemma 4.8 and the hyperconvexity of $\xi$ show that

$$\bigoplus_{i=1}^p \tau_{yx}^{n_i}(w_i) = \bigoplus_{i=1}^p [(\xi^l(x) \oplus \xi^{n_i}(w_i)) \cap \xi^k(y)]$$
$$= \xi^k(y) \cap [\xi^l(x) \oplus \xi^{n_1}(w_1) \oplus \xi^{n_2}(w_2) \oplus ... \oplus \xi^{n_p}(w_p)]$$
$$= \xi^k(y).$$

Suppose next that $x \in \{w_1, ..., w_p\}$ and $y \notin \{w_1, ..., w_p\}$, and without loss of generality take $x = w_1$. To address this case, let $(v_1, ..., v_n)$ be a basis of $\mathbb{R}^n$ adapted to the flag $\xi(x)$, i.e. so that $v_i \in \xi^i(x)$ for $i = 1, ..., n$ (we adopt the convention $\xi^n(x) = \mathbb{R}^n$). For $j = 1, ..., k$, let $W_j(x) = \text{span}(v_{l+1}, ... v_{l+j})$. The point is that $\xi^{l+j}(x) = \xi^l(x) \oplus W_j(x)$ for $j = 1, ..., k$ and

$$\xi^l(x) + W_{n_1}(x) + \xi^{n_2}(w_2) + ... + \xi^{n_p}(w_p) = \xi^{l+n_1}(x) + \xi^{n_2}(w_2) + ... + \xi^{n_p}(w_p) = \mathbb{R}^n,$$

and, in particular, the above sums are direct.

So the General Position Projection Lemma and hyperconvexity of $\xi$ show

$$\bigoplus_{i=1}^p \tau_{yx}^{n_i}(w_i) = [(\xi^l(x) \oplus W_{n_1}(x)) \cap \xi^k(y)] + \left(\bigoplus_{i=2}^p [(\xi^l(x) \oplus \xi^{n_i}(w_i)) \cap \xi^k(y)]\right)$$
$$= \xi^k(y) \cap [\xi^l(x) \oplus W_{n_1}(x) \oplus \xi^{n_2}(w_2) \oplus ... \oplus \xi^{n_p}(w_p)]$$
$$= \xi^k(y).$$

To handle the case where $y \in \{w_1, ..., w_p\}$, we note that for $1 \leq j \leq k$ that $\xi^j(y) = \xi^k(y) \cap (\xi^j(y) \oplus \xi^l(x))$ by hyperconvexity of $\xi$. So when $y \in \{w_1, ..., w_p\}$, the analysis above goes through exactly as before with $y$ in place of, say, $w_p$. So $\tau_{yx}$ is hyperconvex.

We now verify limit compatibility. So let $\sum_{i=1}^p n_i = m < k$ and let $(w_1^n, ..., w_p^n)$ be a sequence of distinct elements of $\partial\Gamma$ that converge to $(w, ..., w)$ for some $w \in \partial\Gamma$. Without loss of generality, if $w_i^n = x$ for some $i$, let $i = 1$. Note that for $n$ sufficiently large, it is impossible for both $x$ and $y$ to be in $\{w_1, ..., w_p\}$. Using the equalities established in the proof of hyperconvexity, whether or not $y \in \{w_1^n, ..., w_p^n\}$,

$$\lim_{n \to \infty} \bigoplus_{i=1}^p \tau_{yx}^{n_i}(w_i^n) = \lim_{n \to \infty} \begin{cases} \xi^k(y) \cap [\xi^l(x) \oplus \xi^{n_1}(w_1^n) \oplus ... \oplus \xi^{n_p}(w_p^n)], & w_1^n \neq x, \\ \xi^k(y) \cap [\xi^{l+n_1}(x) \oplus \xi^{n_2}(w_1^n) ... \oplus \xi^{n_p}(w_p^n)], & w_1^n = x \end{cases}$$
$$= \begin{cases} \xi^k(y) \cap (\xi^l(x) \oplus \xi^m(w)) & w \neq y, x \\ \xi^k(y) \cap \xi^{l+m}(x) & w = x, \\ \xi^m(y) & w = y \end{cases}$$

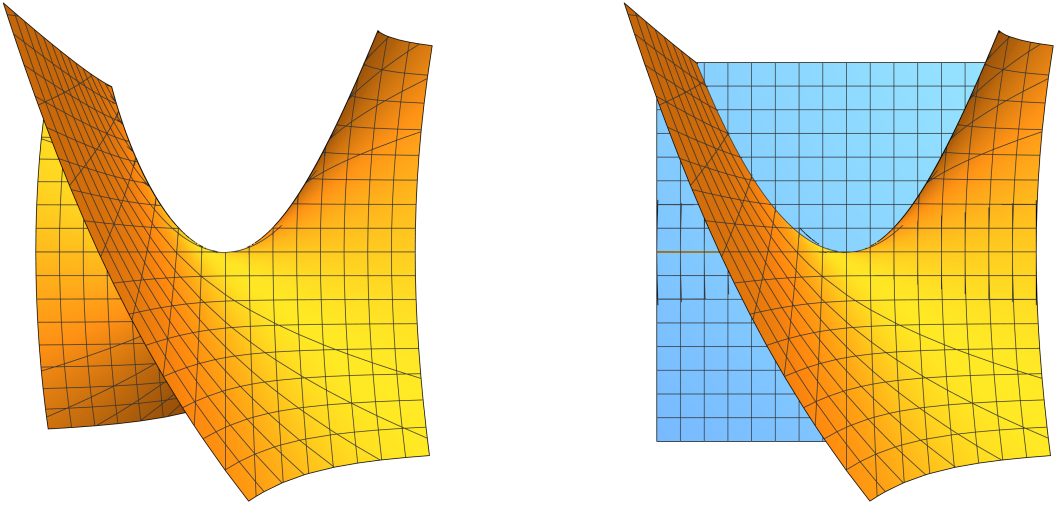

FIGURE 16. Left: rendering of a "top view" of $\partial\Omega_\xi$ for the Fuchsian example in an affine chart for $\mathbb{RP}^3$. Right: the top view of $\partial\Omega_\xi$ together with a plane of the form $\xi^3(x)$ containing the bottom point on $\Xi^1$ seen from this perspective. Note that the vertical projections of $\Xi^1$ lie inside the convex region bounded by the intersection of the plane with $\partial\Omega_\xi$. This is proved in generality in Lemma 4.11.

$$= \tau_{yx}^m(w)$$

In all cases, we used transversality of $\xi^l(x)$ and $\xi^{l+m}(x)$ with $\xi^k(y)$ in taking the limit. In the last case, hyperconvexity is used to deduce that $\xi^k(y) \cap (\xi^l(x) \oplus \xi^m(y)) = \xi^m(y)$. This completes the verification that that $\tau_{yx}$ is a hyperconvex Frenet curve. □

4.4. **The Top View of the Domain.** We analyze the general projections of hyperconvex Frenet curves introduced in §4.3.2 in our setting, and establish some qualitative features of how the images of these curves sit inside the convex domains $C_x$. The analysis gives rise to two foliations of the leaves $C_x$ $(x \in \partial\Gamma)$ of $\mathcal{F}_{\mathrm{pcf}}$, which are a key technical tool in the analysis of the intersection of $(\xi^3(x) \cap \xi^3(z)) \oplus \xi^1(y)$ with $\partial\Omega_\xi$.

4.4.1. *Top View.* In this paragraph, we examine the intersections of lines of the form $\xi^1(y) \oplus \xi^1(w)$, for $y$ fixed, with planes of the form $\xi^3(x)$.

The specialization of the construction of the Frenet Projection Lemma we consider is:

**Definition 4.10.** *For $(x,y) \in \partial\Gamma^{(2)}$, define $\tau_{xy} = (\tau_{xy}^1, \tau_{xy}^2) : \partial\Gamma \to \mathcal{F}(\xi^3(x))$ by*

$$\tau_{xy}^1(w) = \begin{cases} (\xi^1(y) \oplus \xi^1(w)) \cap \xi^3(x) & w \neq y, \\ \xi^2(y) \cap \xi^3(x) & y = w \end{cases}, \quad \tau_{xy}^2(w) = \begin{cases} (\xi^1(y) \oplus \xi^2(w)) \cap \xi^3(x) & w \neq y, \\ \xi^3(y) \cap \xi^3(x) & y = w \end{cases}.$$

*Denote the image of $\tau_{xy}^1$ by $C_{xy}$.*

Of course, when $w = x$ we have $\tau_{xy}^1(x) = \xi^1(x)$ and $\tau_{xy}^2(x) = \xi^2(x)$. See Figure 16.

**Lemma 4.11** (Top View). *In the above notation:*
  (1) *The map $\tau_{xy}$ is a hyperconvex Frenet curve $\partial\Gamma \to \mathcal{F}(\xi^3(x))$,*
  (2) *The image $C_{xy}$ of $\tau_{xy}^1$ is contained in $\Omega_\rho^1$, except for the points $\xi^1(x)$ and $\xi^2(y) \cap \xi^3(x)$.*

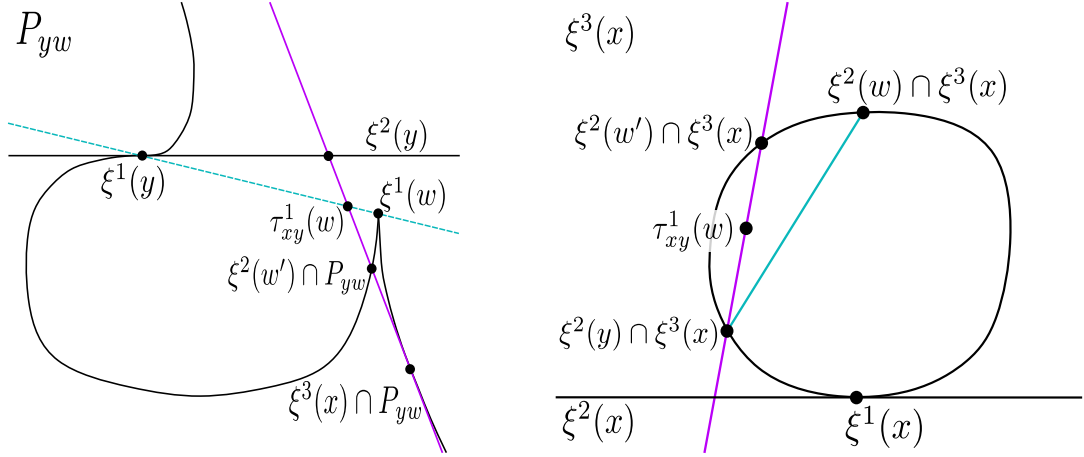

FIGURE 17. Left: the structure of the intersection $\xi^3(x) \cap P_{yw}$, as determined by our analysis of $P_{yw} \cap \partial \Omega_\xi$. Right: the structure in the illustration on the left forces $\tau^1_{xy}(w)$ to lie inside of the half of $C_x$ outside of the teal line segment that does not contain $\xi^1(x)$ in its boundary. This implies that $\tau^1_{xy}$ travels along its image in the "same direction" as the parameterization of $\partial C_x$ by $\xi^2(w) \cap \xi^3(x)$.

(3) For every $w \in \partial \Gamma$, there is a point $w'$ in the connected component of $\partial \Gamma - \{w, y\}$ that does not contain $x$ so that $\tau^1_{xy}(w)$ is in the line segment in $C_y$ with endpoints $\xi^2(w') \cap \xi^3(x)$ and $\xi^2(y) \cap \xi^3(x)$.
(4) Let $\mathcal{U}_1$ and $\mathcal{U}_2$ be the closures of the connected components of $\partial \Gamma - \{x, y\}$. Then maps obtained by restricting $\tau_{xy}$ and $\xi_{x,2}$ to individual components of $U_1$ and $U_2$ are all hyperconvex Frenet curves, in the following sense. If $g_{xy}, f_{xy} \in \{\tau_{xy}, \xi_{x,2}\}$ then the map $h_{xy}$ defined by $h_{xy}|_{\mathcal{U}_1} = f_{xy}$ and $h_{xy}|_{\mathcal{U}_2} = g_{xy}$ is a hyperconvex Frenet curve.

Claim (3) above controls the location of $\tau^1_{xy}(w)$ in $C_y$, which is useful in applications. It is illustrated in Figure 17, Right. Claim (4) shows that the directions of the parameterizations of $C_{xy}$ by $\tau_{xy}$ and $\partial C_x$ by $\xi^1_{x,2}$ agree, as is suggested by Figure 16, Right.

*Proof.* By the Frenet Projection Lemma 4.9, $\tau_{xy}$ is a hyperconvex Frenet curve. The image of $\tau_{xy}$ is contained in $C_y$ except at the endpoints by Prop. 3.7. This gives (1) and (2).

Now let $w \in \partial \Gamma$, and examine how $\xi^3(x)$ intersects the plane $P_{yw} = \xi^2(w) \oplus \xi^1(y)$. See Figure 17, Left. Recall the notation that $\eta_{yw} : \partial \Gamma \to \mathcal{F}(P_{yw})$ parameterizes $\partial \Omega_\xi \cap P_{yw} - (\xi^2(y) - \{\xi^1(y)\})$. Let $\mathcal{C}_1$ be the connected component of $\partial \Gamma - \{y, w\}$ containing $x$ and $\mathcal{C}_2$ the other component. Let $w'$ be the unique element of $\mathcal{C}_2$ supplied by Lemma 4.6.(2f) so that $\eta^2_{yw}(x)$ contains $\eta^1_{yw}(w')$. Then from the structure of $P_{yw} \cap \partial \Omega_\xi$, we see that $\tau^1_{xy}(w)$ is in the line segment in $\xi^3(x) \cap \partial \Omega^1_\xi$ between $\xi^2(w') \cap \xi^3(x)$ and $\xi^2(y) \cap \xi^3(x)$. See Figure 17, Left. Claim (3) then follows from that $w' \in \mathcal{C}_2$. The final claim (4) is a direct consequence of the other claims. □

The above leads to two foliations of the convex domains $C_x$ with some useful features.

**Definition 4.12.** *Define the convex arcs $\mathcal{A}^+_{xy} = \tau_{xy}((x, y))$ and $\mathcal{A}^-_{xy} = \tau_{xy}((y, x))$.*

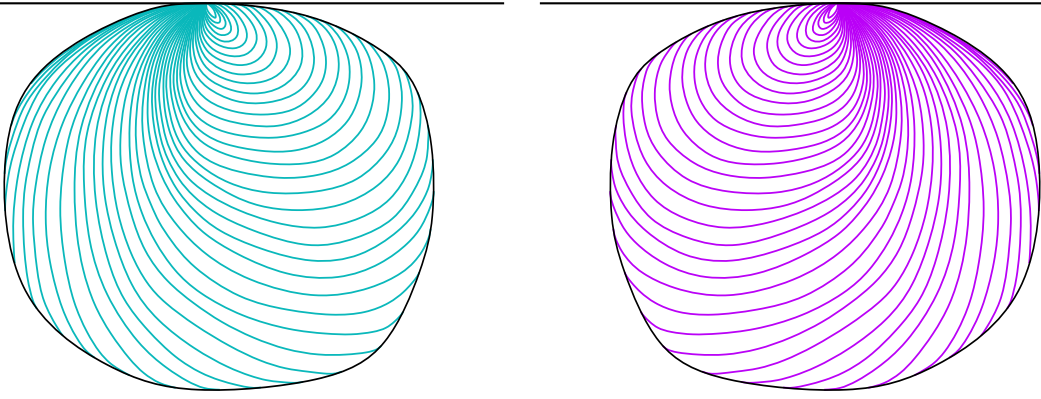

FIGURE 18. The two foliations of the convex domain $C_x$ in Lemma 4.13. For any $y \in \partial\Gamma - \{x\}$, the image of the curve $\tau_{xy}$ is the union of one teal curve and one purple curve, the union of which is the boundary of a convex domain.

The following is now immediate from Prop. 3.7 and Lemma 4.11, using our knowledge of the arrangement of the arcs $\mathcal{A}_{xy}^{\pm}$ from Lemma 4.11.(3)-(4). See Figure 18 for illustrations of the foliations described in the lemma.

**Lemma 4.13** (Arc Foliations). *Let $x \in \partial\Gamma$. The convex domain $C_x$ has foliations by convex arcs specified by $C_x = \bigsqcup_{y \in \partial\Gamma - \{x\}} \mathcal{A}_{xy}^{+} = \bigsqcup_{y \in \partial\Gamma - \{x\}} \mathcal{A}_{xy}^{-}$.*

*For any $y_1, y_2 \in \partial\Gamma$ so that $(x, y_1, y_2)$ is positively oriented,*
$$\mathcal{A}_{xy_1}^{+} \cup \mathcal{A}_{xy_2}^{-} \cup \{\xi^2(w) \cap \xi^3(x) \mid w \in [y_1, y_2]\} \cup \{\xi^1(x)\}$$
*is the boundary of a properly convex, strictly convex domain in $\xi^3(x)$.*

4.5. **The Final Intersection Pattern.** Lemma 4.13 is the additional structure we use to analyze the final type of intersections of planes with $\partial\Omega_\xi$, namely those of the form $(\xi^3(x) \cap \xi^3(z)) \oplus \xi^1(y)$ with $(x, y, z) \in \partial\Gamma^{(3)+}$. We now turn to this. We shall prove that Figure 13, Right is qualitatively accurate.

To this end, let $(x, y, z) \in \partial\Gamma^{(3)+}$, denote $Q_{xyz} = (\xi^3(x) \cap \xi^3(z)) \oplus \xi^1(y)$, and define the map $\nu_{xyz} = (\nu_{xyz}^1, \nu_{xyz}^2) : \partial\Gamma \to \mathcal{F}(Q_{xyz})$ by
$$\nu_{xyz}^1(w) = \xi^2(w) \cap Q_{xyz}, \quad \nu_{xyz}^2(w) = \xi^3(w) \cap Q_{xyz} \quad (w \in \partial\Gamma).$$

**Proposition 4.14** (($\xi^3(x) \oplus \xi^3(z)) \oplus \xi^1(y)$ Pattern). *Let $(x, y, z) \in \partial\Gamma^{(3)+}$ and let $Q_{xyz}$ and $\eta_{xyz}$ be as above. Then:*
  (1) *$Q_{xyz} \cap \partial\Omega_\xi = \nu_{xyz}^1(\partial\Gamma)$ and is entirely contained in an affine chart,*
  (2) *$\nu_{xyz}^1$ is injective and $\nu_{xyz}^2(w_1) = \nu_{xyz}^2(w_2)$ if and only if $w_1 = w_2$ or $\{w_1, w_2\} = \{x, z\}$,*
  (3) *The restrictions of $\nu_{xyz}$ to $[x, y]$, and $[y, z]$, and $(z, x)$ are hyperconvex Frenet curves, as is the restriction of $\nu_{xyz}$ to $[z, w_0]$ or $[w_0, x]$ for any $w_0 \in (x, z)$.*
  (4) *At the points $\xi^2(x) \cap Q_{xyz}$ and $\xi^2(z) \cap Q_{xyz}$ and $\xi^1(y)$,*
      (a) *The restriction of $\nu_{xyz}$ to a sufficiently small neighborhood of $x$ or $z$ is a hyperconvex Frenet curve,*
      (b) *$\nu_{xyz}^1$ has a balanced cusp at $\xi^1(y)$.*

*Proof.* We begin with (2). First note that the injectivity of $\nu_{xyz}^1$ follows from the consequence of hyperconvexity that $\xi^2(x) \cap \xi^2(y) = \emptyset$ for all $(x, y) \in \partial\Gamma^{(2)}$. Next observe that the map

$(\partial\Gamma - \{x,z\}) \to (\xi^3(x) \cap \xi^3(z) - \{\xi^2(z) \cap \xi^3(x), \xi^2(x) \cap \xi^3(z)\})$ given by $w \mapsto \xi^3(w) \cap \xi^3(x) \cap \xi^3(z)$ is a homeomorphism as a consequence of the Frenet Restriction Lemma. This implies that for all distinct $w_1$ and $w_2$ in $\partial\Gamma - \{x,z\}$ that $\nu^2_{xyz}(w_1) \neq \nu^2_{xyz}(w_2)$. The case where exactly one of $w_1$ and $w_2$ is either $x$ or $z$ is similar. This proves (2).

We next address (3), and verify the hyperconvex Frenet curve condition through its characterization in terms of intersections. We prove all general position claims for points in the interiors of our intervals. The endpoints are confirmed analogously.

We first address $[x,y]$; the same argument goes through for $[y,z]$. First let $w_1 < w_2 < w_3$ be points in $(x,y)$. Then we see that

$$(\nu^2_{xyz}(w_1) \cap \nu^2_{xyz}(w_2) \cap \nu^2_{xyz}(w_3))^\perp$$
$$= (\xi^{*1}(w_1) \oplus Q^\perp_{xyz}) + (\xi^{*1}(w_2) \oplus Q^\perp_{xyz}) + (\xi^{*1}(w_3) \oplus Q^\perp_{xyz})$$
$$\supset \xi^{*1}(w_1) + \xi^{*1}(w_2) + \xi^{*1}(w_3) + Q^\perp_{xyz}$$
$$= \xi^{*1}(w_1) + \xi^{*1}(w_2) + \xi^{*1}(w_3) + ((\xi^{*1}(x) \oplus \xi^{*1}(z)) \cap \xi^{*3}(y)).$$

We claim the final sum above is $(\mathbb{R}^4)^*$, and hence the original intersection is trivial. From dimension counting, it suffices to show that $(\xi^{*1}(w_1) \oplus \xi^{*1}(w_2) \oplus \xi^{*1}(w_3)) \cap \xi^{*3}(y)$ does not contain $Q^\perp_{xyz}$.

By the General Position Projection Lemma,

$$(\xi^{*1}(w_1) \oplus \xi^{*1}(w_2) \oplus \xi^{*1}(w_3)) \cap \xi^{*3}(y)$$
$$= [(\xi^{*1}(w_1) \oplus \xi^{*1}(w_2)) \cap \xi^{*3}(y)] \oplus [(\xi^{*1}(w_1) \oplus \xi^{*1}(w_3)) \cap \xi^{*3}(y)].$$

In the notation of §4.4.1, with respect to the dual hyperconvex Frenet curve $\xi^*$,

$$(\xi^{*1}(w_1) \oplus \xi^{*1}(w_2)) \cap \xi^{*3}(y) \in \mathcal{A}^-_{yw_1}, \quad (\xi^{*1}(w_1) \oplus \xi^{*1}(w_3)) \cap \xi^{*3}(y) \in \mathcal{A}^-_{yw_1},$$
$$Q^\perp_{xyz} = (\xi^{*1}(x) \oplus \xi^{*1}(z)) \cap \xi^{*3}(y) \in \mathcal{A}^+_{yx}.$$

By Lemma 4.13, $\mathcal{A}^-_{yw_1}$ and $\mathcal{A}^+_{yx}$ are contained in the boundary of a properly convex, strictly convex domain in $\xi^{*3}(y)$. See Figure 19. The line of interest is now a secant line in the boundary of a properly convex strictly convex domain. It thus meets the boundary of the convex domain, which contains $Q^\perp_{xyz}$, in exactly the two points $\tau_{yw_1}(w_2)$ and $\tau_{yw_1}(w_3)$. As these points are not $Q^\perp_{xyz}$, this gives the desired identity.

The other general position conditions to verify on $[x,y]$ are:

$$\nu^2_{xyz}(w_1) \cap \nu^1_{xyz}(w_2) = \{0\}, \quad \nu^1_{xyz}(w_1) \cap \nu^2_{xyz}(w_2) = \{0\} \quad (w_1 < w_2 \in [x,y]).$$

Slight variations on the above argument resolve these cases. For instance, Figure 19, Middle is the final configuration of points used to show $\nu^1_{xyz}(w_1) \cap \nu^2_{xyz}(w_2) = \{0\}$. So hyperconvexity holds on $[x,y]$.

We next address limit compatibility. We have already seen that all $\xi^3(x)$ and $\xi^2(x)$ ($x \in \partial\Gamma$) are transverse to $Q_{xyz}$, as are $\xi^3(w_1) \cap \xi^3(w_2)$ for $w_1, w_2 \in [x,y]$. Then for any sequence $(w_1^n, w_2^n)$ of pairwise distinct points converging to $(w,w)$ for some $w \in [x,y]$, we see from our observations on transversality and dimension counting that

$$\lim_{n \to \infty} \nu^2_{xyz}(w_1^n) \cap \nu^2_{xyz}(w_2^n) = \lim_{n \to \infty} (\xi^3(w_1^n) \cap Q_{xyz}) \cap (\xi^3(w_2^n) \cap Q_{xyz})$$
$$= \lim_{n \to \infty} (\xi^3(w_1^n) \cap \xi^3(w_2^n)) \cap Q_{xyz}$$
$$= \xi^2(w) \cap Q_{xyz}$$

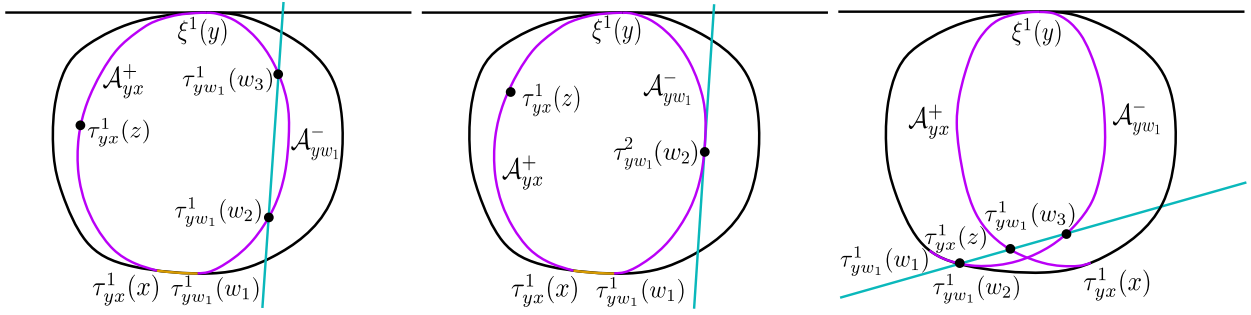

FIGURE 19. Configurations of points in $\xi^{*3}(y)$ with '$*$' suppressed in labels. The point of the general position argument for $[x, y]$ is that the purple and yellow curves pictured left and center bound properly convex and strictly convex domains, which obstructs the teal lines from containing $\tau^1_{xy}(z)$. The right illustration, which does not actually occur, is of the possible configuration that obstructs this argument for $[z, x]$. We pursue a different argument for this interval.

$$= \nu^1_{xyz}(w).$$

So the restriction of $\nu_{xyz}$ to $[x, y]$ is a hyperconvex Frenet curve.

We next address $(z, x)$, which is the connected component of $\partial \Gamma - \{x, z\}$ that does not contain $y$, and the other hyperconvex Frenet curve restriction claims. The arguments made for general position for $[x, y]$ do not directly apply here, see Figure 19, Right.

We once again verify the intersection characterization of the hyperconvex Frenet curve condition, and convert transversality conditions to a containment question for the dual curve. For the first collection of indices under consideration, $(2, 2, 2)$, we must show that for $w_1 < w_2 < w_3$ in $(z, x)$,

$$(\xi^{*1}(w_1) \oplus Q^\perp_{xyz}) + (\xi^{*1}(w_2) \oplus Q^\perp_{xyz}) + (\xi^{*1}(w_3) \oplus Q^\perp_{xyz}) = (\mathbb{R}^4)^*.$$

Arguing as before, it suffices to show that

$$(\xi^{*1}(w_1) \oplus \xi^{*1}(w_2) \oplus \xi^{*1}(w_3)) + Q^\perp_{xyz} = (\mathbb{R}^4)^*,$$

which is equivalent to $Q^\perp_{xyz}$ not being contained in $\xi^{*1}(w_1) \oplus \xi^{*1}(w_2) \oplus \xi^{*1}(w_3)$. Following the notation of §4.2.2, we denote $\xi^{*1}(w_1) \oplus \xi^{*1}(w_2) \oplus \xi^{*1}(w_3)$ by $P^\perp_{w_1 w_2 w_3}$.

Our strategy is to examine where $P^\perp_{w_1 w_2 w_3}$ intersects the projective line $\xi^{*1}(x) \oplus \xi^{*1}(z)$, which contains $Q^\perp_{xyz} = \xi^{*3}(y) \cap (\xi^{*1}(x) \oplus \xi^{*1}(z))$.

A consequence of hyperconvexity is that every projective plane in the connected collection

$$\mathcal{C}^*_{xz} = \{\xi^{*1}(w_1) \oplus \xi^{*1}(w_2) \oplus \xi^{*1}(w_3) \mid w_1 < w_2 < w_3 \in (z, x)\}$$
$$\cup \{\xi^{*1}(w_1) \oplus \xi^{*2}(w_2) \mid w_1, w_2 \in (z, x)\} \cup \{\xi^{*3}(w) \mid w \in (z, x)\}$$

is transverse to $\xi^{*1}(x) \oplus \xi^{*1}(z)$, and contains neither $\xi^{*1}(z)$ nor $\xi^{*1}(z)$. Since the intersections between transverse projective lines and planes vary continuously, there is a single connected component of $\ell^*_{xz} = (\xi^{*1}(x) \oplus \xi^{*1}(z)) - \{\xi^{*1}(x), \xi^{*1}(z)\}$ that contains all points of intersection of planes in $\mathcal{C}^*_{xz}$ and $\xi^{*1}(x) \oplus \xi^{*1}(z)$.

We next note that from the behavior of dev$_{\text{pctf}}$ established in Prop. 3.7 applied to $\xi^*$,

$$\xi^{*1}(x) \oplus \xi^{*1}(z) = \{\xi^{*1}(x), \xi^{*1}(z)\} \sqcup \left( \bigsqcup_{v \in (x,z)} \xi^{*3}(v) \cap \ell^*_{xz} \right) \sqcup \left( \bigsqcup_{v \in (z,x)} \xi^{*3}(v) \cap \ell^*_{xz} \right).$$

Since $\xi^{*3}(w)$ is in $\mathcal{C}^*_{xz}$ for all $w \in (z,x)$ and $y$ is in the oriented interval $(x,z)$ that does not contain $w$, we conclude that $Q \cap (\xi^1(x) \oplus \xi^1(z)) \neq Q^\perp_{xyz}$ for all $Q \in \mathcal{C}^*_{xz}$. This establishes the first general position requirement. The other general position conditions for $(z,x)$ follow from the relevant planes being contained in $\mathcal{C}^*_{xz}$. Limit compatibility in $(z,x)$ follows as in the previous analyses of intersections.

We now turn to the local behavior near endpoints. The main observation is that though projective planes such as $\xi^{*1}(x) \oplus \xi^{*2}(z)$ that are not transverse to $\xi^{*1}(x) \oplus \xi^{*1}(z)$ may be made, their construction requires *both* $x$ and $z$. To be more precise, let $w_0 \in (z,x)$. Then hyperconvexity implies that, as above, every projective plane in the compact connected set $\overline{\mathcal{C}^*_{zw_0}}$ is transverse to $\xi^{*1}(x) \oplus \xi^{*1}(z)$ and intersects this line in the closure of one common component of $\ell^*_{xz}$. The same argument as above then shows $\nu_{xyz}$ is hyperconvex on $[z, w_0]$, and analogously on $[w_0, x]$.

To handle the final claimed hyperconvexity, we note that $\overline{C^*_{zw_0}}$ is compact and contained in the open subset of $\text{Gr}_3((\mathbb{R}^4)^*)$ consisting of projective planes transverse to $Q^\perp_{xyz}$. It then follows for some open interval $(z', w')$ containing $[z_0, w_0]$, by the same argument as above, that $\nu_{xyz}$ restricted to $(z', w')$ is a hyperconvex Frenet curve. This finishes the proofs of (3) and (4a).

We now come to how the endpoints of our convex intervals glue together. From the structure of interections of planes of the form $\xi^2(a) \oplus \xi^1(b)$ with $\partial \Omega_\xi$ proved in §4.2.1, we see that in the projective lines $(\xi^2(x) \oplus \xi^1(y)) \cap Q_{xyz}$ and $(\xi^2(z) \oplus \xi^1(y)) \cap Q_{xyz}$ that $\xi^1(y)$ is an isolated point of intersection with $\partial \Omega_\xi$, a punctured neighborhood of which is entirely contained in $\Omega^1_\xi$. The possible configurations, up to renormalization, are once again shown in Figure 12.

Every case but (3) and (4) is ruled out by component-changing analysis with our above observation. Case (4) is ruled out by counting numbers of intersections of tangent lines on each side of transverse crossings of $\partial \Omega_\xi$. We conclude that $\nu^1_{xyz}(\partial \Gamma)$ has a balanced cusp at $\xi^1(y)$. Local configurations at $\eta^1_{xyz}(x)$ and $\eta^1_{xyz}(z)$ are determined by hyperconvexity of $\nu_{xyz}$ on neighborhoods of $x$ and $z$, finishing (4). The final claim, (1), follows from what we have shown. $\square$

### 4.6. The Front View.
We now explain the view of $\Xi^1$ as seen by a point in $\Xi^2$.

Let us set notation for this subsection. Let $(x, y) \in \partial \Gamma^{(2)}$. Define $P_{xy} = \xi^2(y) \oplus \xi^1(x)$ and $p_{xy} = \xi^2(x) \cap \xi^3(y)$. Define the map $\sigma_{xy} = (\sigma^1_{xy}, \sigma^2_{xy}) : (\partial \Gamma - \{x\}) \mapsto \mathcal{F}(P_{xy})$ by

(4.4) $$\sigma_{xy}(w) = (p_{xy} \oplus \xi^1(w)) \cap P_{xy}, \quad \sigma^2_{xy}(w) = (p_{xy} \oplus \xi^2(w)) \cap P_{xy}.$$

**Lemma 4.15** (Front View). *Let $(x, y) \in \partial \Gamma^{(2)}$. In the notation above,*

(1) $\sigma^1_{xy}$ *and* $\sigma^2_{xy}$ *are continuous injections,*
(2) *The restrictions of $\sigma_{xy}$ to $(x, y]$ and $[y, x)$ are hyperconvex Frenet curves,*
(3) *The image $\sigma^1_{xy}(\partial \Gamma - \{x\})$ is $C^1$ at $\xi^1(y)$ and crosses its tangent line $\xi^2(y)$ at $\xi^1(y)$.*

See Figure 20. We remark that the deduction of the structure of $\sigma_{xy}$ near $\xi^1(y)$ has a considerably different flavor than in the examinations of intersections before.

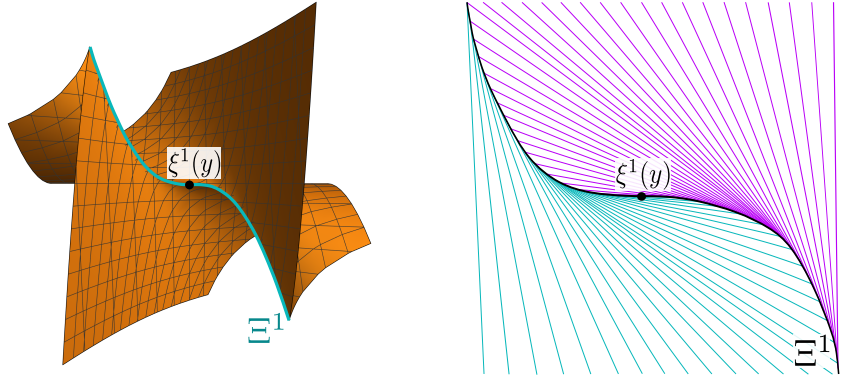

FIGURE 20. Left: annotated rendering of the "Front View" of the Fuchsian domain in our affine chart. Right: Illustration of lines in the canonical foliation of $\partial\Omega_\xi$ from the viewpoint of a point on $\Xi^2$.

*Proof.* Continuity is a consequence of transversality of $p_{xy} \oplus \xi^1(w)$ and $p_{xy} \oplus \xi^2(w)$ to $P_{xy}$ for all $w \neq x$. Injectivity follows from the general position property of $\xi$. This proves (1).

We turn to Claim (2) now. We verify the hyperconvex Frenet curve condition in terms of its formulation in terms of sums. We first address general position on $(x, y)$. So let $w_1, ..., w_k \in (x, y)$ be positively ordered and $n_1, ..., n_k$ be so that $\sum_{i=1}^k n_k = 3$. By Lemma 2.5,

$$p_{xy} + \xi^{n_1}(w_1) + \cdots + \xi^{n_k}(w_k) = [\xi^2(y) \cap \xi^3(x)] + \xi^{n_1}(w_1) + \cdots + \xi^{n_k}(w_k) = \mathbb{R}^4.$$

So by the General Position Projection Lemma 4.8,

$$\bigoplus_{i=1}^k \sigma_{xy}^{n_i}(w_i) = \bigoplus_{i=1}^k [(p_{xy} \oplus \xi^{n_i}(w_i)) \cap P_{xy}]$$
$$= \left[ p_{xy} \oplus \bigoplus_{i=1}^k \xi^{n_i}(w_i) \right] \cap P_{xy}$$
$$= P_{xy}.$$

So the general position condition holds on $(x, y)$. General position for combinations of subspaces in $(x, y]$ containing the endpoint $y$ is a consequence of the synthetic bound given by Lemma 4.11.(3).

For limit compatibility, let $(w_1^n, w_2^n)$ be a sequence of distinct elements of $(x, y]$ converging to $(z, z)$ for some $z \in (x, y]$. Then using the above and transversality of various subspaces,

$$\lim_{n \to \infty} \sigma_{xy}^1(w_1^n) \oplus \sigma_{xy}^1(w_2^n) = \lim_{n \to \infty} [p_{xy} \oplus \xi^1(w_1^n) \oplus \xi^1(w_2^n)] \cap P_{xy}$$
$$= [p_{xy} \oplus \xi^2(z)] \cap P_{xy}$$
$$= \sigma_{xy}^2(z).$$

So the restriction of $\sigma_{xy}$ to $(x, y]$ is a hyperconvex Frenet curve. The other restriction in (2) is symmetric.

In an appropriate affine chart for $P_{xy}$, the local structure of $\sigma_{xy}(\partial\Gamma - \{x\})$ at $\xi^1(y)$ has one of the forms shown once again in Figure 12 above. We will show configuration (2) must occur, which implies Claim (3). In the following, we work entirely in an affine chart $\mathcal{A}$ for $\mathbb{RP}^3$ chosen so that:

(1) $\xi^3(x)$ is the plane at infinity,
(2) $\xi^1(y)$ is at the origin,
(3) $\xi^2(y)$ is the (horizontal) $x$-axis,
(4) $\xi^1(y) \oplus \xi^1(x)$ is the (vertical) $y$-axis,
(5) $\xi^1(y) \oplus p_{xy}$ is the (forward) $z$-axis.

Then the lines in $\mathcal{A}$ containing $p_{xy}$ are exactly those parallel to the $z$-axis, and $\xi^3(y) \cap \mathcal{A}$ is the $(xz)$-plane. For points of the form $\xi^1(w)$ in $\mathcal{A}$, the point $\sigma^1_{xy}(w)$ is described as the projection in $\mathcal{A}$ to the $(xy)$-plane.

To see that there is a tangent line crossing at $\xi^1(y)$, suppose otherwise. As $\lim_{z \to y} \sigma^1_{xy}(z) = \xi^1(y)$, using our description of $\sigma^1_{xy}$ in $\mathcal{A}$ we can then pick $x'_n \in (x, y)$ and $x''_n \in (y, x)$ so that $\xi^1(x'_n)$ and $\xi^1(x''_n)$ have the same vertical coordinate in $\mathcal{A}$. In particular, $(\xi^1(x'_n) \oplus \xi^1(x''_n)) \cap \xi^3(y)$ occurs in the plane at infinity for all $n$. This shows that all limit points of $(\xi^1(x'_n) \oplus \xi^1(x''_n)) \cap \xi^3(y)$ lie in the plane at infinity, and in particular cannot be $\xi^1(y)$. See Figure 21, Left. This contradicts the consequence of Lemma 2.5 that $\lim_{n \to \infty} (\xi^1(x'_n) \oplus \xi^1(x''_n)) \cap \xi^3(y) = \xi^1(y)$, and rules out configurations (3) and (4).

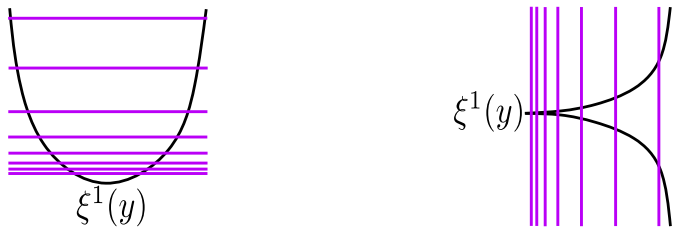

FIGURE 21. We constrain the structure of the "front view" near the point $\xi^1(y)$ by observing that the local configurations above result in incompatibilities of limits of subspaces formed from sums of entries on the hyperconvex Frenet with the limit compatibility condition.

We argue similarly that a horizontal crossing occurs at $\xi^1(y)$. Otherwise, we could take $x'_n$ and $x''_n$ as before, but so that $\xi^1(x'_n)$ and $\xi^1(x''_n)$ have the same $x$-coordinate. Then any limit of $\xi^1(x'_n) \oplus \xi^1(x''_n)$ would need to be vertical in our affine chart, in contradiction to the fact that $\lim_{n \to \infty} \xi^1(x'_n) \oplus \xi^1(x''_n) = \xi^2(y)$. See Figure 21, Right. This rules out configuration (1). We conclude that the local configuration (2) occurs at $\xi^1(y)$, completing the proof. □

4.7. **Deduction of Thm. B and Codimension 1 Theorems.** By Lemma 4.4, we have by now completely classified intersections of projective planes with $\Omega_\xi$. Thm. B is then an immediate corollary of Props. 4.6, 4.7, and 4.14. In particular, a corollary of the structure results above for intersections is that the intersection of a projective plane $P \subset \mathbb{RP}^3$ with $\Omega^1_\rho$ contains a properly embedded convex domain $C$ if and only if $P = \xi^3(x)$, in which case $C = C_x$. The codimension-1 claims in Thms. A and E immediately follow, and the codimension-1 claims in Thm. F follow similarly.

## 5. Foliation Rigidity

We prove the codimension-2 classification results in Thms. A and F in this section. The general setup is to take a $\rho(\Gamma)$-invariant foliation or trifoliation $\mathcal{G}$ of $\Omega^1_\rho$ or $\Omega^2_\rho$ whose leaves are proper convex subsets of projective lines then look at endpoint projections in a foliation chart. The constraints we obtain on the invariant (tri)foliations from applying group elements to these endpoints are substantial but insufficient to prove the classification theorems directly. To make further progress, we use the fine structural features of $\partial \Omega_\rho$ discussed in §4.

We have seen in §4 that the qualitative structure of $\partial \Omega_\rho$ is quite different near points in $\Xi^1$ and $\Xi^2$ and when viewed from points in $\Omega_\rho^1$ and $\Omega_\rho^2$. These differences result in a need for some technical care, book-keeping, and case analysis in the basic setup.

## 5.1. Codimension-2 Foliations of the Foliated Region. We begin with $\Omega_\rho^1$. We prove:

**Theorem 5.1** (Exactly Two: $\Omega_\rho^1$). *Let $\rho : \Gamma \to \mathrm{PSL}_4(\mathbb{R})$ be Hitchin. Let $\mathcal{G}$ be a $\rho(\Gamma)$-invariant foliation of $\Omega_\rho^1$ so that every leaf of $\mathcal{G}$ is a properly embedded properly convex subset of a projective line or a once-punctured projective line. Then $\mathcal{G} = \mathcal{G}_{\mathrm{pcf}}$ or $\mathcal{G} = \mathcal{G}_{\mathrm{tcf}}$.*

5.1.1. *Discontinuities of Endpoints.* In this paragraph, let $\mathcal{G}$ be a foliation of $\Omega_\rho^1$ whose leaves are properly embedded segments in projective lines and let $U \cong [-1, 1] \times \mathbb{D}$ be a foliation chart for $\mathcal{G}$. From the hypothesis of proper embedding, we may define $\phi_F : \mathbb{D} \to \partial \Omega_\rho$ and $\phi_B : \mathbb{D} \to \partial \Omega_\rho$ to be the mappings defined by assigning to $p \in \mathbb{D}$ the endpoint of the leaf of $\mathcal{G}$ through $(p, 0)$ in the direction of $\{1\} \times \mathbb{D}$ or $\{-1\} \times \mathbb{D}$, respectively. We remark that because transversality is an open condition and lines containing two continuously varying uniformly separated points vary continuously, we may take the cross-section $\mathbb{D}$ in $\mathbb{RP}^3$ to be an open subset of a projective plane.

Of course, $\phi_F$ and $\phi_B$ may be discontinuous. We classify their possible discontinuities here. The locations of the endpoints of intervals on $\partial \Omega_\rho$ are salient. Accordingly, we define:

**Definition 5.2.** *Denote:*
  (1) $O = \{x \in \mathbb{D} \mid \phi_F(x) = \phi_B(x)\}$,
  (2) $K_F = \phi_F^{-1}(\Xi^1) - O$ and $K_B = \phi_B^{-1}(\Xi^1) - O$ and $K = K_F \cap K_B$,
  (3) $V_F = \phi_F^{-1}(\Xi^2)$ and $V_B = \phi_B^{-1}(\Xi^2)$,
  (4) $C_F = K_F \cap V_B$ and $C_B = K_B \cap V_F$,
  (5) $\mathscr{L} = K_F \cup K_B \cup O$.

Note that $O$ is characterized as the collection of $x \in \mathbb{D}$ contained in a leaf of $\mathcal{G}$ that is a once-punctured projective line.

**Proposition 5.3** (Discontinuity Control). *With the above notations:*
  (1) *$V_F$ and $V_B$ are open and the restrictions $\phi_F|_{V_F}$ and $\phi_B|_{V_B}$ are continuous,*
  (2) *$C_F$ and $C_B$ are open in $K_F$ and $K_B$, respectively. The restrictions $\phi_F|_{C_F}$ and $\phi_B|_{C_B}$ are continuous,*
  (3) *$K$, $K_F$, and $K_B$ are closed. The restrictions $\phi_F|_K$ and $\phi_B|_K$ are continuous,*
  (4) *$O \cap V_F = O \cap V_B = \emptyset$, the closure $\overline{O}$ of $O$ satisfies $\overline{O} \subset O \cup K$, and $O$ is open in $O \cup K$. The restrictions of $\phi_B$ and $\phi_F$ to $O$ are continuous.*

See Figure 22 for examples of discontinuities of endpoints, to motivate the limitations of the above statement. Roughly, Prop. 5.3 says that endpoints of leaves behave well except when they move from $\Xi^2$ to $\Xi^1$ or change from having two endpoints to having only one.

*Proof.* First observe that the lines $\ell_p \in \mathrm{Gr}_2(\mathbb{R}^4)$ through $p \in \{0\} \times \mathbb{D}$ vary continuously since lines containing two continuously varying and uniformly separated points vary continuously.

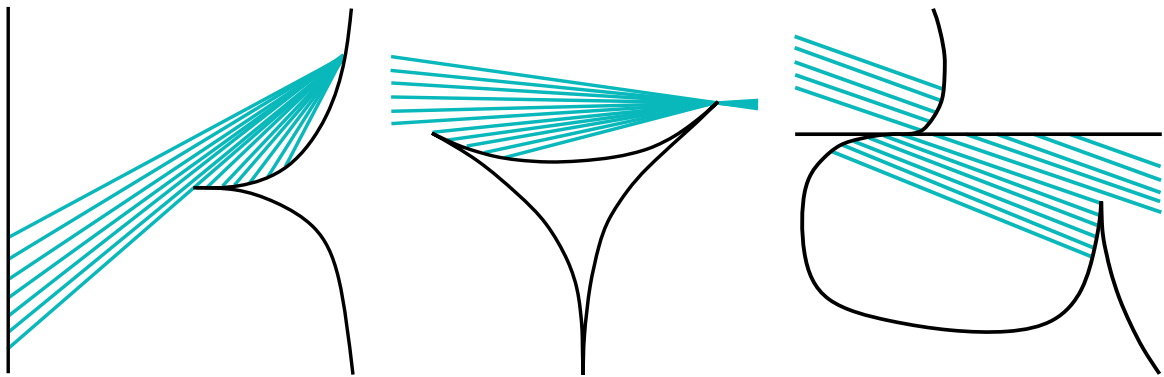

FIGURE 22. Some endpoint discontinuities in local foliations. Endpoints may jump between $\Xi^1$ and $\Xi^2$, or between having one and two endpoints on $\Xi^1$.

Next, suppose $p \in V_F$ and examine $\phi_F(p) \in \Xi^2$. As $\Xi^2$ is locally convex and $C^1$, and since $\Omega^1_\rho$ is on the convex side of $\Xi^2$, the line $\ell_p$ is transverse to $\Xi^2$ at $\phi_F(p)$. So there is a neighborhood $W_1 \subset \Xi^2$ of $\phi_F(p)$ and a neighborhood $W_2 \subset \mathbb{D}$ of $p$ so the forward rays through $q \in W_2$ intersect $W_1$ in continuously varying points. Since $\partial\Omega_\rho - W_1$ is compact and hence uniformly separated from the line segment between $p$ and $\phi_F(p)$, there is a sub-neighborhood $W_3 \subset W_2$ of $p$ so that the forward rays from $q \in W_3$ first intersect $\partial\Omega_\rho$ in $W_1$. This proves $V_F$ is open and $\phi_F|_{V_F}$ is continuous. The analogous claims for $\phi_B$ hold by symmetry.

The relative openness claims for $C_F$ and $C_B$ are now immediate. To see that $\phi_F|_{C_F}$ is continuous, note by hyperconvexity of $\xi$ that any $p \in C_F$ has $\ell_p$ intersect $\Xi^1$ in exactly one point. Since all tangent lines to $\Xi^1$ are contained in $\partial\Omega_\rho$, the line $\ell_p$ is transverse to $\Xi^1$, which is $C^1$ at $\phi_F(p)$. Continuity now follows as before.

To address continuity on $K$, all line segments between distinct points in $\Xi^1$ are contained in the foliation $\mathcal{G}_{\text{tcf}}$, and continuously varying leaves of $\mathcal{G}_{\text{tcf}}$ intersect $\Xi^1$ with continuously varying forward and backward endpoints by inspection of the developing map $\text{dev}_{\text{pctf}}$. Note this also implies that $K$ is closed.

We now show the closedness of $K_F$. The claim for $K_B$ is symmetric. All that is not immediate from the above is that limits of points in $K_F$ and $K_B$ are not in $O$. Let us now prove this.

By closedness of $K$, it suffices to consider sequences of points $x_n \in C_F$ converging to $x \in \mathbb{D}$, with backwards endpoints $\phi_B(x_n) = \xi^2(a_n) \cap \xi^3(b_n) \in \Xi^2$ and $\phi_F(x_n) \in \Xi^1$. By compactness of $\partial\Omega_\rho$, after subsequence there are limiting points $p_{F,\infty} = \xi^1(c_n) \in \Xi^1$ of $\phi_F(x_n)$ and $p_{B,\infty} \in \partial\Omega_\rho$ of $\phi_B(x_n)$. Each is contained in the limiting line, so it now suffices to show that $p_{B,\infty} \neq p_{F,\infty}$.

Suppose otherwise, for contradiction. Because the two-argument map $\xi^1(\cdot)$ is a homeomorphism, this would imply that $c_n, a_n,$ and $b_n$ all have the same limit $c_\infty$. Since $c_n, a_n,$ and $b_n$ are triples of points, the order conditions in Lemma 2.5 are satisfied and this would imply

$$\lim_{n \to \infty} (\xi^2(a_n) \cap \xi^3(b_n)) \oplus \xi^1(c_n) = \xi^2(c_\infty),$$

which is contained in $\partial\Omega_\rho$. This is impossible because the limiting line must be in $O$, and so intersect $\Omega^1_\rho$.

For the statements on $O$, first note that any line that intersects $\Xi^2$ and contains a point in $\Omega_\rho^1$ must intersect $\partial\Omega_\rho$ in at least two points. This is because $\Xi^2$ is $C^1$ and locally convex, with $\Omega_\rho^1$ on the convex side, so that any line $\ell$ that meets $\Xi^2$ at $p$ from the side of $\Omega_\rho^1$ must change connected components of $\Omega_\rho$ at $p$. In order to return to the $\Omega_\rho^1$ side, there must be another point of intersection of $\ell$ and $\partial\Omega_\rho$. This shows that $O \cap V_F = O \cap V_B = \emptyset$.

To see that $\overline{O} \subset O \cup K$, note that the properties of nontrivially intersecting $\Xi^1$ and of being contained entirely in $\overline{\Omega_\rho^1}$ are closed in $\mathrm{Gr}_2(\mathbb{R}^4)$. Since, as mentioned above, every line containing a segment in $\Omega_\rho^1$ and intersecting $\Xi^2$ must intersect $\Omega_\rho^2$, this forces any element of $\overline{O}$ to intersect $\Xi^1$ and not intersect $\Xi^2$, which implies $\overline{O} \subset O \cup K$. The openness claim is trivial because $K$ is closed. Continuity of $\phi_F = \phi_B$ on $O$ now follows by similar analysis to continuity on $C_F$. $\square$

5.1.2. *Basic Leaf Inclusions.* Let us now give the basic means we have to produce leaves of $\mathcal{G}_{\mathrm{pcf}}$ and $\mathcal{G}_{\mathrm{tcf}}$ inside of arbitrary $\rho(\Gamma)$-invariant foliations whose leaves are properly embedded projective line segments.

To disambiguate the leaves of $\mathcal{G}_{\mathrm{tcf}}$ contained in a common projective line, for $(a,b) \in \partial\Gamma^{(2)}$, let $\mathcal{G}_{\mathrm{t}}(a,b) \subset \ell_{ab}$ be the connected component so that $\xi^3(c) \cap \mathcal{G}_t(a,b) \neq \emptyset$ for all $c \in \partial\Gamma$ so that $(a,c,b) \in \partial\Gamma^{(3)+}$. Because the action of $\Gamma$ on $\partial\Gamma$ preserves orientations, $\rho(\gamma)\mathcal{G}_t(a,b) = \mathcal{G}_t(\gamma a, \gamma b)$. Additionally, let $\mathcal{G}_{\mathrm{p}}(a,b)$ denote the leaf of $\mathcal{G}_{\mathrm{pcf}}$ between $\xi^1(a)$ and $\xi^2(b) \cap \xi^3(a)$.

**Lemma 5.4** (Basic Leaf Inclusions). *Suppose $\mathcal{G}$ is a $\rho(\Gamma)$-invariant foliation of $\Omega_\rho^1$ whose leaves are properly embedded, properly convex domains in projective lines or once-punctured projective lines, and let $\gamma \in \Gamma - \{e\}$.*

(1) *If $\ell$ is a leaf of $\mathcal{G}$ and $\ell = \mathcal{G}_{\mathrm{t}}(\gamma^-, a)$ then $\mathcal{G}_{\mathrm{t}}(\gamma^-, \gamma^+)$ is a leaf of $\mathcal{G}$. If $\ell = \mathcal{G}_{\mathrm{t}}(a, \gamma^-)$ then $\mathcal{G}_{\mathrm{t}}(\gamma^+, \gamma^-)$ is a leaf of $\mathcal{G}$.*
(2) *If $\ell$ is a leaf of $\mathcal{G}$ with endpoints $\xi^1(\gamma^-)$ and $\xi^2(a) \cap \xi^3(b)$ with $(\gamma^-, b, a) \in \partial\Gamma^{(3)+}$ then $\mathcal{G}_{\mathrm{t}}(\gamma^-, \gamma^+)$ is a leaf of $\mathcal{G}$. If $(\gamma^-, b, a)$ is negatively oriented, then $\mathcal{G}_{\mathrm{t}}(\gamma^+, \gamma^-)$ is a leaf of $\mathcal{G}$. If instead $b = \gamma^-$ then $\mathcal{G}_{\mathrm{p}}(\gamma^-, \gamma^+)$ is a leaf of $\mathcal{G}$.*
(3) *If $\ell$ is a leaf of $\mathcal{G}$ with an endpoint of the form $\xi^2(\gamma^-) \cap \xi^3(a)$ then $\mathcal{G}_{\mathrm{p}}(\gamma^+, \gamma^-)$ is a leaf of $\mathcal{G}$.*
(4) *If $\ell$ is a leaf of $\mathcal{G}$ with a unique endpoint, which is of the form $\xi^1(\gamma^+)$, then $\mathcal{G}_{\mathrm{t}}(\gamma^+, \gamma^-)$ and $\mathcal{G}_{\mathrm{t}}(\gamma^-, \gamma^+)$ are leaves of $\mathcal{G}$.*

**Remark 5.5.** *No leaf inclusions are given in the case of a leaf with an endpoint of the form $\xi^2(a) \cap \xi^3(\gamma^-)$, in contrast to when the roles of $a$ and $\gamma^-$ are interchanged. The argument we use below in fact completely breaks in this case and gives no leaf inclusions, as mentioned in §1.3.1. The issue is that for $a, b, c \neq \gamma^-$,*

$$\lim_{n\to\infty} \rho(\gamma^n)[(\xi^2(a) \cap \xi^3(\gamma^-)) \oplus (\xi^2(b) \cap \xi^3(c))] = [\xi^2(\gamma^+) \cap \xi^3(\gamma^-)] \oplus \xi^1(\gamma^+) = \xi^2(\gamma^+).$$

*Since $\xi^2(\gamma^+)$ is contained in $\partial\Omega_\rho$, we can hope for no leaf inclusions from this method.*

*This asymmetry of the roles of values of $x$ and $y$ in endpoints of the form $\xi^2(x) \cap \xi^3(y)$ is responsible for much of the difficulty of proving Thm. A.*

*Proof.* For (1), we first note that $\lim_{n\to\infty} \rho(\gamma^n)\mathcal{G}_{\mathrm{t}}(\gamma^-, a) = \mathcal{G}_{\mathrm{t}}(\gamma^-, \gamma^+)$, since by hyperconvexity the limiting line only intersects $\partial\Omega_\rho$ at $\xi^1(\gamma^-)$ and $\xi^1(\gamma^+)$. Then, examining a foliation chart around a point in $\mathcal{G}_{\mathrm{t}}(\gamma^-, \gamma^+)$, we see $\mathcal{G}_{\mathrm{t}}(\gamma^-, \gamma^+)$ is a leaf of $\mathcal{G}$. The other sub-claim is analogous.

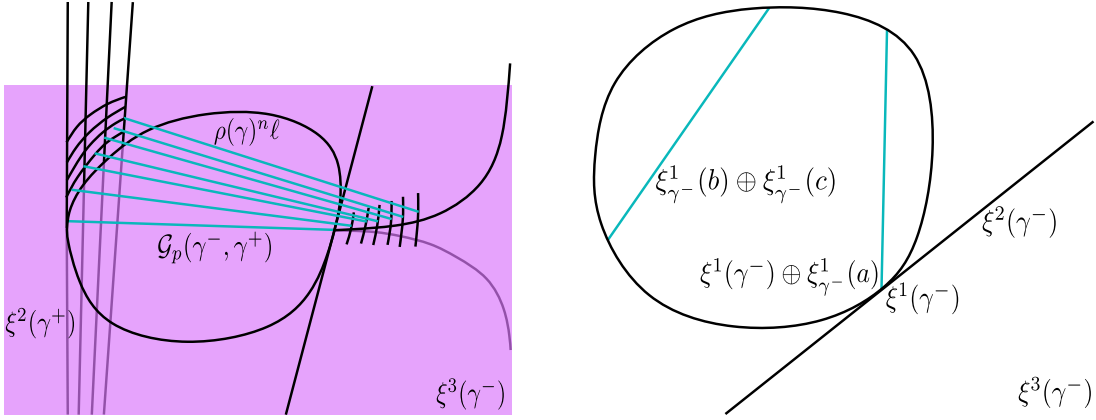

FIGURE 23. Left: in the proof of claim (3) of Lemma 5.4, the local convexity and regularity of $\Xi^2$ ensures that a limiting segment exists in $\xi^3(\gamma^+)$, and the structure of the intersection of $\xi^3(\gamma^+)$ with $\Omega^1_\rho$ causes the limiting leaf to be $\mathcal{G}_p(\gamma^+, \gamma^-)$. Right: segments contained in $\xi^3(\gamma^-) \cap \Omega^1_\rho$ can only intersect $\xi^2(\gamma^-)$ at $\xi^1(\gamma^-)$.

For (2), first note that $a \neq \gamma^-$ because $\xi^2(\gamma^-) \subset \partial \Omega_\rho$. See Figure 23, Right. In the case that $b = \gamma^-$ we have $\rho(\gamma^n)\ell \subset \xi^3(\gamma^-)$ for all $n$ and from the qualitative structure of $\xi^3(\gamma^-) \cap \partial \Omega_\rho$ it follows that $\rho(\gamma^n)\ell$ converges to $\mathcal{G}_p(\gamma^-, \gamma^+)$. The case of claim (2) where $b = \gamma^-$ follows.

Let us next take $b \neq \gamma^-$ and $(\gamma^-, b, a)$ to be positively oriented, with the orientation chosen to disambiguate between the sides of the limiting leaf. The final sub-case is analogous. We begin by remarking that from the qualitative structure of $P_{a\gamma^-} = \xi^2(a) \oplus \xi^1(\gamma^-)$ in Lemma 4.6 that $\ell \cap \xi^3(c)$ is nonempty for all $c \in (\gamma^-, b)$. It follows that, for all $n$, a neighborhood $U \subset (\gamma^-, \gamma^+)$ satisfies that for all $c \in U$, the intersection $\rho(\gamma^n)\ell \cap \xi^3(c)$ is nonempty.

Note next that

$$\lim_{n \to \infty} \rho(\gamma^n)[(\xi^2(a) \cap \xi^3(b)) \oplus \xi^1(\gamma^-)] = \lim_{n \to \infty} [(\xi^3(\gamma^n a) \cap \xi^2(\gamma^n b)) \oplus \xi^1(\gamma^-)]$$
$$= \xi^1(\gamma^+) \oplus \xi^1(\gamma^-)$$

by limit compatibility and transversality of $\xi^1(\gamma^+)$ and $\xi^1(\gamma^-)$. Arguing as in (1), up to subsequence $\rho(\gamma^n)\ell$ converges to $\mathcal{G}_t(\gamma^-, \gamma^+)$ or $\mathcal{G}_t(\gamma^+, \gamma^-)$. Fixing $c \in U$ and using transversality of $\xi^3(c)$ and $\xi^1(\gamma^+) \oplus \xi^1(\gamma^-)$, we see that any subsequential limit must intersect $\xi^3(c)$, and hence be $\mathcal{G}_t(\gamma^-, \gamma^+)$. As in (1), examining a foliation chart around a point $p \in \mathcal{G}_t(\gamma^-, \gamma^+)$ shows $\mathcal{G}_t(\gamma^-, \gamma^+)$ is a leaf of $\mathcal{G}$.

For (3), let the endpoints of $\ell$ be $p = \xi^2(\gamma^-) \cap \xi^3(a)$ and $q$. Note that $q \neq \xi^1(\gamma^-)$ and $a \neq \gamma^-$. We have the possibilities

(a) $q = \xi^1(b)$ and $b \neq \gamma^-$,
(b) $q = \xi^2(b) \cap \xi^3(c)$.

Note in the second case that $b \neq \gamma^-$ and $c \neq \gamma^-$, reasoning as in the argument for Claim (2). In either case, $\lim_{n \to \infty} \rho(\gamma^n)q = \xi^1(\gamma^+)$ and $\lim_{n \to \infty} \rho(\gamma^n)p = \xi^2(\gamma^-) \cap \xi^3(\gamma^+)$. Denote the limit of $\rho(\gamma^n)p$ by $p_\infty$.

Since $\ell_\infty := \xi^1(\gamma^+) \oplus (\xi^2(\gamma^-) \cap \xi^3(\gamma^+))$ is transverse to $\partial \Omega_\rho$ at $p_\infty$, which is a $C^1$ point of $\partial \Omega_\rho$, we see that for $n$ sufficiently large $\rho(\gamma^n)\ell$ has segments of length uniformly bounded below in a reference metric that emanate from $\partial\Omega_\rho$ at $\rho(\gamma^n)p$ on the convex side. Hence the Hausdorff limit of $\rho(\gamma^n)\ell$ is an interior segment in $\Omega_\rho^1$ inside $\ell_\infty$ containing a one-sided neighborhood of $p_\infty$ in $\ell_\infty$. See Figure 23. From the qualitative structure of $\Omega_\rho^1 \cap \xi^3(\gamma^+)$, we must have $\lim_{n\to\infty} \rho(\gamma^n)\ell = \mathcal{G}_p(\gamma^-, \gamma^+)$. It then follows that $\mathcal{G}_p(\gamma^-, \gamma^+)$ is a leaf of $\mathcal{G}$.

For (4), note that $\ell$ can not be contained in $\xi^3(\gamma^+)$ due to the qualitative structure of $\xi^3(\gamma^+) \cap \Omega_\rho$: a once-punctured projective line does not fit inside a properly convex domain. So $\ell$ contains a point $p_2 \notin \xi^3(\gamma^+)$. Then

$$\lim_{n\to\infty} \rho(\gamma^n)\ell = \lim_{n\to\infty} (\rho(\gamma^n)\xi^1(\gamma^+) \oplus \rho(\gamma^n)p_2) = \xi^1(\gamma^+) \oplus \xi^1(\gamma^-).$$

Note that $\rho(\gamma^n)\ell$ is projective line punctured at $\xi^1(\gamma^+)$ for all $n$, and so for any $p$ in $\mathcal{G}_t(\gamma^+, \gamma^-)$, the leaves $\rho(\gamma^n)\ell$ contain segments arbitrarily close to any small segment in $\mathcal{G}_t(\gamma^+, \gamma^-)$ containing $p$. Examining a foliation chart for $\mathcal{G}$ about $p$ shows that $\mathcal{G}_t(\gamma^+, \gamma^-)$ is a leaf of $\mathcal{G}$. The containment of $\mathcal{G}_t(\gamma^-, \gamma^+)$ is identical. □

We shall use the following basic observation to verify that a foliation sharing sufficiently many leaves with $\mathcal{G}_{\text{pcf}}$ or $\mathcal{G}_{\text{tcf}}$ is in fact either $\mathcal{G}_{\text{tcf}}$ or $\mathcal{G}_{\text{pcf}}$.

**Lemma 5.6** (Dense Leaves). *Let $\mathcal{G}$ be a $\rho(\Gamma)$-invariant foliation of $\Omega_\rho^1$ whose leaves are properly embedded properly convex subsets of projective lines or once punctured projective lines. Suppose that there is a neighborhood $U \subset \partial \Gamma^{(2)}$ so that for a dense set $S$ of $U$, for all $(a,b) \in S$ either $\mathcal{G}_p(a,b)$ is a leaf of $\mathcal{G}$ or $\mathcal{G}_t(a,b)$ is a leaf of $\mathcal{G}$. Then $\mathcal{G} = \mathcal{G}_{\text{pcf}}$ or $\mathcal{G} = \mathcal{G}_{\text{tcf}}$.*

*Proof.* This follows from the following observations:
  (1) The union of two nowhere dense sets in $U$ is nowhere dense,
  (2) A foliation $\mathcal{G}$ satisfying the hypotheses with leaves all of $\mathcal{G}_p(a,b)$ or all of $\mathcal{G}_t(a,b)$ for all $(a,b)$ in a dense subset $S$ of a neighborhood $U \subset \partial \Gamma^{(2)}$ must contain all such leaves of the same type for all $(a,b) \in U$. This is proved by examining a foliation chart for $\mathcal{G}$.
  (3) The action of $\Gamma$ on $\partial \Gamma^{(2)}$ is topologically transitive. □

5.1.3. *The Exactly Two Theorem for $\Omega_\rho^1$.* We enter the core of the argument of Thm. 5.1. We retain the notations of the rest of §5.1.

*Proof of the Exactly Two Theorem for $\Omega_\rho^1$.* Let $\mathcal{G}$ be a $\rho(\Gamma)$-invariant foliation of $\Omega_\rho^1$ by properly embedded properly convex subsets of projective lines. Take a foliation chart $B \cong (-1, 1) \times \mathbb{D}$ around an interior point $p \in \Omega_\rho^1$. We do case analysis on the endpoints of leaves through $B$. The two broad cases are if every leaf has an endpoint on $\Xi^1$ or if there is a leaf with both endpoints on $\Xi^2$. The first case has a number of sub-cases based on the possible behaviors of endpoints, each of which is handled by controlling endpoint discontinuities then applying the Basic Leaf Inclusions Lemma 5.4 and the Dense Leaves Lemma 5.6.

The generic case, where there is a leaf with both endpoints on $\Xi^2$, is by far the most involved. The reader may wish to refer to the sketch of this case in §1.3.1.

**Case 1** (Always a $\Xi^1$ Endpoint). *Every leaf of $\mathcal{G}$ has an endpoint on $\Xi^1$.*

This is equivalent to $\mathbb{D} = \mathscr{L}$ for every cross-section $\mathbb{D}$ for a foliation chart for $\mathcal{G}$, in the notation of §5.1.1.

**Step 1.** *If $K$ is not totally disconnected, then $\mathcal{G} = \mathcal{G}_{\text{tcf}}$.*

*Proof.* If $K$ is not totally disconnected, then $K$ contains a connected component $\mathcal{C}$ that is not a point. In this case, one of the two maps $\phi_F$ and $\phi_B$ on $\mathcal{C}$ must surject an open interval in $\Xi^1$ by continuity. Then Basic Leaf Inclusions (1) and the Dense Leaves Lemma implies that $\mathcal{G} = \mathcal{G}_{\text{tcf}}$. □

**Step 2** (One $\Xi^1$ Endpoint). *If $(K_F \cup K_B) - K$ has interior, then $\mathcal{G} \in \{\mathcal{G}_{\text{pcf}}, \mathcal{G}_{\text{tcf}}\}$.*

*Proof.* Without loss of generality assume that $C_F \cap \mathbb{D} \neq \emptyset$. After restriction, we may take $C_F = \mathbb{D}$, so $\phi_F$ and $\phi_B$ are both continuous on $\mathbb{D}$.

Since $\Xi^2$ is locally convex and $\Omega^1_\rho$ is on the convex side of $\Xi^2$, in any projective line $\ell$ that intersects $\partial \Omega_\rho$ transversely at $p \in \Xi^2$ there is a unique properly embedded segment in $\ell \cap \Omega^1_\rho$ with endpoint $p$. So identifying $\Xi^2$ with $\partial \Gamma^{(2)}$ by the map $\xi^2(a) \cap \xi^3(b) \mapsto (a,b)$ and identifying $\Xi^1$ with $\partial \Gamma$ by the map $\xi^1(x) \mapsto x$,

$$\phi : \mathbb{D} \to \partial \Gamma^{(2)} \times \partial \Gamma$$
$$p \mapsto (\phi_B(p), \phi_F(p))$$

is a continuous injection.

By continuity, the image of $\mathbb{D}$ under the projection $p_{\Xi^1}$ to the $\Xi^1$-factor of the image of $\phi$ is either a point or contains an interval. In the case where $p_{\Xi^1}$ has image a point, $\phi_B$ is a continuous injection $\mathbb{D} \to \partial \Gamma^{(2)}$, and so must surject an open set by invariance of domain. Then Basic Leaf Inclusions (3) and the Dense Leaves Lemma show that $\mathcal{G} = \mathcal{G}_{\text{pcf}}$.

On the other hand, if $p_{\Xi^1}$ surjects an interval $U$, then for every repelling fixed-point $\gamma^- \in U$, at least one of $\mathcal{G}_t(\gamma^+, \gamma^-), \mathcal{G}_t(\gamma^-, \gamma^+)$, and $\mathcal{G}_p(\gamma^-, \gamma^+)$ is contained in $\mathcal{G}$ by the Basic Leaf Inclusions Lemma. To handle the order of indices in the $\mathcal{G}_t$ case, we note that the projections of leaf labels of pole-pairs in $U$ to $\partial \Gamma^{(2)}/\mathbb{Z}_2$ given by forgetting the ordering are dense in $(U \times \partial \Gamma)/\mathbb{Z}^2$. Since the finite union of nowhere dense sets is nowhere dense, the unordered endpoints of leaves of $\mathcal{G}_{\text{tcf}}$ in $\mathcal{G}$ with one fixed local order or the endpoints of leaves of $\mathcal{G}_{\text{pcf}}$ in $\mathcal{G}$ are dense in a neighborhood in $\partial \Gamma^{(2)}/\mathbb{Z}_2$. This collection is then dense in a neighborhood of in $\partial \Gamma^{(2)}$ before $\mathbb{Z}_2$-projection. The Dense Leaves Lemma finishes the step. □

The previous steps reduce the case to the following:

**Step 3** (Once-Punctured $\mathbb{RP}^1$). *It is not possible for every leaf of $\mathcal{G}$ to have an endpoint on $\Xi^1$ and $K_B \cup K_F$ to have no interior in any cross-section $\mathbb{D}$.*

This case rules out the possibility that $\mathbb{D} = O \cup K$ with $O$ non-empty in every cross-section.

*Proof.* Assume otherwise. So our cross-section $\mathbb{D}$ for any foliation chart satisfies $\mathbb{D} = O \sqcup K$, and $O$ is open in $\mathbb{D}$. Our proof is motivated by the observation:

**Claim.** *$\phi_F = \phi_B$ is locally constant on $O$.*

*Proof.* From Basic Leaf Inclusions (4), density of attracting fixed-points in $\partial \Gamma$, minimality of the action of $\Gamma$ on $\partial \Gamma$, and the Dense Leaves Lemma, if $\phi_F(O)$ contains an interval then $\mathcal{G} = \mathcal{G}_{\text{tcf}}$. This would violate the contradiction hypothesis as $\mathcal{G}_{\text{tcf}}$ contains no once-punctured projective lines as leaves.

From the Discontinuity Control Lemma, $\phi_F$ is continuous on the open set $O$, and so must be locally constant on $O$. □

The complementary observation here is:

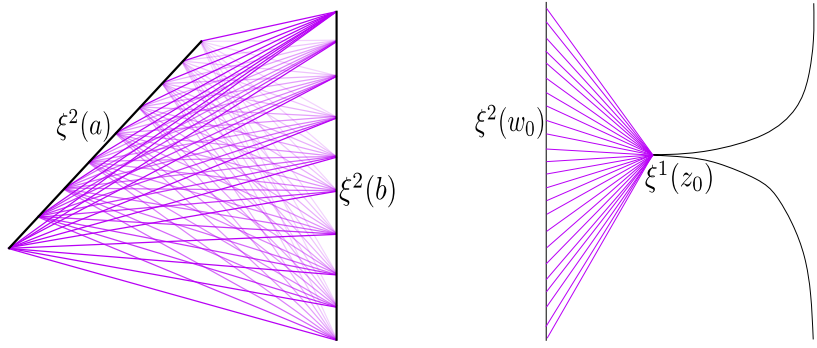

FIGURE 24. These yield no leverage in the proof of the Exactly Two Theorem, due to the phenomenon of Remark 5.5. The left illustration is a piece of a foliation modelled on the convex hull of two lines. The illustration on the right is the portion of a foliation given by the hull of a segment of $\xi^2(w_0)$ and $\xi^1(z_0)$. We call such portions of $\mathcal{G}$ *cusp-fans*.

**Claim.** *On any cross-section of a foliation chart of $\mathcal{G}$, the set $K$ is closed and totally disconnected.*

*Proof.* Closedness is contained in Lemma 5.3, and total disconnectedness is Step 1. □

In finishing the argument, the key point is now the classical (but nontrivial) plane topology theorem that if $C \subset \mathbb{C}$ is closed and separates two points $z_1, z_2 \in \mathbb{C}$, then a connected component of $C$ separates $z_1$ and $z_2$ (e.g. [21, Thm. 14.2]). This fact will also be useful below. We conclude $O \subset \mathbb{D}$ is connected, because its complement $K$ is closed and every connected component of $K$ is a point.

Now, our argument has shown that on *any* cross-section of a foliation chart, $O$ is connected, open, dense, and $\phi_F = \phi_B$ is constant on $O$. We conclude that every leaf of $\mathcal{G}$ that is a once-punctured projective line has a common boundary point $\xi^1(x_0)$. But this is impossible by group-invariance: there is a $\gamma \in \Gamma$ for which $\gamma x_0 \neq x_0$, so that $\rho(\gamma)\ell_0$ has sole endpoint $\xi^1(\gamma x_0)$, which is a contradiction. This completes the step, and with it the case. □

**Case 2** (Generic Properly Convex Case). *There is a leaf of $\mathcal{G}$ with both endpoints in $\Xi^2$.*

This is the main case. As indicated in §1.3.1, we proceed by first producing the existence of a single leaf of $\mathcal{G}_{\mathrm{pcf}}$, i.e. a leaf pointing directly into a singular point of $\partial\Omega_\xi$ along a one-sided tangent to a cusp. We then sufficiently constrain $\mathcal{G}$ to prevent a local model near this leaf built entirely out of the pieces in Figure 24 in such a way that we are unable to apply the Basic Leaf Inclusions Lemma and Dense Leaves Lemma to conclude that $\mathcal{G} \in \{\mathcal{G}_{\mathrm{pcf}}, \mathcal{G}_{\mathrm{tcf}}\}$. In the situation that escapes obstruction for the longest, we produce a possibly different leaf of $\mathcal{G}_{\mathrm{pcf}}$ that is inside of what we call a *cusp-fan*: a 2-dimensional piece of $\mathcal{G}$ formed by the convex hull of a line segment in $\partial\Omega_\rho$ and a point in $\Xi^1$ (Figure 24, Right). We arrange for this cusp-fan to not be the limit of other large cusp-fans in $\mathcal{G}$ on at least one side. The structure of the front-view in Lemma 4.15 gives the final obstruction.

We break the proof into steps.

**Step 1** (Improved Leaf). *$\mathcal{G}$ contains a leaf $\ell_p = \mathcal{G}_{\mathrm{p}}(z_0, w_0)$ of $\mathcal{G}_{\mathrm{pcf}}$ for some $(z_0, w_0) \in \partial\Gamma^{(2)}$.*

*Proof.* Write as $\ell_0$ the leaf assumed to exist with both endpoints in $\Xi^2$, and write the endpoints as $\xi^2(a) \cap \xi^3(b)$ and $\xi^2(c) \cap \xi^3(d)$. As above, $c, d \neq a$.

As $\Gamma$ is a uniform convergence group, there is a sequence $\gamma_n \in \Gamma$ and $w_0 \neq z_0$ in $\partial\Gamma$ so that
$$\lim_{n\to\infty} \gamma_n a = w_0, \qquad \lim_{n\to\infty} \gamma_n x = z_0 \qquad (\text{for all } x \neq z_0 \text{ in } \partial\Gamma).$$
In particular,
$$\lim_{n\to\infty} \rho(\gamma_n)(\xi^2(a) \cap \xi^3(b)) = \lim_{n\to\infty} \xi^2(\gamma_n a) \cap \xi^3(\gamma_n b) = \xi^2(w_0) \cap \xi^3(z_0),$$
$$\lim_{n\to\infty} \rho(\gamma_n)(\xi^2(c) \cap \xi^3(d)) = \lim_{n\to\infty} \xi^2(\gamma_n c) \cap \xi^3(\gamma_n d) = \xi^1(z_0).$$
Reasoning as in Basic Leaf Inclusions, Claim (3), we see $\mathcal{G}_\mathrm{p}(z_0, w_0)$ is a leaf of $\mathcal{G}$. □

Now take a foliation chart centered at a point in the leaf $\mathcal{G}_\mathrm{p}(z_0, w_0)$ produced in the previous step. Let the notations $\mathbb{D}, \phi_F, \phi_B$ and so on from §5.1.1 be for this foliation chart, with $\phi_B(0) = \xi^2(w_0) \cap \xi^3(z_0)$ and $\phi_F(0) = \xi^1(z_0)$.

**Step 2** (Backward Endpoints Rigid). *If $\mathcal{G} \notin \{\mathcal{G}_\mathrm{pcf}, \mathcal{G}_\mathrm{tcf}\}$ then $\phi_B(\mathbb{D}) \subset \xi^2(w_0)$.*

*Proof.* If $\phi_B(\mathbb{D}) \not\subset \xi^2(w_0)$ then continuity, Basic Leaf Inclusions (3) and the Dense Leaves Lemma imply that $\mathcal{G} = \mathcal{G}_\mathrm{pcf}$. □

Step (2) gives a substantial constraint. The next step describes its consequences and a useful book-keeping mechanism for making use of the constraint.

**Step 3** (Foliation Splitting). *Construction of the auxiliary local foliation $\Lambda_{z_0 w_0}$.*

The projective planes containing $\xi^2(w_0)$ are a projective line in $\mathbb{P}((\mathbb{R}^n)^*)$, and any two distinct such projective planes intersect in exactly $\xi^2(w_0)$. In particular, this implies that in our foliation chart, $\mathcal{G}$ restricts to local foliations of the intersections of $\Omega^1_\rho$ with nearby planes of the form $\xi^2(w_0) \oplus \xi^1(x)$ for $x$ near $z_0$ in $\partial\Gamma$. So $\mathcal{G}$ has a local product structure that is adapted to the geometry available to us. In this step we describe the coordinate system on $\mathbb{D}$ adapted to this structure and a useful auxilliary local foliation.

Let $\mathcal{A}$ be an affine chart containing the closure of $\mathcal{G}_\mathrm{p}(z_0, w_0)$. For a neighborhood $V_1$ in $\xi^2(w_0)$ of the point $\xi^2(w_0) \cap \xi^3(z_0)$ and a neighborhood $V_2$ in $\Xi^1$ of $\xi^1(z_0)$, every point in the convex hull in $\mathcal{A}$ of $V_1$ and $V_2$ is contained in a unique line in $\mathcal{A}$ between points in $V_1$ and $V_2$. This is because $\Xi^1$ is $C^1$ with tangent line $\xi^2(z_0)$ at $\xi^1(z_0)$ and the consequence of the hyperconvexity of $\xi$ that $\xi^2(w_0) \oplus \xi^2(z_0) = \mathbb{R}^4$. Call this auxilliary local foliation $\Lambda_{z_0 w_0}$. See Figure 25, Left.

The leaf space of $\Lambda_{z_0 w_0}$ is a product of intervals $I_B \times I_F$, which we identify with a subset of $\partial\Gamma^2$ by identifying the line segment between $\xi^2(w_0) \cap \xi^3(a)$ and $\xi^1(b)$ with $(a, b)$. This product is arranged so that leaves of $\Lambda_{z_0 w_0}$ with constant $I_B$-entry have the same endpoint on $\xi^2(w_0)$, and leaves of $\Lambda_{z_0 w_0}$ have the same endpoint on $\Xi^1$.

Since $\mathcal{G}_\mathrm{p}(z_0, w_0)$ is transverse to $\mathbb{D}$, after further restriction and rearranging by homeomorphism, we may take $\mathbb{D} = I_B \times I_F$ in the above notation. See Figure 25.

The point of using $\Lambda_{z_0 w_0}$ is that it is adapted to the constraint that backward endpoints lie in $\xi^2(w_0)$. In particular, $K_F$ takes a constrained form in the induced coordiantes on $\mathbb{D}$.

The first observation here is that from the qualitative description of $(\xi^2(w_0) \oplus \xi^1(z_0)) \cap \partial\Omega_\rho$ in Lemma 4.6 the line $L_{z_0 w_0}$ containing $\mathcal{G}_\mathrm{p}(z_0, w_0)$ intersects $\partial\Omega_\rho$ only at $\xi^2(w_0) \cap \xi^3(z_0)$ and $\xi^1(z_0)$. It follows that for a neighborhood $W$ of lines through $\xi^2(w_0)$ near $\ell_{w_0 z_0}$ that there is a neighborhood $U$ of $\xi^1(z_0)$ in $\partial\Omega_\rho$ so that every $L \in W$ meets $\partial\Omega_\rho$ only on $\xi^2(w_0) \cup U$. So after restriction, $\phi_F(p)$ is in $U$ for all $p \in I_B \times I_F$.

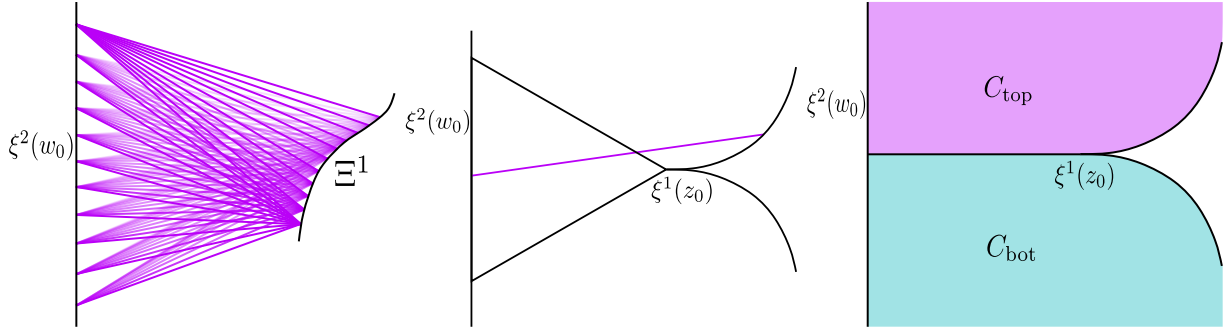

FIGURE 25. Left: the local foliation $\Lambda_{z_0 w_0}$ is made up of all possible nearby cusp-fans. This foliation is useful because it well-adapted to the local product structure of $\mathcal{G}$. Middle: in the argument of Step 4, the triangle traps leaves with endpoints on $\xi^2(w_0)$, which substantially constrains the shape of $K_F$ in $\mathbb{D}$. Right: Step 7 makes the key observation that the structure of the cusp of $\partial \Omega_\rho \cap (\xi^2(w_0) \oplus \xi^1(z_0))$ at $\xi^1(z_0)$ makes certain unions of $\mathcal{G}_p(z_0, w_0)$ with segments of $\xi^2(w_0)$ and $\partial \Omega_\rho \cap (\xi^2(w_0) \oplus \xi^1(z_0))$ bound convex domains.

**Step 4** (Triangle Trap). *For any $x \in I_F$, let $K_x$ be $K_F \cap (I_B \times \{x\})$. Then $K_x$ is a (relatively) closed interval in $I_B \times \{x\}$.*

*Proof.* Note $K_x$ is relatively closed and if $(a, x)$ and $(b, x)$ are in $K_x$, then from the structure of $(\xi^2(w_0) \oplus \xi^1(x)) \cap \partial \Omega_\rho$ these line segments and the segment of $\xi^2(w_0)$ between their $\xi^2(w_0)$-endpoints form a triangle $T$ in the closure of $\Omega_\rho^1$.

By our construction of $I_F \times I_B$ and the local product structure of $\mathcal{G}$, every point $q$ in $(I_F \times I_B) \cap (\xi^2(w_0) \oplus \xi^1(x))$ contained in the triangle $T$ has a leaf $\ell_q$ of $\mathcal{G}$ through $q$ that is entirely contained in $\xi^2(w_0) \oplus \xi^1(x)$, and $\ell_q$ must have its backward endpoint on $\xi^2(w_0)$. The only way this may occur is if $\ell_q$ has forward endpoint $\xi^1(x)$: otherwise, a crossing of leaves would occur between $\ell_q$ and an edge of $T$. See Figure 25, Middle.

We conclude that every point $c \in (a, b)$ has a leaf in $\mathcal{G}$ with endpoints $\xi^2(w_0) \cap \xi^3(c)$ and $\xi^1(x)$. This shows $K_x$ is an interval. $\square$

To simplify statements of later claims, suppose for contradiction that $\mathcal{G}$ is neither $\mathcal{G}_{\text{pcf}}$ nor $\mathcal{G}_{\text{tcf}}$. We will over-constrain $\mathcal{G}$. For $x \in I_F$, let $I_x = I_B \times \{x\}$ and $N = \{x \mid I_x \cap K_F \neq \emptyset\}$ be the entries $x \in I_F$ with some forward endpoint of the leaves of $\mathcal{G}$ through $I_B \times \{x\}$ on $\Xi^1$.

**Step 5** (Cusp Fans Sparse). *$N$ is closed in $I_F$ and has no interior.*

*Proof.* Closedness is from closedness of $K_F$. If $N$ had interior, an argument as in the One $\Xi^1$ Endpoint Step 2 of Case 1 shows $\mathcal{G}$ is either $\mathcal{G}_{\text{pcf}}$ or $\mathcal{G}_{\text{tcf}}$. $\square$

**Step 6** (Forward Endpoints Rigid off Fans). *For every connected component $\mathcal{C}$ of $(I_B \times I_F) - K_F$, there is a point $x_\mathcal{C} \in \partial \Gamma$ so that $\phi_F(\mathcal{C}) \subset \xi^2(x_\mathcal{C})$.*

*Proof.* By construction $\phi_F(\mathcal{C})$ is contained in $\Xi^2$. Now, otherwise an application of continuity, Basic Leaf Inclusions (3) and the Dense Leaves Lemma would show $\mathcal{G} = \mathcal{G}_{\text{pcf}}$. $\square$

**Step 7** (Exceptional Continuity at $z_0$). *The restriction $\phi_F|_{I_B \times \{z_0\}}$ is continuous.*

*Proof.* The structure of the intersection of $\xi^2(w_0) \oplus \xi^1(z_0)$ with $\partial \Omega_\rho$ implies that $\mathcal{G}_p(z_0, w_0)$ divides the nearby portions of $\xi^2(w_0)$, and $\partial \Omega_\rho$ in $\xi^2(w_0) \oplus \xi^1(z_0)$ into two boundary segments of properly convex sets. See Figure 25, Right.

That $\mathcal{G}$ restricts to a local foliation of $\xi^2(w_0) \oplus \xi^1(z_0)$ forces each component of $\Phi_F([(I_B \times \{z_0\}) - \{(z_0, z_0)\}])$ to be entirely contained in one convex set's boundary. Arguing as in the Discontinuity Control Lemma for $V_B$ and $V_F$ gives continuity. □

Let $\widetilde{N} \subset I_F$ be $\{x \in I_F \mid I_B \times \{x\} \subset K_F\}$, i.e. the elements of $N$ containing the largest possible fans in our cross-section. Just as above, $\widetilde{N}$ is closed in $I_F$ and has no interior.

**Step 8** (Full Fan). $z_0 \in \widetilde{N}$.

*Proof.* Corollary of steps 6 and 7. □

**Step 9** (Fan Complement Structure). *For every connected component $\mathcal{C}$ of $I_F - \widetilde{N}$, there is a connected component $\widetilde{\mathcal{C}}$ of $(I_B \times I_F) - K_F$ so that for all $x \in \mathcal{C}$, the segment $[(I_B \times \{x\}) - K_F] \cap \widetilde{\mathcal{C}}$ is nonempty.*

*Conversely, if $\widetilde{\mathcal{C}}$ is a connected component of $(I_B \times I_F) - K_F$, then the projection $\mathcal{C}$ of $\widetilde{\mathcal{C}}$ to $I_F$ is a connected component of $I_F - \widetilde{N}$.*

*Proof.* As in the case of once-punctured projective lines, recall from classical plane topology that if a closed subset $C$ of $\mathbb{C}$ separates $z_1, z_2 \in \mathbb{C}$, then a connected component of $C$ separates $z_1, z_2$ [21, Thm. 14.2].

In our setting, all connected components of $K_F$ are line segments, and only full line segments, i.e. lines over points in $\widetilde{N}$, separate points in our coordinate box $I_B \times I_F$. Both claims are consequences of this. □

**Step 10** (Better Fan). *There is some $z_0' \in I_F$ so that:*
  *(1) $z_0' \in \widetilde{N}$,*
  *(2) $z_0' \in I_B$,*
  *(3) There is an open interval $J \subset I_F - \widetilde{N}$ so that $z_0' \in \partial J$.*

*Proof.* Since $\widetilde{N}$ is a closed subset of an interval with empty interior, endpoints of intervals in the complement of $\widetilde{N}$ are dense in $\widetilde{N}$. The claim now follows from that $z_0 \in \widetilde{N}$ by the previous step and that $z_0 \in I_B$ by construction. □

We now re-center our notation at $z_0'$. The first two requirements on $z_0'$ in Step 10 ensure that $z_0'$ has all of the features of $z_0$ used in the proofs of statements (notably Step 7), so the same statements hold with $z_0'$ in place of $z_0$.

**Step 11** (Backward Stability). *There is a neighborhood $W_2 \subset \partial \Gamma$ of $z_0'$ so that for all $x \in W_2$ there is a point $p \in I_B \times \{x\}$ so that $\phi_B(p) = \xi^2(w_0) \cap \xi^3(z_0')$.*

*Proof.* This follows from the continuity of $\phi_B : I_B \times I_F \to \xi^2(w_0)$ and the consequence of the construction of $z_0'$ that $\phi_B(I_B \times \{z_0'\})$ contains an interval containing $\xi^2(w_0) \cap \xi^3(z_0')$ in its interior. □

**Step 12** (Conclusion). $\mathcal{G}$ *is* $\mathcal{G}_{\text{pcf}}$ *or* $\mathcal{G}_{\text{tcf}}$.

We are now finally ready to conclude. Let $J \subset I_F$ be an interval in the complement of $\widetilde{N}$ with $z_0' \in \partial J$. Step 6 shows there is a $c \in \partial \Gamma$ so that $\phi_F((x, y)) \subset \xi^2(c)$ for all $(x, y) \in I_B \times J$ not in $K_F$. By Step 11 we may take a sequence $(x_n, y_n) \in I_B \times I_F$ so that $y_n$ are in the complement of $N$, so $\lim_{n \to \infty} y_n = z_0'$, and so $\phi_B(x_n, y_n) = \xi^2(w_0) \cap \xi^3(z_0')$ for all $n$.

We see then that $\mathcal{G}$ contains the line segments $s_n$ from the single point $\xi^2(w_0) \cap \xi^3(z_0')$ to $\xi^2(c) \cap (\xi^2(w_0) \oplus \xi^1(y_n))$ for the above $c \in \partial \Gamma$ that is independent of $n$. Observe that the structure of $\xi^3(z_0') \cap \Omega_\rho$ implies $c \neq z_0'$: the boundary of the convex domain in $\xi^3(z_0')$ blocks the relevant segments from meeting $\xi^2(z_0')$. Indeed, the unique line segment between $\xi^2(w_0) \cap \xi^3(z_0')$ and a point in $\xi^2(z_0') - \{\xi^1(z_0')\}$ that is contained entirely in one connected component of $\Omega_\rho \cap \xi^3(z_0')$ is inside of $\Omega_\rho^2$.

Now, from the qualitative structure of line segments from $\xi^2(w_0) \cap \xi^3(z_0')$ to $\partial \Omega_\rho$ near $\mathcal{G}(z_0, w_0)$ (Lemma 4.15), we see that the line segments $s_n$ converge to an embedded segment $s_\infty$ from $\xi^2(w_0) \cap \xi^3(z_0')$ to $\xi^2(c) \cap (\xi^2(w_0) \oplus \xi^1(z_0'))$.

Then $s_\infty$ is a leaf of $\mathcal{G}$. But we have already forced $\mathcal{G}$ to contain the hull of a segment in $\xi^2(w_0)$ containing $\xi^2(w_0) \cap \xi^3(z_0')$ and $\xi^1(z_0)$ (Step 10). We have now forced $\mathcal{G}$ to include distinct leaves that intersect, which gives the desired contradiction. This concludes the case, and hence the proof. □

## 5.2. Codimension-2 Trifoliations.

In this subsection we prove the complementary Exactly Two Theorem for $\Omega_\rho^2$:

**Theorem 5.7** (Exactly Two: $\Omega_\rho^2$). *Let $\rho : \Gamma \to \mathrm{PSL}_4(\mathbb{R})$ be Hitchin. Let $\mathcal{G}$ be a $\rho(\Gamma)$-invariant trifoliation of $\Omega_\rho^2$ so that every leaf of $\mathcal{G}$ is a properly embedded properly convex subset of a projective line or a once-punctured projective line. Then $\mathcal{G} = \mathcal{G}_{\mathrm{ctaf}}$ or $\mathcal{G} = \mathcal{G}_{\mathrm{ctrf}}$.*

The approach is similar in basic outline to that for $\Omega_\rho^1$, but encounters distinct technical matters. The main new complication is that $\Omega_\rho^2$ lies on the non-convex side of $\Xi^2$, and leaves of $\mathcal{G}$ that are tangent to $\Xi^2$ at their endpoints can now occur. This makes the preparation for the main argument require more book-keeping.

On the other hand, the asymmetry of the structure of the sides of cusps in $\partial \Omega_\rho^2$ ends up making the final rigidity argument considerably less complicated, by cutting off the routes of temporary escape that a potential counter-example had available to it in §5.1.

5.2.1. *Discontinuities of Endpoints.* Let $\mathcal{G}$ be a trifoliation of $\Omega_\rho^2$ satisfying the hypotheses of Thm. 5.7. Let $B^\alpha$ ($\alpha \in \{1,2,3\}$) be a trifoliation chart for $\mathcal{G}$ about some $p \in \Omega_\rho^2$. As before, let $\mathbb{D}_\alpha$ be a cross-section and write forward and backwards endpoint projections as $\phi_F^\alpha$ and $\phi_B^\alpha$, respectively.

For $x \in \mathbb{D}_\alpha$, denote the projective line containing the line segment in the local $\alpha$-th trifoliation chart through $x$ by $\ell_x^\alpha$.

**Definition 5.8.** *Denote:*
(1) $O = \{x \in \mathbb{D} \mid \phi_F^\alpha(x) = \phi_B^\alpha(x)\}$,
(2) $K_{B,\overline{\pitchfork}} = \{x \in \mathbb{D}_\alpha \mid \ell_x^\alpha \text{ is tangent to } \Xi^2 \text{ at } \phi_B^\alpha(x)\} - O$,
(3) $K_{F,\overline{\pitchfork}} = \{x \in \mathbb{D}_\alpha \mid \ell_x^\alpha \text{ is tangent to } \Xi^2 \text{ at } \phi_F^\alpha(x)\} - O$,
(4) $K_{B,\pitchfork} = (\phi_B^\alpha)^{-1}(\Xi^1)$, and $K_{F,\pitchfork} = (\phi_F^\alpha)^{-1}(\Xi^1)$,
(5) $K_{\overline{\pitchfork}} = K_{F,\overline{\pitchfork}} \cap K_{B,\overline{\pitchfork}}$ and $K_{F,B,\overline{\pitchfork}} = K_{B,\overline{\pitchfork}} \cup K_{F,\overline{\pitchfork}}$,
(6) $V_F = \{x \in \mathbb{D}_\alpha \mid \ell_x^\alpha \text{ is transverse to } \Xi^2 \text{ at } \phi_F^\alpha(x)\}$,
(7) $V_B = \{x \in \mathbb{D}_\alpha \mid \ell_x^\alpha \text{ is transverse to } \Xi^2 \text{ at } \phi_B^\alpha(x)\}$,
(8) $V_{F,B} = V_F \cup V_B$ and $V = V_F \cap V_B$,
(9) $\mathscr{L} = \mathbb{D}_\alpha - V$, which is equal to $K_{F,B,\overline{\pitchfork}} \cup K_{F,\pitchfork} \cup K_{B,\pitchfork} \cup O$.

A first useful observation is:

**Lemma 5.9.** *If $x \in \mathbb{D}_\alpha$ has a leaf $g_x^\alpha$ of $\mathcal{G}$ through $x$ with an endpoint in $\Xi^1$, then $g_x^\alpha$ is a leaf of $\mathcal{G}_{\mathrm{ctrf}}$. Furthermore, $K_{B,\pitchfork} \cap O = K_{F,\pitchfork} \cap O = \emptyset$. If $x \in K_{B,\overline{\pitchfork}}$ has $\phi_B^\alpha(x) = \xi^2(x) \cap \xi^3(y)$, then $g_x^\alpha \subset \xi^3(y)$. If $x \in K_{\overline{\pitchfork}}$ then $g_x^\alpha$ is a leaf of $\mathcal{G}_{\mathrm{ctaf}}$.*

*Proof.* The first claim is a consequence of noting from our classification of intersections of projective planes with $\partial\Omega_\rho$ that points $\xi^1(y) \in \Xi^1$ always appear as cusps of $P \cap \Omega_\rho^2$ with tangent line $\xi^3(y) \cap P$, and with $\Omega_\rho^2$ on the nonconvex side of each segment of the cusp. This also forces $K_{B,\pitchfork} \cap O = K_{F,\pitchfork} \cap O = \emptyset$.

The second point follows from Lemma 4.2, in particular that all points $p \in \Xi^2$ are $C^1$ points of $\partial\Omega_\rho$ with tangent plane to $\xi^2(x) \cap \xi^3(y)$ given by $\xi^3(y)$. The final claim follows from examination of the qualitative form of $\xi^3(x) \cap \partial\Omega_\rho$. $\square$

The following is the analogue of Lemma 5.3 in this setting.

**Lemma 5.10** (Discontinuity Control). *Let $\gamma \in \Gamma - \{e\}$. In the above notation,*

(1) *$K_{F,\pitchfork}, K_{B,\pitchfork}$, and $K_{F,B,\overline{\pitchfork}}$ are disjoint and $K_{\overline{\pitchfork}}$, $K_{F,\pitchfork}$, $K_{B,\pitchfork}$, and $\mathscr{L}$ are closed.*
(2) *$K_{F,\overline{\pitchfork}} - K_{\overline{\pitchfork}}$ and $K_{B,\overline{\pitchfork}} - K_{\overline{\pitchfork}}$ are open in $\mathscr{L}$. Both $V_B$ and $V_F$ are open,*
(3) *The restrictions of both $\phi_F^\alpha$ and $\phi_B^\alpha$ to $K_{\overline{\pitchfork}}$, $K_{F,\pitchfork}$, and $K_{B,\pitchfork}$ are continuous. The restrictions of $\phi_B^\alpha$ to $K_{F,\overline{\pitchfork}} - K_{\overline{\pitchfork}}$ and $V_B$ and of $\phi_F^\alpha$ to $K_{B,\overline{\pitchfork}} - K_{\overline{\pitchfork}}$ and $V_F$ are continuous.*
(4) *Every $x \in O$ has $\ell_x^\alpha$ tangent to $\partial\Omega_\rho$ at $\phi_F^\alpha(x) = \phi_B^\alpha(x)$. The closure $\overline{O}$ of $O$ satisfies $\overline{O} \subset O \cup K_{\overline{\pitchfork}}$, and $O$ is open in $O \cup K_{\overline{\pitchfork}}$. The restrictions of $\phi_B^\alpha$ and $\phi_F^\alpha$ to $O$ are continuous.*

*Proof.* The closedness of $K_{F,\pitchfork}$, $K_{B,\pitchfork}$ and $K_{\overline{\pitchfork}}$ follow from Lemma 5.9, which demonstrates the leaves through $K_{F,\pitchfork}$ and $K_{B,\pitchfork}$ are leaves of $\mathcal{G}_{\text{ctrf}}$ and the leaves through $K_{\overline{\pitchfork}}$ are leaves of $\mathcal{G}_{\text{ctaf}}$. The disjointnesss assertion in (1) also follows from Lemma 5.9.

As in the last subsection, the main observation on openness is that if the leaf $g_x^\alpha$ through $x \in \mathbb{D}_\alpha$ is transverse to $\phi_F^\alpha(x) \in \Xi^2$ (resp. $\phi_B^\alpha(x)$), then there is an open set $U \subset \mathbb{D}_\alpha$ containing $x$ so that the restriction of $\phi_F^\alpha$ (resp. $\phi_B^\alpha$) to $U$ is continuous and has image contained in $\Xi^2$, with every leaf through $U$ transverse to this endpoint. It follows that $V_F, V_B$ and $V_{F,B}$ are open and that $\phi_F^\alpha|_{V_F}$ and $\phi_B^\alpha|_{V_B}$ are continuous. Hence $\mathscr{L}$ is closed and it follows similarly that $K_{F,\overline{\pitchfork}} - K_{\overline{\pitchfork}}$ and $K_{B,\overline{\pitchfork}} - K_{\overline{\pitchfork}}$ are open in $\mathscr{L}$.

Continuity of $\phi_F^\alpha|_{K_{\overline{\pitchfork}}}$, $\phi_F^\alpha|_{K_{\overline{\pitchfork}}}$, $\phi_F^\alpha|_{K_{F,\pitchfork}}$, $\phi_B^\alpha|_{K_{F,\pitchfork}}$, $\phi_F^\alpha|_{K_{B,\pitchfork}}$, and $\phi_B^\alpha|_{K_{B,\pitchfork}}$ follow from continuity of endpoints of leaves in $\mathcal{G}_{\text{ctrf}}$ and $\mathcal{G}_{\text{ctaf}}$.

We next address continuity of $\phi_B^\alpha|_{K_{F,\overline{\pitchfork}} - K_{\overline{\pitchfork}}}$. The analogue with $\phi_F^\alpha$ is symmetric. The fact that we use for continuity is that if $l_n$ and $m_n$ are sequences of projective lines in $\mathbb{RP}^3$ so $l_n \cap m_n$ is exactly one point $p_n \in \mathbb{RP}^3$ for each $n$, if both $l_n$ and $m_n$ have limits $\lim_{n\to\infty} l_n = l$, $\lim_{n\to\infty} m_n = m$ and if the limits intersect in exactly one point $l \cap m = p$, then $\lim_{n\to\infty} p_n = p$.

Now, take a sequence of points $x_n$ converging to $x$ in $K_{F,\overline{\pitchfork}} - K_{\overline{\pitchfork}}$ with associated leaves $g_x^\alpha \subset \ell_x^\alpha$. Write $\phi_B^\alpha(x_n) = \xi^2(a_n) \cap \xi^3(b_n)$ and $\phi_B^\alpha(x) = \xi^2(a) \cap \xi^3(b)$. As in the discussion for $\Omega_\rho^1$, the lines $\ell_{x_n}^\alpha$ converge to $\ell_x^\alpha$. Since $\phi_B^\alpha|_{V_B}$ is continuous, $a_n$ converges to $a$ and $b_n$ converges to $b$. It follows that $\phi_F^\alpha(x_n) = \ell_{x_n} \cap \xi^2(b_n)$ and $\lim_{n\to\infty} \ell_{x_n} \cap \xi^2(b_n) = \ell_x \cap \xi^2(b)$. Examination of local models at $\phi_B^\alpha$-endpoints, as in §5.1, shows that $g_x^\alpha$ is the segment from $\xi^2(a) \cap \xi^3(b)$ to $\ell_x \cap \xi^2(b)$, and hence $\phi_F^\alpha(x_n)$ converges to $\phi_F^\alpha(x)$.

For (4), begin by considering the collection $\mathscr{C}$ of projective lines $\ell \in \text{Gr}_2(\mathbb{R}^4)$ so that $\ell \cap \partial\Omega_\rho$ is contained in $\Xi^2$ and $\ell$ is tangent to $\Xi^2$ at every point of intersection. Then $\mathscr{C}$ is relatively closed in the collection of lines in $\text{Gr}_2(\mathbb{R}^4)$ that do not intersect $\Xi^1$. Note also the property of not intersecting $\Omega_\rho^1$ is a closed condition among projective lines. As our classification of intersections of projective planes with $\partial\Omega_\rho$ shows that every projective line that intersects $\Xi^1$ also intersects $\Omega_\rho^1$, the closure of $O$ is contained in $O \cap K_{\overline{\pitchfork}}$. The openness assertion follows from closedness of $K_{\overline{\pitchfork}}$. Continuity of $\phi_F^\alpha$ and $\phi_B^\alpha$ on $O$ follows from noting that given a line $\ell$ that is tangent to $\Xi^2$ at $p$ and only meets $\partial\Omega_\rho$ at $p$, for every neighborhood $V$ of $p$ there is a neighborhood $U \in \mathrm{Gr}_2(\mathbb{R}^4)$ of $\ell$ so that every point of intersection of any line $\ell' \in U$ with $\partial\Omega_\rho$ is in $V$. □

5.2.2. *Leaf Inclusions.* We now address basic leaf inclusions. To disambiguate leaves of $\mathcal{G}_{\mathrm{ctrf}}$ contained in the same line, define

$$\mathcal{G}_{\mathrm{ta}}(a_0, b_0) = \{\xi^3(a_0) \cap \xi^3(b_0) \cap \xi^3(c) \mid (a_0, b_0, c) \in \partial\Gamma^{(3)+}\}.$$

We also adopt the notation that the line segment between $\xi^1(x)$ and $\xi^2(y) \cap \xi^3(x)$ in $\Omega_\rho^2$ is $\mathcal{G}_{\mathrm{tr}}(x, y)$. In the following we give precise statements and indicate where our arguments are essentially identical to the case of $\Omega_\rho^1$.

**Lemma 5.11** (Basic Leaf Inclusions). *Let $\mathcal{G}$ be a $\rho(\Gamma)$-invariant trifoliation of $\Omega_\rho^2$ whose leaves are properly embedded properly convex subsets of projective lines or once-punctured projective lines. Let $\gamma \in \Gamma - \{e\}$.*

*(1) If $\ell = \mathcal{G}_{\mathrm{tr}}(\gamma^-, a)$ is a leaf of $\mathcal{G}$, then so is $\mathcal{G}_{\mathrm{tr}}(\gamma^-, \gamma^+)$. If $\mathcal{G}_{\mathrm{tr}}(a, \gamma^-)$ is a leaf of $\mathcal{G}$, then so is $\mathcal{G}_{\mathrm{tr}}(\gamma^+, \gamma^-)$.*

*(2) If $\ell = \mathcal{G}_{\mathrm{ta}}(\gamma^-, a)$ is a leaf of $\mathcal{G}$, then $\mathcal{G}_{\mathrm{ta}}(\gamma^-, \gamma^+)$ is also a leaf of $\mathcal{G}$. If $\ell = \mathcal{G}_{\mathrm{ta}}(a, \gamma^-)$ is a leaf of $\mathcal{G}$, then $\mathcal{G}_{\mathrm{ta}}(\gamma^+, \gamma^-)$ is a leaf of $\mathcal{G}$.*

*(3) If $\ell$ is a leaf of $\mathcal{G}$ with endpoints $\xi^2(a) \cap \xi^3(\gamma^-)$ and $\xi^2(\gamma^-) \cap \xi^3(c)$ with $(\gamma^-, a, c) \in \partial\Gamma^{(3)+}$, then $\mathcal{G}_{\mathrm{ta}}(\gamma^-, \gamma^+)$ is also a leaf of $\mathcal{G}$. If, instead, $(\gamma^-, c, a)$ is positively oriented, then $\mathcal{G}_{\mathrm{ta}}(\gamma^+, \gamma^-)$ is a leaf of $\mathcal{G}$.*

*(4) If $\ell$ is a leaf of $\mathcal{G}$ with endpoint $\xi^2(\gamma^-) \cap \xi^3(a)$ and other endpoint $p \notin \xi^3(\gamma^-)$, then $\mathcal{G}_{\mathrm{tr}}(\gamma^-, \gamma^+)$ is also a leaf of $\mathcal{G}$.*

*(5) If $\ell$ is a leaf of $\mathcal{G}$ with exactly one endpoint, which has the form $\xi^2(\gamma^-) \cap \xi^3(a)$, then $\mathcal{G}_{\mathrm{ta}}(\gamma^-, \gamma^+)$ and $\mathcal{G}_{\mathrm{ta}}(\gamma^+, \gamma^-)$ are leaves of $\mathcal{G}$.*

*Proof.* The strategy in each case, as it was for $\Omega_\rho^1$, is to apply $\rho(\gamma)^n$ to each leaf and take a limit, making use of continuity of the lines containing leaves through a cross-section of $\mathcal{G}$.

Claims (1) and (2) follow by arguments directly analogous to those for $\Omega_\rho^1$. For claim (3), let $\widetilde{\ell}$ be the line containing $\ell$. We first note that

$$\lim_{n\to\infty} \rho(\gamma)^n \widetilde{\ell} = \lim_{n\to\infty} \rho(\gamma)^n [((\xi^2(a) \cap \xi^3(\gamma^-)) \oplus (\xi^2(\gamma^-) \cap \xi^3(c))]$$
$$= (\xi^2(\gamma^+) \cap \xi^3(\gamma^-)) \oplus (\xi^2(\gamma^-) \cap \xi^3(\gamma^+))$$
$$= \xi^3(\gamma^+) \cap \xi^3(\gamma^-).$$

The disambiguation of leaves in the limiting line is an analysis of intersections with other planes in the image of $\xi^3$, similar to the analogous analysis for $\Omega_\rho^1$, and shows that the limiting segment is $\mathcal{G}_{\mathrm{ta}}(\gamma^-, \gamma^+)$. Continuity of leaf-lines in a cross-section then gives the desired conclusion.

For (4), by Lemma 5.9 the endpoint $p$ must be in $\Xi^2$ and have the form $p = \xi^2(c) \cap \xi^3(d)$ with $d \neq \gamma^-$. We must also have $c \neq \gamma^-$ for the leaf not to be contained in $\partial\Omega_\rho$. With $\widetilde{\ell}$ the line containing $\ell$,

$$\lim_{n\to\infty} \rho(\gamma^n) \widetilde{\ell} = \lim_{n\to\infty} \rho(\gamma)^n [(\xi^2(\gamma^-) \cap \xi^3(a)) \oplus (\xi^2(c) \cap \xi^3(d))]$$
$$= (\xi^2(\gamma^-) \cap \xi^3(\gamma^+)) \oplus \xi^1(\gamma^+),$$

using transversality of the direct summands of the limit. The analysis of verifying the desired segment is contained in the limit exactly follows the proof of (3) in the Basic Leaf Inclusions Lemma for $\Omega^1_\rho$.

For the final claim, note that $\ell$ must be tangent to $\Xi^2$ at $p = \xi^2(\gamma^-) \cap \xi^3(a)$, and hence contained in $\xi^3(\gamma^-)$ by Lemma 4.2. As $\ell$ is not $\xi^2(\gamma^-)$, there must be a point $p' \in \ell - \xi^2(\gamma^-)$. We then note that $\lim_{n \to \infty} \rho(\gamma^n)p' = \xi^2(\gamma^+) \cap \xi^3(\gamma^-)$. On the other hand, $\lim_{n \to \infty} \rho(\gamma^n)p = \xi^2(\gamma^-) \cap \xi^3(\gamma^+)$. Examining the limiting projective lines, as for once-punctured line segments in the $\Omega^1_\rho$ case, shows $\mathcal{G}_{\mathrm{ta}}(\gamma^-, \gamma^+)$ and $\mathcal{G}_{\mathrm{ta}}(\gamma^+, \gamma^-)$ are leaves of $\mathcal{G}$. □

The Dense Leaves Lemma holds essentially verbatim in this setting, with the same proof:

**Lemma 5.12** (Dense Leaves). *Let $\mathcal{G}$ be a $\rho(\Gamma)$-invariant trifoliation of $\Omega^2_\rho$ whose leaves are properly embedded projective line segments or once-punctured projective lines. Suppose that there is a neighborhood $U \subset \partial \Gamma^{(2)}$ so that for a dense set $S$ of $U$, for every $(a, b) \in S$ either $\mathcal{G}_{\mathrm{ta}}(a, b)$ or $\mathcal{G}_{\mathrm{tr}}(a, b)$ is a leaf of $\mathcal{G}$. Then $\mathcal{G} \in \{\mathcal{G}_{\mathrm{ctaf}}, \mathcal{G}_{\mathrm{ctrf}}\}$.*

5.2.3. *The Exactly 2 Theorem for $\Omega^2_\rho$.* We now enter the main proof for $\Omega^2_\rho$.

*Proof of the Exactly Two Theorem for $\Omega^2_\rho$.* Let $\mathcal{G}$ be a $\rho(\Gamma)$-invariant trifoliation of $\Omega^2_\rho$ satisfying the hypotheses of Thm. 5.7. We split the argument into the two cases that $\mathscr{L} = \mathbb{D}_\alpha$ and when there is a leaf with both endpoints transverse to $\Xi^2$.

Take a foliation chart $B_\alpha \cong \mathbb{D}_\alpha \times (-1, 1)$ ($\alpha \in \{1, 2, 3\}$) around an interior point $p$. We do case analysis on the endpoints of leaves through $\mathbb{D}_\alpha$, as before.

**Case 1** (No Double-Transverse Leaves). *Every leaf of $\mathcal{G}$ either is tangent to $\Xi^2$ at one endpoint or has an endpoint in $\Xi^1$.*

This is equivalent to $\mathbb{D}_\alpha = \mathscr{L}$ for every trifoliation chart.

**Step 1** (Always a $\Xi^1$ Endpoint). *If either $K_{F,⋔}$ or $K_{B,⋔}$ is not totally disconnected in $\mathbb{D}_\alpha$, then $\mathcal{G} \in \{\mathcal{G}_{\mathrm{ctrf}}, \mathcal{G}_{\mathrm{ctaf}}\}$.*

*Proof.* Assume $K_{F,⋔}$ is not totally disconnected, the other case is of course symmetric. Then let $\mathcal{C}$ be a connected component of $K_{F,⋔}$ that is not a point. Then by Lemma 5.10, $\phi^\alpha_F$ and $\phi^\alpha_B$ are continuous, so write $\phi^\alpha_F(x) = \xi^1(a(x))$ and $\phi^\alpha_B(x) = \xi^2(b(x)) \cap \xi^3(a(x))$ for $x \in \mathcal{C}$. As $x \mapsto (\phi^\alpha_F(x), \phi^\alpha_B(x))$ is an injection on $\mathbb{D}_\alpha$ the restriction to $\mathcal{C}$ of $x \mapsto (a(x), b(x))$ is a continuous injection of a connected subset of $\mathbb{D}_\alpha$ that is not a point, and hence either $a(x)$ or $b(x)$ surjects an open interval $U$ of $\partial \Gamma$.

In either case, the Basic Leaf Inclusions Lemma (1) and the Dense Leaves Lemma show that $\mathcal{G} = \mathcal{G}_{\mathrm{ctaf}}$. This finishes the step. □

**Step 2** (All Endpoints Tangent). *If $K_{⌢}$ is not totally disconnected in $\mathbb{D}_\alpha$, then $\mathcal{G}$ is $\mathcal{G}_{\mathrm{ctaf}}$ or $\mathcal{G}_{\mathrm{ctrf}}$.*

*Proof.* We argue as in the previous case. On $K_{⌢}$, write $\phi^\alpha_F(x) = \xi^2(a(x)) \cap \xi^3(b(x))$ and $\phi^\alpha_B(x) = \xi^2(b(x)) \cap \xi^3(a(x))$. Then, as discussed above, the map $x \mapsto (a(x), b(x))$ is continuous and locally injective on $K_{⌢}$. If $K_{⌢}$ is not totally disconnected, then on any non-point connected component $\mathcal{C}$ of $K_{⌢}$, either $a(x)$ or $b(x)$ surjects an interval. The Basic Leaf Inclusions Lemma Claim (2) together with the Dense Leaves Lemma then implies that $\mathcal{G} = \mathcal{G}_{\mathrm{ctaf}}$. □

**Step 3** (One Tangent Endpoint). *Under the hypotheses of the case, if $V_F \cap K_{B,⌢}$ or $V_B \cap K_{F,⌢}$ is nonempty, then $\mathcal{G} \in \{\mathcal{G}_{\mathrm{ctaf}}, \mathcal{G}_{\mathrm{ctrf}}\}$.*

*Proof.* Then, up to a symmetric argument with $K_{B,\overline{\sqcap}}$, we may take $K_{F,\overline{\sqcap}} - K_{\overline{\sqcap}}$ to be a nonempty open subset of $\mathbb{D}_\alpha$. After restricting we may take $\mathbb{D}_\alpha = K_{F,\overline{\sqcap}} - K_{\overline{\sqcap}}$. Then the Discontinuity Control Lemma shows that $\phi_F^\alpha$ and $\phi_B^\alpha$ are continuous. For $x \in \mathbb{D}_\alpha$, write $\phi_F^\alpha(x) = (\xi^2(a(x)) \cap \xi^3(b(x))$ and $\phi_B^\alpha(x) = \xi^2(b(x)) \cap \xi^3(c(x))$.

Then the map $\Phi : \mathbb{D}_\alpha \to \partial\Gamma^3$ given by $x \mapsto (a(x), b(x), c(x))$ is continuous. The image of $\Phi$ composed with a projection to a $\partial\Gamma$ factor of $\partial\Gamma^3$ must contain open sets in at least two projections: otherwise, we would obtain a continuous local injection from an open subset of $\mathbb{R}^2$ into $\mathbb{R}$.

If the projection to the $a(x)$ factor contains an open set $U$, then claim (4) of Lemma 5.11 shows $\mathcal{G}$ contains $\mathcal{G}_{\text{tr}}(\gamma^-, \gamma^+)$ as a leaf for all pole-pairs $(\gamma^-, \gamma^+)$ with $\gamma^- \in U$. Then the Dense Leaves Lemma forces $\mathcal{G} = \mathcal{G}_{\text{ctrf}}$. If the projection to the $b(x)$ factor contains an open set $U$, an analogous argument, using Lemma 5.11.(3) allows us to conclude $\mathcal{G} = \mathcal{G}_{\text{ctaf}}$.

As at least two projections contain a nonempty open set, we are finished with this case. □

The previous steps have now reduced the case to proving the following.

**Step 4** (Once-Punctured $\mathbb{RP}^1$). *With the hypotheses of the case it is impossible for all of $K_{F,\text{ft}}, K_{B,\text{ft}}$, and $K_{\overline{\sqcap}}$ to have empty interior and for both $V_F \cap K_{B,\overline{\sqcap}}$ and $V_B \cap K_{F,\overline{\sqcap}}$ to be empty in every cross-section $\mathbb{D}_\alpha$ of a foliation chart for $\mathcal{G}$.*

*Proof.* The method in this case is rather similar to the case for $\Omega^1_\rho$, with some minor modifications to the setting. Suppose for contradiction that the situation we wish to show is impossible occurs. By our hypotheses in this case and the resolutions of the above cases, every endpoint of a leaf of $\mathcal{G}$ on $\Xi^2$ is tangent to $\Xi^2$.

So our cross-section $\mathbb{D}_\alpha$ for any trifoliation chart has $\mathbb{D}_\alpha = O \cup K$. On $O$, write $\phi_F^\alpha(x) = \xi^2(a(x)) \cap \xi^3(b(x))$, with $a(x)$ and $b(x)$ continuous. The analogue to the motivating claim in the case of $\Omega^1_\rho$ here is:

**Claim.** *$a(x)$ is locally constant on $O$.*

The argument is identical to the corresponding claim in the case of $\Omega^1_\rho$. The complementary claim is also analogous:

**Claim.** *The intersection of $O$ with any cross-section $\mathbb{D}_\alpha$ is connected.*

*Proof.* The complement of $O$ in $\mathbb{D}_\alpha$ is $K_{\overline{\sqcap}} \cup K_{F,\text{ft}} \cup K_{B,\text{ft}}$, which is the union of three closed, disjoint, and totally disconnected plane sets, whose union is hence connected by [21, Thm. 14.2] and an elementary covering space argument. □

Combining the two claims, it follows that every leaf of $\mathcal{G}$ that is a once-punctured projective line has its unique endpoint contained in a common line $\xi^2(x_0)$ for some $x_0 \in \partial\Gamma$, which is incompatible with $\rho(\Gamma)$-invariance of $\mathcal{G}$. □

This exhausts the possibilities when all leaves of $\mathcal{G}$ have at least one endpoint either tangent to $\Xi^2$ or in $\Xi^1$. This completes the step and the case.

**Case 2** (Generic Concave Case). $\mathbb{D}_\alpha \neq K_{F,\text{ft}} \cup K_{B,\text{ft}} \cup K_{F,B,\text{ft}}$.

In this case, $V_F \cap V_B \neq \emptyset$, and so we may take $\mathbb{D}_\alpha = V_F \cap V_B$. Write the line containing the leaf of some $x \in \mathbb{D}_\alpha$ by $\ell_x$, and the leaf by $g_x^\alpha$. Then there are $a, b, c, d \in \partial\Gamma$ so that $\ell_x = (\xi^2(a) \cap \xi^3(d)) \oplus (\xi^2(c) \cap \xi^3(d))$ with $a \neq b$ and $c \neq d$. Note also that $a \neq c$, lest $\ell_x \subset \partial\Omega_\rho$. Furthermore, $b \neq d$ by inspecting the qualitative shape of $\xi^3(b) \cap \partial\Omega_\rho$. Also, $b \neq c$ or $g_x^\alpha$ would have an endpoint tangent to $\xi^2(c)$, again using the qualitative shape of $\xi^3(c) \cap \partial\Omega_\rho$.

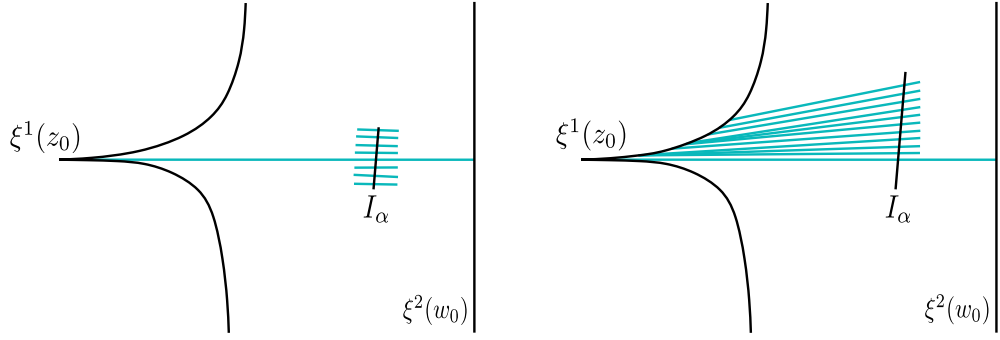

FIGURE 26. Left: the local picture near the limiting leaf. Right: as for $\Omega^1_\rho$, the presence of the limiting leaf prevents discontinuities from occurring on the cross-section $I_\alpha$ within $\xi^2(w_0) \cap \xi^1(z_0)$. The argument in this case is from foliation implying monotonicity of endpoints on the convex arc near $\xi^1(z_0)$.

To obtain a better leaf to work with again, since $\Gamma$ is a uniform convergence group, there is a sequence $\gamma_n \in \Gamma$ and $z_0 \neq w_0$ in $\partial\Gamma$ so that

$$\lim_{n\to\infty} \gamma_n a = w_0, \qquad \lim_{n\to\infty} \gamma_n y = z_0 \text{ for all } y \neq a \text{ in } \partial\Gamma.$$

Applying $\rho(\gamma_n)$ to $\ell_x$ gives us

$$\lim_{n\to\infty} \rho(\gamma_n)\ell_x = \lim_{n\to\infty} (\xi^2(\gamma_n a) \cap \xi^3(\gamma_n b)) \oplus (\xi^2(\gamma_n c) \oplus \xi^3(\gamma_n d)) = (\xi^2(w_0) \cap \xi^3(z_0)) \oplus \xi^1(z_0).$$

An analysis as in the generic case of the $\Omega^1_\rho$ argument shows that $\lim_{n\to\infty} \gamma_n g_x^\alpha = \mathcal{G}_{\mathrm{tr}}(z_0, w_0)$ and that $\mathcal{G}_{\mathrm{tr}}(z_0, w_0)$ is a leaf of $\mathcal{G}$.

Let $p \in \mathcal{G}_{\mathrm{tr}}(z_0, w_0)$ and take a foliation chart for $\mathcal{G}$ containing a segment in $\mathcal{G}_{\mathrm{tr}}(z_0, w_0)$, and adopt the notations we have repeatedly used for this foliation chart. Arrange so that $\phi_B^\alpha(0) = \xi^2(w_0) \cap \xi^3(z_0)$ and $\phi_F^\alpha(0) = \xi^1(z_0)$. After restriction, we may arrange for $\mathbb{D}_\alpha = V_B$ and for $\phi_B^\alpha$ to be continuous.

If the image $\phi_B^\alpha(\mathbb{D}_\alpha)$ is not contained in $\xi^2(w_0)$ then Lemma 5.11.(4) and the Dense Leaves Lemma complete the argument as in the analysis for $\Omega^1_\rho$. So we may restrict to the case where every $x \in \partial\mathbb{D}_\alpha$ has $\phi_B^\alpha(x) \in \xi^2(w_0)$.

We make use of this major local structure constraint by observing, as in the case of $\Omega^1_\rho$, that it forces leaves of $\mathcal{G}$ near $g_0^\alpha$ to locally foliate $\xi^2(w_0) \oplus \xi^1(z_0)$ in a neighborhood of $p$, with leaves transverse to $I_\alpha := \mathbb{D}_\alpha \cap (\xi^2(w_0) \oplus \xi^1(z_0))$. See Figure 26, Left.

Let $I_\alpha^+$ be one (closed) side of $I_\alpha$ about 0. We claim:

**Claim.** *The image $\phi_F^\alpha(I_\alpha^+)$ contains an interval of $\partial\Omega_\rho^2 \cap (\xi^2(w_0) \oplus \xi^1(z_0))$.*

*Proof of Claim.* After restricting $I_\alpha^+$ in an affine chart $\mathcal{A} \cong \mathbb{R}^2$ containing $g_0^\alpha$ we arrange for:
  (1) The point $\xi^1(z_0)$ to be at $(0,0)$,
  (2) The line $(\xi^2(w_0) \cap \xi^3(z_0)) \oplus \xi^1(z_0)$ to be the horizontal axis,
  (3) $I_\alpha^+$ to be contained in the first quadrant, and all points in $I_\alpha^+$ to have $x$-value in $[1, 2]$.
  (4) The line $\xi^2(x) \cap \mathcal{A}$ to be a vertical line $x = c > 2$,
  (5) Every leaf through $I_\alpha^+$ to have slope in $[-1, 1]$.

To get the last point, we are using continuity of the slopes of leaves on $\mathbb{D}_\alpha$.

Parametrize $I_\alpha^+$ as $[0,1]$ with $0$ at $p$. We note that because leaves cannot cross and $\Omega_\rho^2$ is on the nonconvex side of the cusp in $\xi^2(w_0) \oplus \xi^1(z_0)$, the function $\phi_F^\alpha(x)$ is monotone in $x \in I_\alpha^+$ and $0$ is the only element of $I_\alpha^+$ for which $\phi_F^\alpha(x) = (0,0)$.

Consequently, every $x \neq 0$ in $I_\alpha^+$ is a continuity point of $\phi_F^\alpha$ by the Discontinuity Control Lemma. The claim follows if we show $\lim_{x \to 0} \phi_F^\alpha(x) = (0,0)$. To see this, suppose otherwise. Then there is a sequence $x_n$ converging to $0 \in I_\alpha^+$ so that $\lim_{n \to \infty} \phi_F^\alpha(x_n) = (a,b)$ with $a, b > 0$. This forces the lines through $x_n$ to have limiting slope $b/a > 0$, in contradiction to the continuity of lines through a local cross-section of $\mathcal{G}$. This is the desired contradiction. □

Now note that we have obtained leaves with endpoints in $\Xi^2$ in continuously and non-trivially varying lines of the ruling of $\partial\Omega_\rho$. This is enough to apply the Basic Leaf Inclusions Lemma 5.11.(4) and the Dense Leaves Lemma to conclude the case and with it the theorem. □

## 6. Equivalences of Properly Convex Foliated Projective Structures

We resolve the remaining matters claimed in the introduction and not yet recorded.

First, we claimed that the codimension-2 conclusions in Thm. A fail if either the hypothesis of proper containment of leaves in projective lines or of group invariance is removed. Let us give sketches of counter-examples that occur when these hypotheses are omitted.

A counter-example to Thm. A in the absence of group invariance is given by fixing $x_0 \in \partial\Gamma$ and noting the line segments in $C_x$ between $\xi^2(x_0) \cap \partial C_x$ and all other points in $\partial C_x$ ($x \in \partial\Gamma$) foliate $\Omega_\rho^1$.

For the necessity of disallowing entire projective lines as leaves, it is known that for any $\mathrm{PSp}_4(\mathbb{R})$-Hitchin representation $\rho$ the domain $\Omega_\rho^1$ admits a $\rho(\Gamma)$-invariant foliation by projective lines, as a consequence of [4]. The case of Fuchsian representations may be done explicitly in the polynomial model of [11, §3]. This is enough for the claim. It is also interesting to note that in contrast to $\rho(\Gamma)$-invariant foliations whose leaves are properly embedded segments, such foliations are *not* rigid. It is an exercise to find deformations of such a foliation that are supported in a tube over a small neighborhood in a cross-section. Such deformations can then be propagated with $\rho(\Gamma)$-invariance to obtain continuous families of $\rho(\Gamma)$-invariant foliations of $\Omega_\rho^1$ whose leaves are all projective lines.

Finally, let us deduce Cor. D on the rigidity of equivalences of properly convex foliated projective structures.

A relevant fact proved by Guichard-Wienhard is that the holonomy $\mathrm{hol}_A : \pi_1(\mathrm{T}^1 S) \to \mathrm{PSL}_4(\mathbb{R})$ of any properly convex foliated projective structure $(A, \mathcal{F}, \mathcal{G})$ on $\mathrm{T}^1 S$ vanishes on the kernel of the canonical projection $\pi_1(\mathrm{T}^1 S) \to \pi_1(S)$ [11, Cor. 5.29]. We write the induced map $\Gamma \to \mathrm{PSL}_4(\mathbb{R})$ as $\mathrm{hol}_{*A}$. The map given by induced holonomies $\mathrm{hol}_{*A}$ from $\mathcal{P}_{\mathrm{pcf}}$ to the $\mathrm{PSL}_4(\mathbb{R})$-character variety is the map Guichard-Wienhard show is a homeomorphism onto $\mathrm{Hit}_4(S)$.

*Proof of Thm. D.* Let $(A_1, \mathcal{F}_1, \mathcal{G}_1)$ and $(A_2, \mathcal{F}_2, \mathcal{G}_2)$ be properly convex foliated projective structures and $\varphi \in \mathrm{Homeo}_0(S)$ be an isomorphism of the projective structures $A_1$ and $A_2$.

Now $\varphi$ lifts to $\mathrm{T}^1 \widetilde{S}$ because $\varphi \in \mathrm{Homeo}_0(S)$, and so gives a projective map $\Omega^1_{\mathrm{hol}_{*A_1}} \to \Omega^1_{\mathrm{hol}_{*A_2}}$, which is the restriction of some $B \in \mathrm{PSL}_4(\mathbb{R})$.

From our classification theorems, we must have $\mathrm{dev}_{A_1}(\mathcal{F}) = \mathcal{F}_{\mathrm{pcf}}$, $\mathrm{dev}_{A_1}(\mathcal{G}) = \mathcal{G}_{\mathrm{pcf}}$, and $\mathrm{dev}_{A_2}(\mathcal{F}) = \mathcal{F}_{\mathrm{pcf}}$, $\mathrm{dev}_{A_2}(\mathcal{G}) = \mathcal{G}_{\mathrm{pcf}}$ for the respective domains. Thm. E forces preservation of $\mathcal{F}_{\text{pcf}}$. The Exactly Two Theorem for $\Omega^1_\rho$ and the nesting of leaves of $\mathcal{F}_{\text{pcf}}$ and $\mathcal{G}_{\text{pcf}}$ implies preservation of $\mathcal{G}_{\text{pcf}}$, so that $\varphi$ is an isomorphism of properly convex foliated $(\mathbb{RP}^3, \text{PSL}_4(\mathbb{R}))$-structures. □

## Appendix A. Proof of General Position Projection Lemma

We record a proof of the General Position Projection Lemma 4.8 here. Let us recall the statement:

**Lemma A.1** (General Position Projection)**.** *Let $V$ be an $n$-dimensional vector space. Let $X, Y, W_1, ..., W_k$ be subspaces of $V$ so that $V = X \oplus Y = X \oplus W_1 \oplus W_2 \oplus \cdots \oplus W_k$.*

*Then for $j = 1, ..., k$,*

$$(A.1) \qquad Y \cap (X \oplus W_1 \oplus \cdots \oplus W_j) = \bigoplus_{i=1}^{j} ((X \oplus W_i) \cap Y),$$

*and the sum on the right is direct. In particular, $Y = \bigoplus_{i=1}^{k}((X \oplus W_i) \cap Y)$, and the sums in Eq. A.1 have dimension $\dim(Y \cap (X \oplus W_1 \oplus \cdots \oplus W_j)) = \sum_{i=1}^{j} \dim W_j$ for $j = 1, ..., k$.*

*Proof.* We first establish Eq. (A.1), without its statement of directness on the right side, by induction on $j$. For $j = 1$ the claim is vacuous. Now suppose the claim holds for $j$ and examine $X, W_1, ..., W_{j+1}$. Then by the induction hypothesis,

$$(A.2) \qquad \sum_{i=1}^{j+1} ((X \oplus W_i) \cap Y) = (Y \cap (X \oplus W_1 \oplus \cdots \oplus W_j)) + ((X \oplus W_{j+1}) \cap Y).$$

We next observe that the terms on the right of Eq. (A.2) have dimension $\sum_{i=1}^{j} \dim(W_i)$ and $\dim(W_{j+1})$, respectively. Indeed, the dimension of $X \oplus W_1 \oplus \cdots \oplus W_j$ is $(n - \dim Y) + \sum_{i=1}^{j} \dim W_i$ by the directness of the sums $X \oplus Y$ and $X \oplus W_1 \oplus \cdots \oplus W_k$. This forces the subspace $Z_j$ defined by $Y \cap (X \oplus W_1 \oplus \cdots \oplus W_j)$ to have dimension at least $\sum_{i=1}^{j} \dim W_i$. To see that equality holds for $\dim Z_j$, note that if $\dim Z_j > \sum_{i=1}^{j} \dim W_i$ then there would be a nonzero vector $v \in X \cap Z_j$, which would force the impossible conclusion $v \in X \cap Y$. The other dimension computation that $\dim[(X \oplus W_{j+1}) \cap Y] = \dim W_{j+1}$ is similar.

Each term on the right of Eq. (A.2) is contained in $Y \cap (X \oplus W_1 \oplus \cdots \oplus W_{j+1})$, which has dimension $\sum_{i=1}^{j+1} \dim(W_j)$. To show that the desired equality holds, it hence suffices to show that the sum on the right of Eq. (A.2) is direct, i.e. any

$$v \in (Y \cap (X \oplus W_1 \oplus \cdots \oplus W_j)) \cap (Y \cap (X \oplus W_{j+1}))$$

must be 0. For any such $v$, we have $v \in Y$ and $v \in (X \oplus W_1 \oplus \cdots \oplus W_j) \cap (X \oplus W_{j+1}) = X$, by hypothesis. So $v \in X \cap Y$ and hence $v = 0$, which completes the induction. The directness of the sum in Eq. (A.1) is now a dimension count using the above. □